%% file: SGFEMmanu.tex
\documentclass{article}
\usepackage{amsmath,amssymb}
\usepackage[dvips]{graphicx}
\usepackage{theorem}
\usepackage{esint}
\usepackage{color}
\usepackage{booktabs}
\usepackage{multirow}
\usepackage{lineno}
\newtheorem{assump}{Assumption}
\newtheorem{mytheorem}{Theorem}[section]
\theoremstyle{plain}\newtheorem{myprop}[mytheorem]{Proposition}
\theoremstyle{plain} \newtheorem{mylemma}[mytheorem]{Lemma}
\theoremstyle{plain}
\theoremstyle{plain}\newtheorem{myremark}[mytheorem]{Remark}
\numberwithin{equation}{section}

\newcommand{\ih}{\mathcal{I}_h}
\newcommand{\intp}[1]{\mathcal{I}_{#1}}

\newcommand{\myend}{{\huge$\centerdot$\vspace{0.1ex}}}
\newcommand{\mcal}[1]{\mathcal{#1}}
\newcommand{\locphi}[2]{\varphi_{#1}^{[#2]}}
\newcommand{\locphibar}[2]{\overline{\varphi}_{#1}^{[#2]}}
\newcommand{\myqed}{ \raisebox{.8ex}{\framebox[2mm]{}}}
\newcommand{\open}[1]{\stackrel{o}{#1}}

\newcommand{\captionfonts}{\small}

\makeatletter
\long\def\@makecaption#1#2{%
    \vskip\abovecaptionskip
    \sbox\@tempboxa{{\captionfonts #1: #2}}%
    \ifdim \wd\@tempboxa >\hsize
        {\captionfonts #1: #2\par}
    \else
        \hbox to\hsize{\hfil\box\@tempboxa\hfil}%
    \fi
    \vskip\belowcaptionskip}
\makeatother

\begin{document}

\title{Stable Generalized Finite Element Method\\(SGFEM)}

\author{I. Babu\v{s}ka \thanks{
ICES, University of Texas at Austin,
Austin, TX.}
\and U. Banerjee
\thanks{
Department of Mathematics, 215 Carnegie, Syracuse University,
Syracuse, NY 13244. E-mail address: banerjee@syr.edu. This research was partially
supported by IMA, University of Minnesota, Minneapolis, MN and J. T. Oden Faculty Fellowship, ICES, University of Texas at Austin, Austin, TX.}
}

\date{}

\maketitle

\begin{abstract}
The Generalized Finite Element Method (GFEM) is a Partition of Unity Method (PUM), where the  trial space of standard Finite Element Method (FEM) is augmented with non-polynomial shape functions with compact support. These shape functions, which are also known as the enrichments, mimic the local behavior of the unknown solution of the underlying variational problem. GFEM has been successfully used to solve a variety of problems with complicated features and microstructure. However, the stiffness matrix of GFEM is badly conditioned (much worse compared to the standard FEM) and there could be a severe loss of accuracy in the computed solution of the associated linear system. In this paper, we address this issue and propose a modification of the GFEM, referred to as the Stable GFEM (SGFEM). We show that the conditioning of the stiffness matrix of SGFEM is not worse than that of the standard FEM. Moreover, SGFEM is very robust with respect to the parameters of the enrichments. We show these features of SGFEM on several examples.
\end{abstract}

{\bf Keywords:} Generalized finite element method (GFEM); partition of unity (PU);  Extended Finite Element Method (XFEM); approximation; condition number, loss of accuracy, linear system; Validation and Verification

\input sec1_intro.tex

\input sec2_modelprob.tex

\input sec3_GFEM.tex
\input newSec4_SGFEMa.tex
\input newSec4_SGFEMb.tex

\input sec5p1_appl.tex
\input sec5p2_appl.tex

\input sec5p3_appl.tex
\input sec_conclu.tex
\input sec6_appendix.tex


\bibliography{gfem}
\bibliographystyle{plain}

\end{document}

%% file: sec1_intro.tex
\section{Introduction}\label{intro}

During the last decade, the Generalized Finite Element Method (GFEM) and the eXtended Finite Element Method (XFEM) -- two approaches based on the Partition of Unity Method (PUM) -- were developed independently and have been widely used to solve various types of problems. Only recently, it was clearly recognized that these two methods are same and were referred to as XFEM/GFEM (\cite{FriesBelRev}). Hence we believe that it is interesting to briefly describe the early development of these methods. It was also recognized that, though these methods have  excellent convergence properties, the stiffness matrices associated with these methods could be ill-conditioned. In this paper, we especially address this issue and propose an easy modification, which we call the Stable Generalized Finite Element Method (SGFEM), that does not have the above mentioned conditioning problem and is very robust.

We start with a brief history of the early development of the methods based on PUM. Since this is history and not a survey, it is important to provide not only the publication date, but in addition, the submission dates for various papers, and pay careful attention to nomenclature for the various methods.

\smallskip

\noindent \textbf{\emph{Brief early history:}} The idea of adding non-polynomial basis functions into the trial space of the FEM started in 1970's (\cite{Benzley,Byskov,FixGulWak}). However, these basis functions had global support and the associated stiffness and mass matrices lost their local structure.

Three Special FEMs, which used non-polynomial shape functions, were proposed in  \cite{BabCalOsb}(1994, sub: Mar.1992) to solve second order problems with rough coefficients. In particular, the shape functions used in the Special FEM \#3 have compact supports and are products of piecewise linear FE hat-functions and a non-polynomial function that mimic the special features of the unknown solution. This idea was further generalized with detailed mathematical theory and applications in the Ph.D. dissertation of J. M. Melenk \cite{MelPhd}(1995), where it was shown that the hat-functions could be replaced by any PU (with compact support). This method was referred to as PUM and PUFEM in \cite{MelBab}(1996, sub: Apr.1996) and \cite{BabMel}(1997, sub: Jul.1995), respectively; these papers contain major results on the method and its application to the problems with highly oscillatory solutions, problems with solutions with boundary layers, differential equations with rough coefficients, etc.

The PUM was referred to as the GFEM in \cite{StrBabCop}(2000, sub: Jul.1998), \cite{StrCopBab1}(2000, sub: Nov.1998), \cite{StrCopBab}(2001, sub: Jul.2000), where the hat-functions were used as the PU (similar to the Special FEM \#3 in \cite{BabCalOsb}). In these papers, GFEM is interpreted as an FEM augmented with non-polynomial shape functions with compact support, and it is shown that the use of only a few of these shape functions is enough to address the problems with singular solutions. Moreover, the idea of obtaining the non-polynomial shape functions by solving certain local problems is also introduced in these papers in the context of the analysis of a perforated plate.

In a parallel development, but independent of \cite{StrBabCop}, \cite{StrCopBab1}, \cite{StrCopBab}, PUM with hat-functions serving as the PU was also investigated in \cite{BelyBla}(1999, sub: Jul.1998) and \cite{MoeDolBel},(1999, sub: Feb.1999) in the context of crack propagation problems. This method is similar to GFEM  as it also uses the standard FE trial space augmented with non-polynomial shape functions. However, the major contributions of these papers is to show that the method does not need remeshing  as the crack propagates. Also shape functions with jump discontinuities were used in these papers. The method was first referred to as the XFEM in the Ph.D. thesis of J. Dolbow \cite{DolbowPhd}(1999) and almost simultaneously in \cite{SukMoeEtAl}(2000, sub: Sep.1999), \cite{DauxMoesEtAl}(2000, sub: Sep.1999], and \cite{DolMoeBel1}(2000, sub: Sep.1999).
We mention that XFEM was employed in \cite{SukMoeEtAl, DauxMoesEtAl, DolMoeBel1} to address crack propagation problems in 3-d, problems with branched cracks, and fracture in Reissner-Mindlin plates.

Another idea similar to the PUM was used in the $h$-$p$ Cloud method in \cite{DuarOden}(1996, sub: Jun.1995), \cite{DuarOden2}(1996, sub: Apr.1996), where the shape functions were the products of a PU and polynomials. The goal of this method is to obtain $h$-$p$ FEM like approximation without using a FE mesh, in the spirit of meshless methods. The use of the ``customized function'' (which mimicked the exact solution) for crack problems was also suggested in \cite{OdenDuar}(1997, sub: Dec.1996), under this framework. Later, the hat-functions were also used as the PU in the $h$-$p$ Cloud method in \cite{OdenDuarZien}(1998, sub: Dec. 1996).

Lot of work has been done in the area of GFEM and XFEM since these early work, described above. We will comment on some of the recent developments near the end of this section.

\smallskip

\noindent \textbf{\emph{GFEM and the problem with conditioning:}} PUM is a flexible framework to design Galerkin methods that accurately approximate solutions of variational problems. The framework involves (a) accurately approximating the solution, locally, using functions in a local approximation space, and (b) gluing the local approximations, using a PU, to construct a globally conforming approximate solution. The GFEM, which is a PUM with FE hat functions serving as the PU, retains the important flexibility of choosing the local approximation space. The efficiency of GFEM lies in the fact that is requires only modifying an existing FE code to incorporate special shape functions with compact support. The GFEM, with appropriate choice of special shape functions, leads to excellent convergence properties. However, the use of hat-functions as PU may result into almost linearly dependent shape functions in GFEM, and the stiffness matrix could be severely ill-conditioned; the ill-conditioning could be much worse than the conditioning of the stiffness matrix of the FEM. This results into the loss of accuracy in the solution of the linear system associated with the GFEM. In fact, the shape functions could be linearly dependent yielding a singular stiffness matrix.

Various ad-hoc approaches have been developed in the literature to address this issue. For example, the stiffness matrix of GFEM was perturbed by an identity matrix of size $\epsilon$ (small) in \cite{StrBabCop, StrCopBab} and an iterative method was used to solve the perturbed linear system. Preconditioning of the stiffness matrix, based on domain decomposition, have been recently suggested in \cite{MenkBordas} to address the conditioning problem. In \cite{GrieSchw6, Schweit}, a flat-hat PU (modified FE hat functions with flattened top) was used in the PUM instead of hat-functions. The use of flat-hat PU certainly avoids the problem of loss of accuracy in the linear system, but it requires developing a code from the scratch.

Naturally, it is pertinent to ask if GFEM could be modified so that it retains the excellent convergence properties of the GFEM, and the loss of accuracy in the computed solution of the linear system of the modified GFEM is of the same order as that of the standard FEM. In this paper, we will show that the SGFEM has both of these features. We have chosen a 1-d problem to present the idea of the SGFEM primarily for the clarity of exposition and not to obscure the analysis with details that are not directly related to the SGFEM. However, the ideas and the associated analysis (including the notational machinery) could be easily generalized to higher dimensions and will be reported in a future publication.

\medskip

\noindent \textbf{\emph{Indicator of the loss of accuracy in computed solution of the linear system:}} Consider the linear system $Ax=b$, associated with FEM, GFEM, or SGFEM, where $A$ is an $n \times n$ sparse symmetric positive definite matrix. Let $\hat{x}$ be the computed solution of the linear system, obtained from an elimination method encoded in a linear algebra package and the computations follow the IEEE standard for floating point arithmetic (with guard digits). Set $\eta :=\|x-\hat{x}\|_2/\|x\|_2$ -- the relative error that measures the loss of accuracy in the computed solution.  $\eta$ depends on the round-off, but in general, it also depends on the elimination algorithm and its implementation in the package, the compiler, the processor, and the computing platform with single or multiple processors. $\eta$ is related to the relative error in the approximate solution due to round-off.

We seek an indicator that reliably indicates the loss of accuracy in the computed solution, characterized by $\eta$, and is practically independent of  other factors mentioned above. Let $H = DAD$, where $D$ is a diagonal matrix with $D_{ii}=A_{ii}^{-1/2}$. Define the \emph{scaled condition number} $\mathfrak{K}(A)$ of $A$ by
$
\mathfrak{K}(A):= \kappa_2(H),
$
where $\kappa_2(H)=\|H\|_2\|H^{-1}\|_2$ is the condition number of $H$ based on the $\|\cdot \|_2$ vector norm. We hypothesize that $\mathfrak{K}(A)$ is the indicator, which we formalize as follows:
\begin{quote}
\noindent \textbf{Hypothesis H:}
\begin{equation}
\eta \approx Cn^\beta \mathfrak{K}(A) \epsilon; \ \beta \approx 0, \label{HypH}
\end{equation}
where $\epsilon$ is the machine precision. \end{quote}


\noindent  We will elaborate on the precise meaning of the hypothesis and \emph{validate} it in the Appendix, borrowing the ideas from the area of \emph{Validation and Verification}. The indicator $\mathfrak{K}(A)$ will be used to compare various GFEMs with respect to the loss of accuracy in the computed solution, which will allow us to choose a preferable GFEM. In particular, we will show in this paper that
$\mathfrak{K}(A^{SGFEM}) \le \mathfrak{K}(A^{GFEM})$, where $A^{SGFEM}$ and $A^{GFEM}$ are the stiffness matrices of SGFEM and GFEM, respectively, and therefore the SGFEM is preferable over the GFEM.



%

\medskip

\noindent \textbf{\emph{Some current work in GFEM/XFEM:}} These methods have been used in a variety of applications. For example, XFEM has been used recently to address two-phase fluid flow problems (\cite{EsserEtal}), mechanical behavior of nano-structures (\cite{FarsadEtal}), and heterogeneous material with random interfaces (\cite{NouyClem}); GFEM has been used to address heat transfer problems with sharp thermal gradient (\cite{OharaDuarteEason}), grain boundary in polycrystals (\cite{SimDuarEtal}), and electromagnetic problems (\cite{LuShan}). Special shape functions for problems with locally periodic coefficients are constructed in \cite{MatBabSch} that yield exponential order of convergence. Also local problems to compute the shape functions for problems with rough coefficients are constructed in \cite{BabLip}, and it has been proved that GFEM yields exponential order of convergence. For an extensive collection of references in XFEM/GFEM, we refer to \cite{FriesBelRev}.

\medskip

\noindent \textbf{\emph{Organization of the paper:}} In Section \ref{MP}, we give the model problem in 1-d. We intentionally chose the problem in 1-d so that we could communicate the main ideas of SGFEM, when applied to this problem, in a fairly general fashion, without the notational and other technical complexity associated with higher dimensions. We describe the PUM and GFEM, together with the approximation results, in Section \ref{GFEM} and show the conditioning problem in GFEM on an example.
In Section \ref{SGFEM}, we first describe the SGFEM in a simpler setting, show that SGFEM retains the convergence properties of GFEM, and establish that the scaled condition numbers of the stiffness matrices of the SGFEM and FEM are of the same order. We chose the simpler setting primarily to communicate the main idea of the method and the associated analysis. We then describe the SGFEM and provide the analysis in full generality. We note that some of the ideas presented here could have been presented in a simpler fashion by using 1-d arguments. However we did not take this approach; the notations and framework of the analysis, developed in this section, could be easily generalized to higher dimensions. In Section \ref{APPLIC}, we applied SGFEM to three specific examples, namely, interface problems, problems with singular solutions, and problems with discontinuous solutions. In the Appendix, we discuss the validation of Hypothesis H and present many validation experiments. We note that the Appendix is a very important part of this paper

%% file: sec2_modelprob.tex
\section{Model problem} \label{MP}


Let $\Omega = (0,1)$ and, for an integer $k \ge 0$, we denote the standard Sobolev spaces by $H^k(\Omega)$ with the norm $\|\cdot \|_{H^k(\Omega)}$ and seminorm $|\cdot |_{H^k(\Omega)}$; for $k=0$, $H^0(\Omega)=L^2(\Omega)$. We would also use the spaces $H^k(A)$, where $A$ is a sub-domain of $\Omega$. Consider the variational problem
\begin{equation}
u \in H^1(\Omega),\ \ B(u,v) = F(v), \quad \forall \ v \in H^1(\Omega), \label{VarProb}
\end{equation}
where
\begin{equation}
B(u,v):= \int_\Omega a u^\prime v^\prime \, dx  \quad \mbox{and }F(v):= \int_\Omega fv\, dx\label{Buv}
\end{equation}
such that $F(1)= \int_\Omega f\, dx = 0$. We assume that the function $a(x)$ is bounded, i.e., there are constants $\alpha,\, \beta$ such that
\begin{equation}
0 < \alpha \le a(x) \le \beta,\quad \forall \ x \in \Omega \label{BndOfax}
\end{equation}
We note that  $a(x)$ could be smooth, but it also could be rough. It is well known that the problem \eqref{VarProb} has a unique solution, up to an additive constant.

We define the \emph{Energy norm}, $\|v\|_{\mcal{E}(A)}$, of $v \in H^1(A)$, where $A$ is a sub-domain of $\Omega$, by
\[
\|v\|_{\mcal{E}(A)}^2:= B_A(v,v), \quad \mbox{where } B_A(w,z):= \int_A a w^\prime z^\prime \, dx.
\]

It is well known that the solution $u$ of \eqref{VarProb} is also the solution of a boundary value problem (BVP), posed in the strong form as
\begin{equation}
- [a(x)u^\prime]^\prime = f, \quad au^\prime (0) = au^\prime (1) = 0 \label{BVP}
\end{equation}
provided $au^\prime$ is differentiable.

%% file: sec3_GFEM.tex
\section{Generalized Finite Element Method (GFEM):} \label{GFEM}

Let $\mcal{S}$ be a finite dimensional subspace of $H^1(\Omega)$. The Ritz-Galerkin method to approximate the solution $u$ of \eqref{VarProb} is given by
\begin{equation}
u_h \in \mcal{S},\ \ B(u_h,v) = F(v), \quad \forall \ v \in \mcal{S}. \label{RitzMet}
\end{equation}
The solution $u_h$ is unique up to an additive constant. We can obtain a unique solution by imposing a \emph{natural constraint} on $u_h$, namely, $u_h(0)=0$.

A \emph{Partition of Unity method} (PUM) is a Ritz-Galerkin method, where $\mcal{S}$ is constructed employing a (a) \emph{Partition of Unity} (PU) and (b) \emph{Local approximating spaces}. A Generalized Finite Element method (GFEM) is a PUM with special PU. We first briefly described the PUM.

For a parameter $h>0$, Let $I^h:= \{i \in \mathbb Z: 0 \le i \le N\}$, where $N=N(h)$ is an integer. For $i \in I^h$, let $\omega_i^h:= (a_i^h,b_i^h) \subset \Omega$ such that (i) $\Omega = \cup_{i\in I^h} \omega_i^h$, and (ii) any $x \in \Omega$ belongs to at most $\kappa$ of the open intervals $\omega_i^h$; $\kappa$ is independent of $i,\, h$. The open interval $\omega_i^h$ is called a \emph{patch}. Subordinate to the cover $\{\omega_i^h\}_{i\in I^h}$, let $\{N_i^h\}_{i\in I^h}$ be a $C^0$ PU satisfying
\[
\sum_{i\in I^h}N_i^h(x)=1,\ \forall \ x\in \Omega,\ \
 \|N_i^h\|_{L^\infty(\Omega)} \le C,\ \  \mbox{diam}\{\omega_i^h\} \|(N_i^h)^\prime\|_{L^\infty(\Omega)} \le C,
\]
where $C>0$ is independent of $i$ (for details, see \cite{BabMel,MelBab,BBO7}).

On each patch $\omega_i^h$, $i \in I^h$, we consider an $(n_i+1)$-dimensional space $V_i^h$ -- the \emph{local approximating space}, namely
\begin{equation}
V_i^h = \mbox{span}\{\varphi_{j}^{[i],h}\}_{j=0}^{n_i}, \ \varphi_{j}^{[i],h} \in H^1(\omega_i) \mbox{ and } \varphi_{0}^{[i],h}=1, \label{DefLocSp}
\end{equation}
where $n_i$s are non-negative integers. The functions $\varphi_{j}^{[i],h}$, $j >0$,  are carefully chosen such the functions in $V_i^h$ mimic the the exact solution $u$, locally in $\omega_i^h$. We will further comment on this issue later. In the rest of the paper, we will write $I,\, \omega_i,\, N_i, \, V_i,\, \locphi{j}{i}$ in place of $I^h,\, \omega_i^h,\, N_i^h,\, V_i^h$, $\varphi_{j}^{[i],h}$, respectively, with an understanding that they depend on the parameter $h$. The PUM is precisely \eqref{RitzMet}, with the finite dimensional space $\mcal{S}$ is given by
\begin{equation}
\mcal{S} = \sum_{i\in I} N_i V_i = \mbox{span}\{N_i \, \varphi_j^{[i]},\, 0\le j \le n_i,\, i \in I \}
:= \mcal{S}_1 + \mcal{S}_2, \label{GFEMspace}
\end{equation}
where
\begin{equation}
\mcal{S}_1=\{ \zeta :\, \zeta = \sum_{i\in I} y_0^{[i]}\varphi_0^{[i]} N_i\},\quad
\mcal{S}_2 =\{ \zeta :\, \zeta = \sum_{i\in I} \sum_{j=1}^{n_i} y_j^{[i]}\varphi_j^{[i]}N_i\}, \label{GFEMspace2}
\end{equation}
and $y_0^{[i]},\, y_j^{[i]} \in \mathbb R$. The functions $\locphi{j}{i}$, $j \ge 1$, and the associated spaces $V_i$ are sometimes referred to as \emph{enrichments} and \emph{enrichment spaces} respectively in the literature. We will refer to $\mcal{S}_1$ as the \emph{basic part of $\mcal{S}$}; $\mcal{S}_2$ will be referred to as the \emph{enrichment part} of $\mcal{S}$. Moreover, we will refer to the Galerkin method with $\mcal{S}=\mcal{S}_1$ as the \emph{basic part of PUM}. Thus every PUM has a basic part based only on the PU.

We now present the main approximation result of PUM in the Energy norm (see \cite{BabMel,MelBab,BBO7}).

\begin{mytheorem} \label{GFEMApproxThm} Suppose $u \in H^1(\Omega)$. Suppose for $i\in I$, there exists $\xi^i \in V_i$ and $C_1 >0$, independent of $i$, such that
\[
\|u-\xi^i\|_{L^2(\omega_i)} \le C_1 \mbox{diam}(\omega_i)\,\|u-\xi^i\|_{\mcal{E}(\omega_i)} \ \mbox{ and }\ \|u-\xi^i\|_{\mcal{E}(\omega_i)} \le \epsilon_i.
\]
Then there exists $v \in \mcal{S}$ such that
\begin{equation}
\|u - v\|_{\mcal{E}(\Omega)} \le C \big[\textstyle \sum_{i\in I} \epsilon_i^2 \big]^{1/2}, \label{GFEMApprox}
\end{equation}
where the positive constant $C$ depends on $\kappa, C_1,\beta/\alpha$.
\end{mytheorem}

It is immediate from Theorem \ref{GFEMApproxThm} that the PUM solution $u_h \in \mcal{S}= \mcal{S}_1+\mcal{S}_2$ of \eqref{RitzMet} satisfies
\begin{equation}
\|u - u_h\|_{\mcal{E}(\Omega)} \le \inf_{v \in \mcal{S}}\|u - v\|_{\mcal{E}(\Omega)} \le C \big[ \textstyle \sum_{i\in I}  \epsilon_i^2 \big]^{1/2}, \label{ErrEstGFEM}
\end{equation}
where $u$ is the solution of \eqref{VarProb}. It is clear from above that the global accuracy of the PUM solution $u_h$ depends on how accurately the solution $u$ of \eqref{VarProb} can be approximated by the functions in $V_i$, locally on the patches $\omega_i$.

We mention that in higher dimensions, the patches $\omega_i$ are subdomains, which can have quite general shape. Theorem \ref{GFEMApproxThm}, as presented above, is also true is higher dimensions.

We now describe the GFEM. Recall that the choice of PU in PUM is arbitrary. The GFEM is a PUM, where
(a) the patches $\omega_i$ are ``FE stars'' relative to a finite element (FE) triangulation of $\Omega$, and (b) the piecewise linear FE hat-functions $N_i$, associated with the vertices of FE triangulation, serve as the PU.

Let $N=1/h$ and recalling $I=\{i:\, 0 \le i \le N\}$, let $\mcal{T}:=\{x_i=ih:\, i\in I$. Let $\{\tau_k\}_{k \in I\backslash \{0\}}$ be the \emph{uniform mesh} on $\Omega$, where $\tau_k:=[x_{k-1},x_k]$ are the \emph{elements}; $\open{\tau}_k:= (x_{k-1},x_k)$ is the interior of $\tau_k$. The points $x_i$ are called the \emph{vertices} of the mesh. The patches $\{\omega_i\}_{i\in I}$ are defined as $\omega_i :=(x_{i-1},x_{i+1})$, $i=1,2,\cdots,N-1$; also $\omega_0 := (x_0,x_1)$ and $\omega_N := (x_{N-1},x_N)$. For $i \in I$, let $N_i$ be the standard hat-functions associated with the vertex $x_i$; the support of $N_i$ is $\overline{\omega}_i$. Note that $\overline{\omega}_0=\tau_1$, $\overline{\omega}_N=\tau_N$ and $\overline{\omega}_i=\tau_i \cup \tau_{i+1}$ for $i=1,2,\dots,N-1$. $\overline{\omega}_i$ is the FE star associated with the vertex $x_i$. Clearly, $\{N_i\}_{i\in I}$ form a PU subordinate to the patches $\{\omega_i\}_{i\in I}$.  The associated GFEM is the Galerkin method \eqref{RitzMet} with $\mcal{S} = \mcal{S}_1 + \mcal{S}_2$ (see \eqref{GFEMspace}). Clearly $\mcal{S}_1$ is the standard FE space of piecewise linear functions, and consequently, the \emph{basic part of the GFEM} is the standard finite element method (FEM). Thus the trial space $\mcal{S}$ of the GFEM is precisely the standard FE trial space, augmented with the space $\mcal{S}_2$. Thus GFEM could be implemented by incorporating enrichments into an existing FE code. The name GFEM was first used in \cite{StrBabCop,StrCopBab1} to highlight exactly this point. The description of GFEM is exactly same in higher dimensions; it is based on the standard FE triangulation of $\Omega$.


\begin{myremark} \upshape
We note that we have considered a uniform mesh only for the simplicity of exposition; in fact, the ideas and theory in this paper could also be presented for \emph{locally quasi-uniform meshes}, i.e., when $C_1\le |\tau_{k+1}|/|\tau_k| \le C_2$ for $k=1,\cdots,N-1$, with $C_1, C_2 >0$  independent of $k$. \myend
\end{myremark}



The accuracy of the GFEM (also PUM) solution depends on the choice of $V_i$, as mentioned before (see Theorem \ref{GFEMApproxThm}). The functions $\varphi_j^{[i]} \in V_i$ (see \eqref{DefLocSp},\eqref{GFEMspace}) are carefully chosen based on the available information on the unknown solution $u$ of \eqref{VarProb} to mimic the unknown solution locally in $\omega_i$. Examples of $V_i$, suitable for specific applications are available in the literature (e.g., see \cite{BBO7}). We briefly mention some of the examples that we will consider in this paper:

\medskip

\noindent $\bullet$ If the unknown solution $u$ is smooth in $\omega_i$, then the $\varphi_j^{[i]}$s are usually chosen to be polynomials in $\mcal{P}^j(\omega_i)$ and the associated spaces $V_i$ are spaces of polynomials of degree $n_i$. We note that $n_i$ could could be different for different $i$, based on the available information on $u$.

\smallskip

\noindent $\bullet$ When $a(x)$ is a piecewise smooth and discontinuous function (interface problems), $\varphi_j^{[i]}$s are chosen such that $a[\varphi_j^{[i]}]^\prime$ is smooth on $\omega_i$. Clearly, $\varphi_j^{[i]}$ are continuous piecewise smooth functions with derivatives that are discontinuous at the discontinuities of $a(x)$.

\smallskip

\noindent $\bullet$ If the unknown solution $u$ is singular, then $\varphi_j^{[i]}$ should be chosen as singular functions, mimicking the singularity of $u$.

\noindent

\noindent $\bullet$ If $u$ is discontinuous at $x=c$ in the domain, then $\varphi_j^{[i]}$s are chosen to be discontinuous functions on those $\omega_i$s that contain $x=c$. We note however that problems with  discontinuous solutions cannot be cast as \eqref{VarProb}; we will address these problems in Section \ref{DiscontProb} of this paper.

\begin{myremark}
\upshape GFEM provides a flexible framework to obtain various Galerkin methods. Many classical methods could be cast in this framework. For example, with $n_i=0$ for $i \in I$ in \eqref{DefLocSp}, GFEM (with $\mcal{S}=\mcal{S}_1$) yields the classical FEM. Moreover, let $s(x)$ be a function (could be singular) defined on $\Omega$. Consider $n_i=1$  and $\locphi{1}{i} = s(x)|_{\omega_i}$ in the definition of $\mcal{S}$ in \eqref{GFEMspace2}. Then GFEM, with $y_1^{[i]} = b$ (a constant) for $i\in I$, yields the classical ``singular FEM'' (see \cite{StrangFixBook, FixGulWak, BlumDob, Byskov, Benzley,RaoRajuMur}), where the standard finite element trial space is augmented by the global function $s(x)$. Moreover, we note that \emph{$n_i$ in \eqref{DefLocSp} could be different for different values of $i$}. In fact, one may use $n_i>0$ only for a few patches $\omega_i$, as needed for accuracy, based on the available information; for other patches, $n_i=0$, i.e., $V_i = \mbox{span}\{1\}$. This idea was also discussed and implemented in the original GFEM papers \cite{StrBabCop, StrCopBab1}.
\end{myremark}

\subsection{Scaled condition number of the stiffness matrix of GFEM}

The stiffness matrix $\mathbf{A}$ of the GFEM is positive semi-definite.
Even when the GFEM solution is naturally constrained with $u_h(0)=0$, i.e., when $\mathbf{A}$ is positive definite, the condition number $\kappa_2(\mathbf{A})$ can be extremely large, specifically larger than the condition number of standard FE stiffness matrix, which is $O(h^{-2})$ for second order problems. However, according to Hypothesis H, the scaled condition number $\mathfrak{K}(\mathbf{A})= \kappa_2(H)$ is a reliable indicator of the loss of accuracy in the computed solution of $\mathbf{A}x=b$. Recall $H=D\mathbf{A}D$, where $D$ is a diagonal matrix with $D_{ii}=\mathbf{A}_{ii}^{-1/2}$. We now present an example where $\mathfrak{K}(\mathbf{A})$ is much larger than the scaled condition number of the standard FE stiffness matrix, which is again $O(h^{-2})$.

Suppose $a(x)=1$ in \eqref{Buv} and let $u=x^\alpha$ with $1/2 < \alpha < 3/2$, $\alpha \ne 1$. Note that $x^\alpha \in H^1(\Omega)$ but $x^\alpha \notin H^2(\Omega)$. We consider a GFEM with $n_i=1$ and $\locphi{1}{i}:= x^\alpha|_{\omega_i},\, i\in I$. From the definition of $\mcal{S}$, any $v \in \mcal{S}$ is of the form
\begin{equation}
v(x) = \sum_{i\in I\backslash \{0\}} a_i N_i(x) + \sum_{i\in I} b_i N_i(x) x^\alpha;\quad a_i,b_i \in \mathbb R . \label{unscaledv}
\end{equation}
We have set $a_0=0$ to impose the constraint $u_h(0)=0$ on the GFEM solution $u_h$. It can be easily shown that $u_h = u$, i.e. \emph{there is no approximation error}.

We let $\eta:=[a_1,\dots,a_N,b_0,\cdots,b_N]^T \in \mathbb R^{2N+1}$. Then $B(v,v) = \eta^T \mathbf{A} \eta$, where $\mathbf{A}$ is the $(2N+1)\times(2N+1)$ positive definite stiffness matrix of the GFEM. We note that $\mathbf{A}_{ii}= |N_i|_{H^1(\omega_i)}^2$ for $1 \le i \le N$ and $\mathbf{A}_{N+1+j,N+1+j}= |N_j x^\alpha|_{H^1(\omega_j)}^2$ for $0\le j \le N$. Therefore by considering $v \in \mcal{S}$ of the form
\begin{equation}
v(x) = \sum_{i\in I\backslash \{0\}} a_i \frac{N_i(x)}{|N_i|_{H^1(\Omega)}} + \sum_{i\in I} b_i \frac{N_i(x) x^\alpha}{|N_i x^\alpha|_{H^1(\Omega)}},\quad a_i,b_i \in \mathbb R , \label{scaledv}
\end{equation}
it is easy to see that $B(v,v) = \eta^T H \eta$, where $H$ is as mentioned before.

We consider a $v \in \mcal{S}$ of the form \eqref{scaledv} with $a_i=0$ for $1 \le i \le N-1$, $a_N = 1$, and $b_i=0$ for $i \in I$. Then
\[
B(v,v) = \int_\Omega {v^\prime}^2 dx = 1 \quad \mbox{ and } \|\eta\|^2 := \sum_{i\in I\backslash \{0\}} a_i^2 + \sum_{i\in I} b_i^2 = 1.
\]
Therefore,
\begin{equation}
\frac{B(v,v)}{\|\eta\|^2} = 1 \le \lambda_{M}, \label{ev.eq1}
\end{equation}
where $\lambda_M$ is the largest eigenvalue of $H$.

Let $g(x) \in H^2(\Omega)$ be a non-decreasing function with $g(x)=0$ for $0 \le x \le 1/4$ and $0 < C \le g(x_i) \le 1$ for $i \ge \lceil N/2 \rceil$. For $h$ small enough, let $1/8 < x_k \le 1/4$ be the vertex closest to $x=1/4$. Clearly $x^\alpha$ and $g x^\alpha$ are in $H^2(\widehat{\Omega})$, where $\widehat{\Omega} := (1/8,1)$. We now consider a $v \in \mcal{S}$ of the form \eqref{scaledv} with $a_i=-g(x_i)x_i^\alpha|N_i|_{H^1(\Omega)}$ and $b_i = g(x_i)|N_i x^\alpha|_{H^1(\Omega)}$. Then
\[
v(x) = - \sum_{i=k}^N g(x_i)x_i^\alpha N_i(x) + \sum_{i=k}^N g(x_i) N_i(x) x^\alpha .
\]
Thus $v=0$ on $[0,x_k]$. Moreover, on $\tau_i$, $i \ge k+1$, we have
\[
v|_{\tau_i} = - \mcal{I}_h^i (gx^\alpha) + x^\alpha \mcal{I}_h^i(g),
\]
where $\mcal{I}_h^i(f)$ is the linear interpolant of $f$ on $\tau_i$, interpolating at $x_{i-1}$ and $x_i$.
Therefore,
\begin{eqnarray*}
|v|_{H^1(\open{\tau}_i)} &=& \big|gx^\alpha - \mcal{I}_h^i(gx^\alpha) - gx^\alpha + x^\alpha \mcal{I}_h^i(g)\big|_{H^1(\open{\tau}_i)} \\
&\le& Ch\big[ \, |gx^\alpha|_{H^2(\open{\tau}_i)} + \|x^\alpha\|_{H^2(\open{\tau}_i)} |g|_{H^2(\open{\tau}_i)}\big],
\end{eqnarray*}
where we have used standard interpolation estimates. Thus recalling that $v = 0$ on $[0,x_k]$, we have
\begin{eqnarray}
&&\hspace{-0.5cm}B(v,v) = |v|_{H^1(\Omega)}^2 = \sum_{i=k+1}^N |v|_{H^1(\open{\tau}_i)}^2 \nonumber \\
&&\hspace{0.0cm}\le Ch^2 \big[ |gx^\alpha|_{H^2(\widehat{\Omega})}^2 + \|x^\alpha\|_{H^1(\widehat{\Omega})}^2 \|g\|_{H^1(\widehat{\Omega})}^2    \big]:=Ch^2 |||gx^\alpha |||^2. \label{eqn2}
\end{eqnarray}
Also,
\[
\|\eta\|^2 \ge  \sum_{i=k}^N [g(x_i)]^2 |N_i x^\alpha|_{H^1(\Omega)}^2 \ge \sum_{i=\lceil N/2 \rceil}^N [g(x_i)]^2 |N_i x^\alpha|_{H^1(\Omega)}^2 \ge  \frac{C}{h^2},
\]
where we have used that $|N_i x^\alpha|_{H^1(\Omega)}^2 \ge C/h$ for $i\ge \lceil N/2 \rceil$. Thus using \eqref{eqn2}, we have
\[
\frac{B(v,v)}{\|\eta\|^2} \le Ch^4 |||gx^\alpha|||^2,
\]
and hence,
\[
\lambda_m \le Ch^4 |||gx^\alpha|||^2
\]
where $\lambda_m$ is the smallest eigenvalue of $H$. Finally, from \eqref{ev.eq1}, we get
\[
\mathfrak{K}(\mathbf{A}) = \kappa_2(H) = \frac{\lambda_M}{\lambda_m} \ge \frac{Ch^{-4}}{|||gx^\alpha|||^2},
\]
which is much bigger than the scaled condition number of the stiffness matrix of the standard FEM; we recall that the standard FEM is basic part of the GFEM. Thus from Hypothesis H, there will be severe loss of accuracy in the computed solution of $\mathbf{A}x=b$. We will show this feature in the Appendix.

It is interesting to note that using $v \in \mcal{S}$ of the form \eqref{unscaledv} and following the same arguments as before, we can also show that the condition number $\kappa_2(\mathbf{A}) \ge Ch^{-4}/|||gx^\alpha|||^2$. We stated this property at the beginning of this subsection. 

%% file: newSec4_SGFEMa.tex
\section{Stable Generalized Finite Element Method \\(SGFEM):} \label{SGFEM}

A GFEM will be referred to as an SGFEM if the GFEM satisfies the following property: the scaled condition number $\mathfrak{K}(\mathbf{A})$ of the associated stiffness matrix $\mathbf{A}$ is of the same order with respect to $h$ as of the stiffness matrix of the basic part of the GFEM. Since the basic part of any GFEM is the standard FEM, therefore a GFEM is an SGFEM provided $\mathfrak{K}(\mathbf{A})= O(h^{-2})$ for second order problems. As mentioned before, we will present the analysis for uniform meshes. However, the analysis is valid for locally quasi-uniform meshes.

We first present a particular example highlighting the ideas and results related to SGFEM in a simpler setting.

\subsection{An example of the SGFEM:} \label{example}

Let $a(x)=1$ in \eqref{VarProb} and suppose the solution of \eqref{VarProb} is smooth, in particular let $u \in H^3(\Omega)$. Since the solution is unique up to an additive constant, we seek $u$ with $u(0)=0$. It is well known that a function in $H^3(\Omega)$ could be accurately approximated, locally in $\omega_i$, by polynomials of degree 2; recall that the patches $\omega_i$ have been defined in Section \ref{GFEM}. Based on this information, we consider $V_i = \mbox{span}\{\varphi_j^{[i]}\}_{j=0}^2$ (i.e., $n_i=2$), where $\varphi_1^{[i]} =(x-x_i)$ and $\varphi_2^{[i]} =(x-x_i)^2$, for $0 \le i \le N$. Recall that $\varphi_0^{[i]} =1$. Thus $V_i = \mcal{P}^2(\omega_i)$.

We let
\[
\locphibar{j}{i} := \locphi{j}{i} - \intp{\omega_i}(\locphi{j}{i}),
\quad \mbox{where }
\intp{\omega_i}(\locphi{j}{i}) := \sum_{1-1 \le k \le i+1} \locphi{j}{i}(x_k) N_k\, \big|_{\omega_i};
 \]
$\intp{\omega_i}(\locphi{j}{i})$ is the piecewise linear interpolant of $\locphi{j}{i}$ on the patch $\omega_i$. We adjust the operators $\intp{\omega_0}$ and $\intp{\omega_1}$; they interpolate at $\{x_0,x_1\}$ and $\{x_{N-1},x_N\}$ respectively. We define a \emph{modified local approximation space} $\overline{V}_i = \mbox{span}\{\locphibar{j}{i}\}_{j=0}^{2}$, associated with $V_i$. Clearly, $\locphibar{j}{i}=0$ for $j=0,1$ and thus $\overline{V}_i = \mbox{span}\{\locphibar{2}{i}\}$.

It is well known (see \cite{MelBab,StrCopBab}) that the scaled condition number of the stiffness matrix of the GFEM, with $V_i$ as the local approximation spaces, could be extremely large or even unbounded. We will use the GFEM with $\overline{V}_i$ precisely to address this issue, and show that the GFEM based on the approximation space
\[
\mcal{S}=\mcal{S}_1 + \overline{\mcal{S}}_2, \quad \mbox{with }\
\mcal{S}_1=\sum_{i \in I\backslash \{0\}}a_iN_i \ \mbox{ and }\
\overline{\mcal{S}}_2 = \sum_{i\in I} N_i \, \overline{V}_i
\]
is an SGFEM. Note that $v(0)=0$ for all $v \in \mcal{S}$. We have chosen $a_0=0$ in the definition of $\mcal{S}_1$ to impose the constraint $u_h(0)=0$ to obtain a unique GFEM solution $u_h\in \mcal{S}$.

It is easy to check that the assumptions in Theorem \ref{GFEMApproxThm} hold; in fact, there exists $\xi_i \in V_i$ such that $\|u - \xi_i\|_{\mcal{E}(\omega_i)} \le Ch^2 |u|_{H^3(\omega_i)}$. Therefore it is clear from \eqref{ErrEstGFEM} that $\|u - u_h\|_{\mcal{E}(\Omega)} = \mcal{O}(h^2)$, where $u_h$ is the GFEM solution, based on $\mcal{S}=\mcal{S}_1 + \mcal{S}_2$ (recall $\mcal{S}_2=\sum_{i\in I} N_i  V_i$). We first show that the GFEM based on $\mcal{S}=\mcal{S}_1 + \overline{\mcal{S}}_2$ also yields the same optimal order of convergence.

\begin{myprop} \label{Approx1dEx2} There exists a $v \in \mcal{S}=\mcal{S}_1 + \overline{\mcal{S}}_2$, such that
\[
\|u-v\|_{\mcal{E}(\Omega)} \le Ch^2 |u|_{H^3(\Omega)},
\]
where the positive constant $C$ independent of $h$.
\end{myprop}

\emph{Proof:} Since  $u \in H^3(\omega_i)$ for $0 \le i \le N$, it is well known that there exists $\xi^i \in V_i = \mcal{P}^2(\omega_i)$ such that
\begin{equation}
\|u-\xi^i\|_{\mcal{E}(\omega_i)} \le Ch^{2}|u|_{H^3(\omega_i)} . \label{App1dEx2:1}
\end{equation}

Let $\ih u = \sum_{i \in I}u(x_i)N_i$. It is clear that $\ih u = \intp{\omega_i}u$ on $\omega_i$, Therefore using standard interpolation results, we have
\begin{eqnarray*}
\|(u - \ih u) - (\xi^i - \intp{\omega_i}\xi^i)\|_{L^2(\omega_i)} &=&
\|(u-\xi^i) - \intp{\omega_i}(u-\xi^i)\|_{L^2(\omega_i)} \\
&\le& C\mbox{diam}(\omega_i) \, \|(u-\xi^i)\|_{\mcal{E}(\omega_i)},
\end{eqnarray*}
and similarly,
\begin{eqnarray*}
\|(u - \ih u) - (\xi^i - \intp{\omega_i}\xi^i)\|_{\mcal{E}(\omega_i)}
&\le& C\|u-\xi^i\|_{\mcal{E}(\omega_i)} \le Ch^{2}|u|_{H^3(\omega_i)}.
\end{eqnarray*}
Let $w:=u-\ih u$; clearly $w \in H^1(\Omega)$. From above, $\xi^i - \intp{\omega_i}\xi^i \in \overline{V}_i$ approximates $w$ locally in $\omega_i$.
Therefore, from the Theorem \ref{GFEMApproxThm}, there is $\overline{v} \in \overline{\mcal{S}}_2$ such that
\begin{eqnarray}
\|w -\overline{v}\|_{\mcal{E}(\Omega)}^2
&\le& C^2\sum_{i\in I} h^4 |u|_{H^3(\omega_i)}^2 \le C^2 h^4 |u|_{H^3(\Omega)}^2. \label{App1dEx2:2}
\end{eqnarray}
Let $v = \ih u -u(x_0) + \overline{v}$. Since $\{N_i\}_{i\in I}$ is a PU, we have $\ih u - u(x_0) = \sum_{i\in I\backslash \{0\}}[u(x_i)-u(x_0)]N_i \in \mcal{S}_1$. Thus $v \in \mcal{S}$ and using \eqref{App1dEx2:2}, we get
\[
\|u -v\|_{\mcal{E}(\Omega)} = \|w -\overline{v}\|_{\mcal{E}(\Omega)} \le Ch^2 |u|_{H^3(\Omega)},
\]
which is the desired result. \myqed

Using Proposition \ref{Approx1dEx2}, we immediately get that $\|u-u_h\|_{\mcal{E}(\Omega)} = \mcal{O}(h^2)$, where $u_h$ is the GFEM solution based on $\mcal{S}=\mcal{S}_1 + \overline{\mcal{S}}_2$. We also note that we approximated $u - \ih u$ by the functions in $\overline{V}_i$ in the proof of Proposition \ref{Approx1dEx2} -- this is the main idea of SGFEM. Later, we will further comment on this issue.

We now address the scaled condition number of the stiffness matrix $\mathbf{A}$ associated with the GFEM based on $\mcal{S}=\mcal{S}_1 + \overline{\mcal{S}}_2$. With a suitable ordering of the shape function of $\mcal{S}$, the matrix $\mathbf{A}$ is of the form
\begin{equation}
\mathbf{A} = \left[ \begin{array}{cc}
        \mathbf{A}_{11} & \mathbf{A}_{12} \\ \mathbf{A}_{21} & \mathbf{A}_{22} \end{array} \right], \label{MGFEMstiff}
\end{equation}
where $\mathbf{A}_{ij}$ are block matrices. The matrix $\mathbf{A}_{11} = \{B(N_i,N_j)\}_{ i,j \in I\backslash \{0\} }$, which is the  stiffness matrix of the basic part of GFEM, is the standard $N\times N$ FE stiffness matrix. The $(N+1)\times (N+1)$ matrix $\mathbf{A}_{22}$ is of the form $\mathbf{A}_{22} = \{B(N_i\locphibar{2}{i},N_j\locphibar{2}{j})\}_{ i,j \in I }$. Also $\mathbf{A}_{21}=\mathbf{A}_{12}^T$. For the clarity of notation, we will write $\mathbf{A}_{22} = \{(\mathbf{A}_{22})_{ij}\}_{i,j=1}^M$, where $M=N+1$. Note that $(\mathbf{A}_{22})_{jj}$ are associated with the vertices $x_{j-1}$, respectively, and the GFEM solution $u_h$ is computed by postprocessing. We remark that, in general, $M$ will vary based on the application and $(\mathbf{A}_{22})_{jj}$ will be associated with some vertex $x_{i(j)}$.


We first note that $\overline{\varphi}_2^{[i]}(x_j) = 0$ for $j=i-1,i,i+1$. Therefore it is easy to show that $\mcal{S}_1$ and $\overline{\mcal{S}}_2$ are orthogonal in the inner product $B(\cdot,\cdot)$, i.e.,
\begin{equation}
B(v_1,v_2) = 0, \quad \forall \ v_1 \in \mcal{S}_1, \ v_2 \in \overline{\mcal{S}}_2. \label{Assump1Ex2}
\end{equation}
Thus it is immediate that $\mathbf{A}_{12}$ and $\mathbf{A}_{21}$ in \eqref{MGFEMstiff} are ``zero-matrices".

The matrix $\mathbf{A}_{11}$ is tridiagonal and is constructed by the assembly process from the element stiffness matrices $A_{11}^{(k)}$, for the element $\tau_k=[x_{k-1},x_k]$, $k=1,2,\cdots,N$. the matrices $A_{11}^{(k)}$ are given by
\begin{equation}
A_{11}^{(k)} = \frac{1}{h} \hat{A}_{11}^{(k)}; \quad
\hat{A}_{11}^{(k)}:= \left[ \begin{array}{rr} 1 & -1 \\ -1 & 1 \end{array} \right],\ 2 \le k \le N, \mbox{ and }  \hat{A}_{11}^{(1)}:= [1]. \label{ElmStifFE}
\end{equation}
Similarly, the matrix $\mathbf{A}_{22}$, which is also tridiagonal,  is constructed by the assembly process from the element matrices
\begin{equation}
A^{(k)}_{22} = \left[ \begin{array}{cc}
B_{\tau_k}(N_{k-1} \overline{\varphi}_2^{[k-1]},N_{k-1} \overline{\varphi}_2^{[k-1]}) & B_{\tau_k}(N_{k} \overline{\varphi}_2^{[k]},N_{k-1} \overline{\varphi}_2^{[k-1]})\\[1.8ex]
B_{\tau_k}(N_{k-1} \overline{\varphi}_2^{[k-1]},N_{k} \overline{\varphi}_2^{[k]}) & B_{\tau_k}(N_{k} \overline{\varphi}_2^{[k]},N_{k} \overline{\varphi}_2^{[k]})
\end{array}
\right], \label{ElemMatEx}
\end{equation}
for the element $\tau_k$, $k=1,2,\cdots, N$, and $B_{\tau_k}(w,v) := \int_{\tau_k}a u^\prime v^\prime\, dx$. A direct computation yields
\[
A_{22}^{(k)}= h^3 \hat{A}_{22}^{(k)};\quad  \hat{A}_{22}^{(k)} = \left[ \begin{array}{rr} \frac{2}{15} & \frac{1}{30} \\[1.0ex] \frac{1}{30} & \frac{2}{15} \end{array} \right].
\]
It is easy to check that the matrix $\hat{A}_{22}^{(k)}$ is positive definite (the eigenvalues are $\frac{1}{10}$ and $\frac{1}{6}$) and thus
\begin{equation}
\frac{h^3}{10} \|y\|^2 \le y^T A_{22}^{(k)} y \le \frac{h^3}{6} \|y\|^2,\quad \forall \ y=(y_1,y_2) \in \mathbb R^2. \label{Assump2Ex2}
\end{equation}

We now consider the diagonal matrix $\mathbf{D}=diag(\mathbf{D}_1,\mathbf{D}_2)$ with
\begin{eqnarray*}
&&\mathbf{D}_1=diag(\mathbf{d_1}),\ \mathbf{d_1}= m_1^{-1/2} (2^{-1/2},\cdots,2^{-1/2}, 1)^T \in \mathbb R^{N},\\
&&\mathbf{D}_2=diag(\mathbf{d_2}),\ \mathbf{d_2}= m_2^{-1/2} (1,2^{-1/2},\cdots,2^{-1/2},1)^T \in \mathbb R^{N+1},
\end{eqnarray*}
where $m_1 = 1/h$ and $m_2=2h^3/15$.

We next define
\[
\widehat{\mathbf{A}}:=\mathbf{D}\mathbf{A}\mathbf{D} = \left[ \begin{array}{cc}
       \mathbf{D}_1  \mathbf{A}_{11} \mathbf{D}_1 & \mathbf{0} \\ \mathbf{0} & \mathbf{D}_2\mathbf{A}_{22}\mathbf{D}_2 \end{array} \right]
        =  \left[ \begin{array}{cc}
        \widehat{\mathbf{A}}_{11} & \mathbf{0} \\ \mathbf{0} & \widehat{\mathbf{A}}_{22} \end{array} \right],
\]
where $\widehat{\mathbf{A}}_{11} = \mathbf{D}_1\mathbf{A}_{11}\mathbf{D}_1$ and $\widehat{\mathbf{A}}_{22} = \mathbf{D}_2\mathbf{A}_{22}\mathbf{D}_2$. Clearly $\widehat{\mathbf{A}}_{11}$ and $\widehat{\mathbf{A}}_{22}$ are $N \times N$ and $(N+1)\times (N+1)$ tri-diagonal matrices, respectively. The diagonal elements of $\widehat{\mathbf{A}}_{11}$ and $\widehat{\mathbf{A}}_{22}$ are equal to $1$. Consequently, diagonal elements of $\widehat{\mathbf{A}}$ are equal to $1$ and the scaled condition number $\mathfrak{K}(\mathbf{A})$ of $\mathbf{A}$ is $\kappa_2(\widehat{\mathbf{A}})$.

\begin{myprop} \label{CondNoEx2} Suppose $\mathfrak{K}(\mathbf{A})$ be the scaled condition number of $\mathbf{A}$ and let $\lambda_{min}(\widehat{\mathbf{A}}_{11})$, $\lambda_{max}(\widehat{\mathbf{A}}_{11})$ be the smallest and largest eigenvalue of $\widehat{\mathbf{A}}_{11}$, respectively. Then
\[
\mathfrak{K} (\mathbf{A}_{11}) \le \mathfrak{K}(\mathbf{A}) \le \mathfrak{K} (\mathbf{A}_{11}) \frac{\max \{1,C_2/\lambda_{max}(\widehat{\mathbf{A}}_{11}) }
{\min \{1,C_1/\lambda_{min}(\widehat{\mathbf{A}}_{11})\}},
\]
where $C_1=3/4$ and $C_2=5/4$.
\end{myprop}

\emph{Proof:} Let $\mathbf{z} = (\mathbf{z}_1,\mathbf{z}_2)^T \in \mathbb R^{2N+1}$, where $\mathbf{z}_1 \in \mathbb R^N$ and $\mathbf{z}_2 \in \mathbb R^{N+1}$. Then
\begin{eqnarray}
\mathbf{z}^T \widehat{\mathbf{A}} \mathbf{z} &=& (\mathbf{D}_1\mathbf{z}_1)^T \mathbf{A}_{11}(\mathbf{D}_1\mathbf{z}_1) + (\mathbf{D}_2\mathbf{z}_2)^T \mathbf{A}_{22}(\mathbf{D}_2\mathbf{z}_2) \nonumber \\
&=& \mathbf{z}_1^T \widehat{\mathbf{A}}_{11} \mathbf{z}_1 + \mathbf{z}_2^T \widehat{\mathbf{A}}_{22} \mathbf{z}_2. \label{Prop2Ex1:eq1}
\end{eqnarray}

Let $\mathbf{z}_2 = (y_1,y_2,\cdots,y_{N+1})^T$, then $\mathbf{D}_2 \mathbf{z}_2 = m_2^{-\frac{1}{2}} (y_1,2^{-\frac{1}{2}}y_2,\cdots,\\2^{-\frac{1}{2}}y_N,y_{N+1})^T$, where $m_2=2h^3/15$. We define $\mathbf{z}_{2,k} := (y_k,y_{k+1})^T$, and
\begin{eqnarray*}
&&\mathbf{\overline{z}}_{2,1} := m_2^{-1/2}(y_1,2^{-1/2}y_2)^T, \quad
\mathbf{\overline{z}}_{2,N}:=m_2^{-1/2}(2^{-1/2}y_N,y_{N+1})^T ,\\
&&\mathbf{\overline{z}}_{2,k}:= (2m_2)^{-1/2}(y_k,y_{k+1})^T ,\quad \mbox{for } k=2,\cdots,N-1.
\end{eqnarray*}
Recalling that $\mathbf{A}_{22}$ could be obtained from the element matrices $A_{22}^{(i)}$ through the assembly process, using \eqref{Assump2Ex2} we get,
\begin{eqnarray*}
\mathbf{z}_2^T \widehat{\mathbf{A}}_{22} \mathbf{z}_2 &=& (\mathbf{D}_2\mathbf{z}_2)^T \mathbf{A}_{22} (\mathbf{D}_2\mathbf{z}_2) = \sum_{k=1}^N {\mathbf{\overline{z}}_{2,k}}^T A_{22}^{(k)}\mathbf{\overline{z}}_{2,k} \\
&\le& \frac{h^3}{6}\sum_{k=1}^N\|\mathbf{\overline{z}}_{2,k}\|^2 = \frac{h^3}{6}\sum_{i=1}^{N+1} m_2^{-1}y_i^2 = \frac{5}{4} \|\mathbf{z}_2\|^2.
\end{eqnarray*}
Similarly from \eqref{Assump2Ex2}, we also get
\[
\frac{3}{4} \|\mathbf{z}_2\|^2 \le \mathbf{z}_2^T \widehat{\mathbf{A}}_{22} \mathbf{z}_2,
\]
and therefore from \eqref{Prop2Ex1:eq1},
\begin{equation}
\mathbf{z}_1^T \widehat{\mathbf{A}}_{11} \mathbf{z}_1 +C_1 \|\mathbf{z}_2\|^2 \le \mathbf{z}^T \widehat{\mathbf{A}} \mathbf{z} \le \mathbf{z}_1^T \widehat{\mathbf{A}}_{11} \mathbf{z}_1 +C_2 \|\mathbf{z}_2\|^2 . \label{Prop2Ex1:eq2}
\end{equation}
It is clear from above that
\begin{eqnarray*}
\mathbf{z}^T \widehat{\mathbf{A}} \mathbf{z} &\ge& \mathbf{z}_1^T \widehat{\mathbf{A}}_{11} \mathbf{z}_1 +C_1 \|\mathbf{z}_2\|^2 \\
&\ge& \lambda_{min}(\widehat{\mathbf{A}}_{11}) \|\mathbf{z}_1\|^2 + C_1 \|\mathbf{z}_2\|^2 \\
&\ge& \min \{C_1, \lambda_{min}(\widehat{\mathbf{A}}_{11})\} \|\mathbf{z}\|^2,
\end{eqnarray*}
where $C_1:=\frac{3}{4}$. Therefore,
\[
\lambda_{min}(\widehat{\mathbf{A}}) \ge \min \{C_1, \lambda_{min}(\widehat{\mathbf{A}}_{11})\} = \lambda_{min}(\widehat{\mathbf{A}}_{11}) \min \{1,C_1/\lambda_{min}(\widehat{\mathbf{A}}_{11})\}.
\]
Similarly from the upper bound of \eqref{Prop2Ex1:eq2}, we can show that
\[
\lambda_{max}(\widehat{\mathbf{A}}) \le \lambda_{max}(\widehat{\mathbf{A}}_{11}) \max \{1,C_2/\lambda_{max}(\widehat{\mathbf{A}}_{11}) \},
\]
where $C_2 = \frac{5}{4}$. Thus
\begin{eqnarray*}
\mathfrak{K}(\mathbf{A})= \kappa_2(\widehat{\mathbf{A}}) = \frac{ \lambda_{max}(\widehat{\mathbf{A}}) }
{ \lambda_{min}(\widehat{\mathbf{A}}) }
&\le& \frac{ \lambda_{max}(\widehat{\mathbf{A}}_{11}) \max \{1,C_2/\lambda_{max}(\widehat{\mathbf{A}}_{11}) \} }
{ \lambda_{min}(\widehat{\mathbf{A}}_{11}) \min \{1,C_1/\lambda_{min}(\widehat{\mathbf{A}}_{11})\} } \\
&=& \mathfrak{K} (\mathbf{A}_{11}) \frac{\max \{1,C_2/\lambda_{max}(\widehat{\mathbf{A}}_{11}) \} }
{\min \{1,C_1/\lambda_{min}(\widehat{\mathbf{A}}_{11})\}},
\end{eqnarray*}
where $\mathfrak{K} (\mathbf{A}_{11}) = \kappa_2(\widehat{\mathbf{A}}_{11})$. Thus we have the required upper bound of $\mathfrak{K}(\mathbf{A})$.

Now let $\mathbf{z}_1$ be an eigenvector of $\widehat{\mathbf{A}}_{11}$ associated with $\lambda_{max}(\widehat{\mathbf{A}}_{11})$. Also let $\mathbf{z}_2 = \mathbf{0}$. Then from \eqref{Prop2Ex1:eq1}, we have
\[
\mathbf{z}^T \widehat{\mathbf{A}} \mathbf{z} = \lambda_{max}(\widehat{\mathbf{A}}_{11})\, \|\mathbf{z}_1\|^2  = \lambda_{max}(\widehat{\mathbf{A}}_{11})\, \|\mathbf{z}\|^2,
\]
and therefore,
\[
\lambda_{max}(\widehat{\mathbf{A}}) \ge \lambda_{max}(\widehat{\mathbf{A}}_{11}).
\]
Similarly, considering $\mathbf{z}_1$ to be an eigenvector of $\widehat{\mathbf{A}}_{11}$ associated with $\lambda_{min}(\widehat{\mathbf{A}}_{11})$ and $\mathbf{z}_2=\mathbf{0}$, we have
\[
\mathbf{z}^T \widehat{\mathbf{A}} \mathbf{z} = \lambda_{min}(\widehat{\mathbf{A}}_{11})\, \|\mathbf{z}_1\|^2 =  \lambda_{min}(\widehat{\mathbf{A}}_{11})\, \|\mathbf{z}\|^2,
\]
and therefore,
\[
\lambda_{min}(\widehat{\mathbf{A}}) \le \lambda_{min}(\widehat{\mathbf{A}}_{11}).
\]
Now,
\[
\mathfrak{K}(\mathbf{A}) = \frac{ \lambda_{max}(\widehat{\mathbf{A}}) }
{ \lambda_{min}(\widehat{\mathbf{A}}) }
\ge \frac{ \lambda_{max}(\widehat{\mathbf{A}}_{11}) }{ \lambda_{min}(\widehat{\mathbf{A}}_{11}) }
= \mathfrak{K} (\mathbf{A}_{11}),
\]
which is the required lower bound of $\mathfrak{K}(\mathbf{A})$. \myqed

The Proposition \ref{CondNoEx2} establishes that $\mathfrak{K}(\mathbf{A}) \approx \mathfrak{K} (\mathbf{A}_{11})$, i.e., the scaled condition numbers of the stiffness matrices for the GFEM and the basic part of the GFEM are of same order. Thus the GFEM with $\mcal{S}=\mcal{S}_1+\overline{\mcal{S}}_2$ is indeed an SGFEM.

\begin{myremark} \upshape We note that the orthogonality of the spaces $\mcal{S}_1$ and $\overline{\mcal{S}}_2$ was essential in proving Proposition \ref{CondNoEx2}. This property does not hold in general. Later we will define a notion of ``almost orthogonality'' of $\mcal{S}_1$ and $\overline{\mcal{S}}_2$, which will address this issue. \myend
\end{myremark}

\begin{myremark} \upshape The inequality \eqref{Assump2Ex2} played an important role in obtaining Proposition \ref{CondNoEx2}. This property depends on the functions in $\overline{V}_i$. For general approximation spaces $\overline{V}_i$, we need an assumption that will be presented later. \myend
\end{myremark}

\begin{myremark} \upshape SGFEM uses $\overline{V}_i$ as the enrichment space, which is a modification of $V_i$. Other modifications of $V_i$ have been reported in different contexts. For example, the \emph{shifting} modification, namely,  $\locphibar{j}{i}(x) = \locphi{j}{i}(x) - \locphi{j}{i}(x_i)$, $j>0$, is used in XFEM in the context of approximation error as well as enforcement the \emph{Kronecker delta property} (see \cite{FriesBelRev}).\myend
\end{myremark} 

%% file: newSec4_SGFEMb.tex
\subsection{SGFEM and its analysis:}

We now present the SGFEM  for \eqref{VarProb}, with $a \in L^\infty(\Omega)$ and $0<\alpha \le a(x) \le \beta$. Moreover, suppose it is a priori known that $u(0)=0$ (we will further comment on a priori information later).
We consider the uniform mesh $\{\tau_k\}_{k\in I\backslash \{0\}}$ with the set of vertices $\mcal{T}$, as described in Section \ref{GFEM}. Recall that the hat function $N_i$ and the patch $\omega_i$ is associated with each $x_i \in \mcal{T}$. We will refer to $\{x_{i-1},x_i,x_{i+1}\}$ as the \emph{vertices} of $\omega_i$; the vertices of $\omega_0$, $\omega_N$ are $\{x_0,x_1\}$, $\{x_{N-1},x_N\}$, respectively.
Let
\[
\mcal{T}_1,\, \mcal{T}_2 \subset \mcal{T}; \quad \zeta_1 := \mbox{card}(\mcal{T}_1),\ \zeta_2 := \mbox{card}(\mcal{T}_2);\quad \zeta_1,\, \zeta_2 \le N+1.
\]
We define $\mcal{S}_1 = \sum_{x_i \in \mcal{T}_1}a_i N_i$, $a_i \in \mathbb R$; $\mcal{T}_1$ will be referred to as \emph{$\mcal{S}_1$-relevant set of vertices}. We consider $\mcal{T}_1=\{x_i\in \mcal{T}: 1\le i \le N\}$ as in the example in Section \ref{example}. For other choices of $\mcal{T}_1$, we refer to Remark \ref{RemT1T2}.


For $x_i \in \mcal{T}$, let $V_i = \mbox{span}\{\locphi{j}{i}\}_{j=0}^{n_i} \subset H^1(\omega_i)$ such that there exists $\xi^i \in V_i$ satisfying $\|u - \xi^i\|_{\mcal{E}(\open{\tau}_k)} \le \epsilon_i$ for all $\tau_k \subset \overline{\omega}_i$. Clearly, $\|u - \xi^i\|_{\mcal{E}(\omega_i)} \le 2\epsilon_i$. We consider the modified space
$\overline{V}_i = \mbox{span}\{\locphibar{j}{i}\}_{j=1}^{n_i}$, where
\[
\locphibar{j}{i} = \locphi{j}{i} - \intp{\omega_i}\locphi{j}{i};\quad
\intp{\omega_i}\locphi{j}{i}:= \sum_{i-1 \le k\le i+1}
\locphi{j}{i}(x_k)N_k \big|_{\omega_i}.
\]
$\intp{\omega_i}v$ is the piecewise linear interpolant of $v \in H^1(\omega_i)$ on the patch $\omega_i$ based on the vertices of $\omega_i$; we adjust $\intp{\omega_0}$ and $\intp{\omega_N}$ accordingly as before. It is important to note that if for some $x_i \in \mcal{T}$, $V_i = \{\xi \in H^1(\omega_i): \, \xi|_{\tau_k} \in \mcal{P}^1(\tau_k) \mbox{ for all } \tau_k \subset \overline{\omega}_i\}$, then $\overline{V}_i = \{0\}$.  Also $\overline{\xi}^i(x_k)=0$ with $ k=i-1,i,i+1$ for all $\overline{\xi}^i \in \overline{V}_i$. We refer to a patch $\omega_i$ as \emph{enriched} if $\overline{V}_i \ne \{0\}$. Let
$\mcal{T}_2:=\{x_i \in \mcal{T}:\, \omega_i \mbox{ is enriched}\}$ and define
$\overline{\mcal{S}}_2 = \sum_{x_i \in \mcal{T}_2}N_i\overline{V}_i$; $\mcal{T}_2$ will be referred to as the \emph{$\overline{\mcal{S}}_2$-relevant set of vertices}. In Section \ref{example}, we chose $\mcal{T}_2 = \mcal{T}$. We will present examples with $\zeta_2 << N+1$ (i.e., only few patches enriched) later in the paper.


\begin{myremark} \label{RemT1T2} \upshape
The sets $\mcal{T}_1,\, \mcal{T}_2 \subset \mcal{T}$ provide a framework to address numerical treatment of many applications. Selection of both sets depends on a priori information on the problem and its solution. Selection of $\mcal{T}_2$ will be apparent from the examples in Section \ref{APPLIC}. Suppose $\mcal{T}_0 \subset \mcal{T}$ contains all the vertices $x_j\in \mcal{T}$, where it is known a priori that $u(x_j) =0$. We choose $\mcal{T}_1 = \mcal{T}\backslash \mcal{T}_0$. Typically, $\mcal{T}_1$ will not contain any boundary vertex with homogeneous Dirichlet condition. However $\mcal{T}_1$ may exclude other vertices in $\mcal{T}$ based on a priori information. For example, let $f(x) = \sum_{k=0}^\infty c_k\cos [2\pi(2k+1) x]$, $a(x) = 1$, and suppose it it known that $u(0)=0$. Then $u(1/4)=u(3/4)=0$, and the vertices $x_j \notin \mcal{T}_1$ if $x_j=1/4$ or $x_j = 3/4$. Thus we can accommodate many a priori information in this framework. Only for simplicity, we have considered $\mcal{T}_0 = \{x_0\}$ in this section. \myend
\end{myremark}

We now consider a GFEM with
\begin{equation}
\mcal{S} = \mcal{S}_1 + \overline{\mcal{S}}_2 = \sum_{x_i \in \mcal{T}_1} a_i N_i + \sum_{x_i \in \mcal{T}_2}N_i \overline{V}_i \,.\label{SGFEMspaceMain}
\end{equation}
Note that $v(0)=0$ for all $v \in \mcal{S}$. We will show that this GFEM is an SGFEM, under certain assumptions on the space $\overline{\mcal{S}}_2$, which we will present later. We mention that $\mcal{T}_1$ and $\mcal{T}_2$ are called $\mcal{S}_1$ and $\overline{\mcal{S}}_2$ relevant vertices, respectively, since the degrees of freedom associated only with these vertices appear in the GFEM.

We first present an approximation result for the GFEM with $\mcal{S} = \mcal{S}_1 + \overline{\mcal{S}}_2$.

\begin{mytheorem} \label{ApproxThmSGFEM1d} Let $u\in H^1(\Omega)$ be the solution of \eqref{VarProb}. Suppose for each $x_i \in \mcal{T}_2$, there exists $\bar{\xi}^i \in \overline{V}_i$ and $C_1>0$, independent of $i$, such that
\[
\|u -\intp{\omega_i}u -\bar{\xi}^i\|_{L^2(\omega_i)} \le C_1 \mbox{diam}(\omega_i)\|u -\intp{\omega_i}u -\bar{\xi}^i\|_{\mcal{E}(\omega_i)},
\]
and
$
\|u -\intp{\omega_i}u -\bar{\xi}^i\|_{\mcal{E}(\omega_i)} \le \epsilon_i.
$
Then there exists $v \in \mcal{S}=\mcal{S}_1 + \overline{\mcal{S}}_2$ such that
\begin{equation}
\|u - v \|_{\mcal{E}(\Omega)} \le C \big\{ \sum_{x_i\in \mcal{T}\backslash \mcal{T}_2} \|u -\intp{\omega_i}u\|_{\mcal{E}(\omega_i)}^2 + \sum_{x_i\in \mcal{T}_2} \epsilon_i^2 \big\}^{1/2}. \label{SGFEMapprox}
\end{equation}
\end{mytheorem}

\emph{Proof:} Let $\ih u= \sum_{x_i \in \mcal{T}}u(x_i)N_i$ be the piecewise linear interpolant of $u$. We note that $\ih u = \intp{\omega_i}u$ on $\omega_i$. Define $w:=u-\ih u$ and let $\overline{v}:= \sum_{x_i \in \mcal{T}_2} N_i \bar{\xi}^i \in \overline{\mcal{S}}_2$. Then recalling that $\{N_i\}_{x_i \in \mcal{T}}$ is a PU, we have
\[
w - \overline{v} = \sum_{x_i \in \mcal{T}}N_iw - \sum_{x_i \in \mcal{T}_2}N_i \bar{\xi}^i = \sum_{x_i \in \mcal{T}\backslash \mcal{T}_2}N_iw + \sum_{x_i \in \mcal{T}_2}N_i(w-\bar{\xi}^i).
\]
Therefore
\begin{equation}
\|w - \overline{v}\|_{\mcal{E}(\Omega)}^2 \le C\big[ \big\|\sum_{x_i \in \mcal{T}\backslash \mcal{T}_2}N_iw \, \big\|_{\mcal{E}(\Omega)}^2 + \big\|\sum_{x_i \in \mcal{T}_2}N_i(w-\bar{\xi}^i) \, \big\|_{\mcal{E}(\Omega)}^2 . \label{SGFEMConvProof:1}
\end{equation}
We first address the last term of \eqref{SGFEMConvProof:1}. Using the fact that $x \in \Omega$ is in at most two patches $\omega_i$, $\omega_{i+1}$, we see that the sum $\sum_{x_i \in \mcal{T}_2}[N_i(w-\bar{\xi}^i)]^\prime$ has at most two terms for any $x \in \Omega$. Using this observation, the assumption that $\|N_i^\prime\|_{L^\infty(\Omega)} \le C[\mbox{diam}\{\omega_i\}]^{-1}$, and the hypothesis of the Theorem, we can show that
\begin{eqnarray}
\big\|\sum_{x_i \in \mcal{T}_2}N_i(w-\bar{\xi}^i) \, \big\|_{\mcal{E}(\Omega)}^2 &\le& C \Big[ \sum_{x_i \in \mcal{T}_2} \frac{\|w-\bar{\xi}^i\|^2_{L^2(\omega_i)}}{\mbox{diam}\{\omega_i\}^2} \nonumber \\
&&\hspace{1cm}+  \sum_{x_i \in \mcal{T}_2} \|w-\bar{\xi}^i\|^2_{\mcal{E}(\omega_i)} \Big] \nonumber \\
&\le& \sum_{x_i \in \mcal{T}_2} \|w-\bar{\xi}^i\|^2_{\mcal{E}(\omega_i)} \le \sum_{x_i \in \mcal{T}_2} \epsilon_i^2. \label{SGFEMConvProof:2}
\end{eqnarray}
(We refer to the proof of Theorem 3.2 in \cite{BBO7} for details of the argument leading to \eqref{SGFEMConvProof:2}). Using exactly same argument and the interpolation estimate $\|w\|_{L^2(\omega_i)}=\|u - \intp{\omega_i}u\|_{L^2(\omega_i)} \le Ch\|u - \intp{\omega_i}u\|_{\mcal{E}(\omega_i)}$, we get
\[
\big\|\sum_{x_i \in \mcal{T}\backslash \mcal{T}_2}N_iw \, \big\|_{\mcal{E}(\Omega)}^2 \le C \sum_{x_i \in \mcal{T}\backslash \mcal{T}_2} \|u-\intp{\omega_i}u\|_{\mcal{E}(\omega_i)}^2.
\]
Therefore, from \eqref{SGFEMConvProof:1} and \eqref{SGFEMConvProof:2}, we have
\[
\|w - \overline{v}\|_{\mcal{E}(\Omega)}^2 \le C\big[\sum_{x_i \in \mcal{T}\backslash \mcal{T}_2} \|u-\intp{\omega_i}u\|_{\mcal{E}(\omega_i)}^2 + \sum_{x_i \in \mcal{T}_2} \epsilon_i^2 \big].
\]
Finally, writing $w = u - \ih u$ and setting $v = \ih u + \overline{v} \in \mcal{S}_1 + \overline{\mcal{S}}_2$, we get the desired result. \myqed

We mention that unlike in Theorem \ref{Approx1dEx2}, we did not assume $\mcal{T}_2 = \mcal{T}$ in Theorem \ref{ApproxThmSGFEM1d}.
We further note that $\ih u$ for $u\in H^1(\Omega)$ is not defined in higher dimensions, since the point values of $u$, in general, do not exist in higher dimensions (in contrast to 1-d). However, using a generalized interpolant based on the average of $u$ in a ball around the vertices $x_i$, the proof of the above result can be easily generalized to higher dimensions.


\begin{myremark} \upshape From the proof of Proposition \ref{Approx1dEx2}, it is clear that accurate local approximation of $u - \ih u$ by functions in $\overline{V}_i$ is crucial to obtain the desired result. This is the main idea of SGFEM -- the spaces $\overline{V}_i$ are constructed such that the functions in $\overline{V}_i$ accurately approximate $u - \ih u$ in $\omega_i$. This is in contrast to the standard GFEM, where the functions in local approximating spaces $V_i$ accurately approximate $u$ in $\omega_i$. \myend
\end{myremark}

\begin{myremark} \upshape We note that $\overline{V}_i = \{0\}$ for $x_i \in \mcal{T}\backslash \mcal{T}_2$. If $u \in H^1(\Omega)$ is locally smooth, namely, $u \in H^2(\omega_i)$ for $x_i \in \mcal{T}\backslash \mcal{T}_2$, then $\|u-\intp{\omega_i}u\|_{\mcal{E}(\omega_i)} \le Ch|u|_{H^2(\omega_i)}$ for $x_i \in \mcal{T}\backslash \mcal{T}_2$ and \eqref{SGFEMapprox} could be written as
\begin{equation}
\|u - v \|_{\mcal{E}(\Omega)} \le C \big\{ h^2\sum_{x_i\in \mcal{T}\backslash \mcal{T}_2} |u|_{H^2(\omega_i)}^2 + \sum_{x_i\in \mcal{T}_2} \epsilon_i^2 \big\}^{1/2}. \label{RemApprox:1}
\end{equation}
By incorporating the available information on the solution $u$ in $V_i$, for $x_i \in \mcal{T}_2$, we can have $\epsilon_i = O(h)$, and consequently, $\|u - v \|_{\mcal{E}(\Omega)} = O(h)$. The set $\mcal{T}_2$ can be chosen adaptively with respect to a prescribed tolerance, which we do not elaborate in this paper. \myend
\end{myremark}

\begin{myremark} \upshape
A rate of convergence of $O(h)$ for various problems have been reported for the Corrected XFEM (which is also a GFEM); see e.g., \cite{TFries}.
However, for the crack propagation problems, the enrichment spaces $V_i$ in XFEM requires the use of a \emph{ramp-function} to obtain the $O(h)$ rate of convergence. In contrast, the GFEM based on $\mcal{S} = \mcal{S}_1+\overline{\mcal{S}}_2$ does not require the use of a ramp-function to obtain the rate of convergence of $O(h)$.
\end{myremark}

We now address the scaled condition number of the stiffness matrix of the GFEM with $\mcal{S} = \mcal{S}_1 + \overline{\mcal{S}}_2$. For clarity of the exposition, we will present the analysis for the case when $n_i = 1$  i.e.,  $\overline{V}_i = \mbox{span}\{\locphibar{1}{i}\}$ . The analysis for general $n_i$ is similar.

As in the example presented in Section \ref{example}, the stiffness matrix $\mathbf{A}$ is of the form
$
\mathbf{A} = \textstyle{ \left[ \begin{array}{cc}
                \mathbf{A}_{11} & \mathbf{A}_{12}\\
                \mathbf{A}_{21} & \mathbf{A}_{22} \end{array} \right] },
$
where $\mathbf{A}_{11} = \{B(N_i,N_j)\}_{x_i,x_j\in \mcal{T}_1}$ is the $\zeta_1 \times \zeta_1$ stiffness matrix of the basic part of the GFEM. Let $\mathbf{D}_1$ be a diagonal matrix with $(\mathbf{D}_1)_{ii} = (\mathbf{A}_{11})_{ii}^{-1/2}$. Clearly, the diagonal elements of
\begin{equation}
\widehat{\mathbf{A}}_{11}:= \mathbf{D}_1 \mathbf{A}_{11} \mathbf{D}_1 \label{A11hat}
\end{equation}
are equal to 1.

The matrix $\mathbf{A}_{22}$ plays a central role in our analysis and depends on elements that have been enriched. We will refer to an element $\tau_k=[x_{k-1},x_k]$ as \emph{enriched} if (a) $x_{k-1}\in \mcal{T}_2$ and $\locphibar{1}{k-1}|_{\tau_k} \not \equiv 0$, or (b) $x_k\in \mcal{T}_2$ and $\locphibar{1}{k}|_{\tau_k} \not \equiv 0$. Let
\[
\mcal{K}_{enr} := \{ \tau_k :\, \tau_k \mbox{ is enriched} \}.
\]
The matrix $\mathbf{A}_{22}$ is constructed by the assembly process using the element stiffness matrices $A_{22}^{(k)}$ defined only on $\tau_k \in \mcal{K}_{enr}$.

We now address the structure of the element matrices $A_{22}^{(k)}$ in detail and set up some notions and notations that will be used in the analysis. We denote the vertices of the element $\tau_k$ as $x_1^{(k)} := x_{k-1}$ and $x_2^{(k)}:= x_k$; we consider only $\tau_k \in \mcal{K}_{enr}$. The element stiffness matrix $A_{22}^{(k)}$  is of the form
\begin{equation}
A_{22}^{(k)}= \left[ \begin{array}{cc} b_{11}^{(k)} & b_{12}^{(k)} \\[.8ex] b_{12}^{(k)} & b_{22}^{(k)} \end{array} \right], \label{LocStiffMat}
\end{equation}
where $b_{ij}^{(k)} = B_{\tau_k} \big( N_{k-2+i}\locphibar{1}{k-2+i}, N_{k-2+j}\locphibar{1}{k-2+j} \big)$, $ 1 \le i,j \le 2$.

If $b_{11}^{(k)}, b_{22}^{(k)} >0$, then $A_{22}^{(k)}$ is $2 \times 2$ and we say that the local stiffness matrix $A_{22}^{(k)}$ is \emph{associated with the vertices} $x_1^{(k)}=x_{k-1}$ and $x_2^{(k)}=x_k$. We define a diagonal matrix $D^{(k)} = \mbox{diag}\{\delta_1^{(k)}, \delta_2^{(k)}\}$ with $\delta_1^{(k)}, \delta_2^{(k)}>0$, such that the diagonal elements of
\[
\hat{A}_{22}^{(k)}:= D^{(k)} A_{22}^{(k)} D^{(k)}
\]
are of equal to 1 or $O(1)$, independent of $h$. Clearly, $\delta_1^{(k)}, \delta_2^{(k)}$ are associated with vertices $x_{k-1},x_k$ respectively.

On the other hand, if $b_{22}^{(k)} = 0$ (i.e., $\locphibar{1}{k}|_{\tau_k} \equiv 0$  and consequently $b_{12}^{(k)}=b_{21}^{(k)} =0$) in \eqref{LocStiffMat}, then the local stiffness matrix $A_{22}^{(k)} = [b_{11}^{(k)}]$ is of size $1 \times 1$ and is associated only with the vertex $x_1^{(k)} = x_{k-1}$. We define $D^{(k)} = [\delta_1^{(k)}]$, where $\delta_1^{(k)} = \{b_{11}^{(k)}\}^{-1/2}$; $\delta_1^{(k)}$ is associated with the vertex $x_1^{(k)}=x_{k-1}$. Similarly, if $b_{11}^{(k)} = 0$ in \eqref{LocStiffMat}, then the local stiffness matrix $A_{22}^{(k)} = [b_{22}^{(k)}]$ is associated only with the vertex $x_2^{(k)} = x_{k}$. Also $D^{(k)} = [\delta_2^{(k)}]$ with $\delta_1^{(k)} = \{b_{22}^{(k)}\}^{-1/2}$ associated with the vertex $x_2^{(k)}=x_{k}$. Let $\varsigma^{(k)}$ be the number of vertices associated with the local stiffness matrix $A_{22}^{(k)}$. Thus the size of $A_{22}^{(k)}$ is $\varsigma^{(k)} \times \varsigma^{(k)}$; note that $\varsigma^{(k)}$ is either $1$ or $2$ with our assumption $n_i =1$.

Recall that $\mathbf{A}_{22}$ is obtained by the assembly process using the element stiffness matrices $A_{22}^{(k)}$; the size of $\mathbf{A}_{22}$ is $\zeta_2 \times \zeta_2$. Let $c=(c_1,c_2,\cdots, c_{\zeta_2})$, then
\begin{equation}
c^T\mathbf{A}_{22}c = \sum_{\tau_k \in \mcal{K}_{enr}} [c^{(k)}]^T A_{22}^{(k)}c^{(k)}, \label{assembly}
\end{equation}
where $c^{(k)} \in \mathbb R^{\varsigma^{(k)}}$ . Moreover, the components of $c^{(k)}$ are also the components of $c$ that correspond to those vertices of $\tau_k$  that are associated with $A_{22}^{(k)}$. For example, if $b_{11}^{(k)}, b_{22}^{(k)} >0$ in $A_{22}^{(k)}$, then as mentioned before, the vertices $x_1^{(k)}=x_{k-1}$, $x_2^{(k)}=x_k$ are associated with $A_{22}^{(k)}$. Suppose the components $c_{j(k)-1},c_{j(k)}$ of $c$ are associated with the vertices $x_{k-1},x_k$, respectively,  of $\tau_k$. Then $c^{(k)} = [c_{j(k)-1},c_{j(k)}]^T$. Similarly, if $A_{22}^{(k)} = [b_{11}^{(k)}]$, then $A_{22}^{(k)}$ is associated with $x_1^{(k)}$ and $c^{(k)} = [c_{j(k)-1}]$ -- a vector with one component. Later in our analysis, we will use \eqref{assembly} with a particular vector $c$ and $c^{(k)}$ as defined above.

Next we note that each vertex $x_i$ of the FE mesh is associated with a FE star -- union of all elements $\tau_k\subset \overline{\omega}_i$ (equivalently, union of all elements $\tau_k$ with common vertex $x_i$). For $x_i \in \mcal{T}_2$, we define
$
\mcal{K}_i := \{\tau_k \in \mcal{K}_{enr}:\, \tau_k \subset \overline{\omega}_i \}
$.
For $x_i \in \mcal{T}_2$ and $\tau_k \in \mcal{K}_i$, we set the index $1 \le l(i,k) \le 2$ as follows. We first note that $k \in \{i,i+1\}$. For $k=i$, we set $l(i,k)=l(i,i) = 2$ and for $k=i+1$, we set $l(i,k)=l(i,i+1)=1$. Thus $l(i,k)$ is the index such that $x_{l(i,k)}^{(k)} = x_i$; note $x_{l(i,k)}^{(k)}$ may not be associated with $A_{22}^{(k)}$. We define
\[
\mcal{K}_i^*:= \{\tau_k \in \mcal{K}_i:\, x_{l(i,k)}^{(k)} \mbox{ is associated with } A_{22}^{(k)}\}.
\]
Thus $\mcal{K}_i^*$ is the set of  $\tau_k \in \mcal{K}_i$ such that $\locphibar{1}{i}|_{\tau_k} \not\equiv 0$. For $x_i \in \mcal{T}_2$, we define
\begin{equation}
\Delta_i := \sum_{\tau_k \in \mcal{K}_i^*} [\delta_{l(i,k)}^{(k)}]^{-2}, \label{Ai}
\end{equation}
which will be used later in our analysis.

Each diagonal element of $\mathbf{A}_{22}$ is associated with a vertex in $\mcal{T}_2$. Let $(\mathbf{A}_{22})_{j_ij_i}$ be associated with $x_i \in \mcal{T}_2$. Moreover, we note that $(\mathbf{A}_{22})_{j_ij_i} = \sum_{\tau_k \in \mcal{K}_i^*} b^{(k)}_{l(i,k),l(i,k)}$, where $b^{(k)}_{pq}$ was defined in \eqref{LocStiffMat}. Thus $(\mathbf{A}_{22})_{j_ij_i}>0$ for all $x_i \in \mcal{T}_2$ (i.e., all the diagonal elements of $\mathbf{A}_{22}$ are positive). We now define the diagonal matrix $\mathbf{D}_2 =\mbox{diag}\{d_1,d_2,\cdots,d_{\zeta_2}\}$ with $d_{j} = (\mathbf{A}_{22})_{jj}^{-1/2}$, $1 \le j \le \zeta_2$. Note that $d_{j_i}=(\mathbf{A}_{22})_{j_ij_i}^{-1/2}$ is associated with $x_i \in \mcal{T}_2$. Clearly, the diagonal elements of
\begin{equation}
\widehat{\mathbf{A}}_{22}:= \mathbf{D}_2 \mathbf{A}_{22}  \mathbf{D}_2 \label{A22hat}
\end{equation}
are equal to 1. Define the diagonal matrix $\mathbf{D} := \mbox{diag}\{\mathbf{D}_1,\mathbf{D}_2\}$. Since the diagonal elements of
$\widehat{\mathbf{A}}_{11}$, $\widehat{\mathbf{A}}_{22}$ (see \eqref{A11hat}, \eqref{A22hat}) are equal to 1, the diagonal elements of
\begin{equation}
\widehat{\mathbf{A}} := \mathbf{D} \mathbf{A} \mathbf{D} =
\left[ \begin{array}{cc} \widehat{\mathbf{A}}_{11} & \widehat{\mathbf{A}}_{12} \\[.8ex] \widehat{\mathbf{A}}_{21} & \widehat{\mathbf{A}}_{22} \end{array} \right] \label{ScldStiff}
\end{equation}
are also equal to 1. Also $\widehat{\mathbf{A}}_{12} = \mathbf{D}_1 \mathbf{A}_{12} \mathbf{D}_2$ and $\widehat{\mathbf{A}}_{21} = \widehat{\mathbf{A}}_{12}^T$.

We will show that the GFEM with $\mcal{S}=\mcal{S}_1 + \overline{\mcal{S}}_2$ is an SGFEM, under the following assumptions on the local approximation spaces $\overline{V}_i$ and the enrichment part of $\mcal{S}$, namely $\overline{S}_2$.

\begin{assump} \label{asp1} The spaces $\mcal{S}_1$ and $\overline{\mcal{S}}_2$ are \emph{almost orthogonal} with respect to the inner product $B(\cdot.\cdot)$, i.e., there exist constants $0<L_1, U_1< \infty$, independent of $h$, such that
\[
L_1 \big\{ \|v_1\|_{\mcal{E}(\Omega)}^2 + \|v_2\|_{\mcal{E}(\Omega)}^2 \big\}\le |B(v_1+v_2,v_1+v_2)| \le U_1 \big\{ \|v_1\|_{\mcal{E}(\Omega)}^2 + \|v_2\|_{\mcal{E}(\Omega)}^2 \big\},
\]
for all $v_1 \in \mcal{S}_1$ and $v_2 \in \overline{\mcal{S}}_2$.
\end{assump}

\begin{assump} \label{asp2} For $\tau_k \in \mcal{K}_{enr}$, there exist constants $0< L_2,U_2 < \infty$, independent of $k$ and $h$ such that
\[
 L_2\|[D^{(k)}]^{-1} \mathbf{x}\|^2 \le \mathbf{x}^T A_{22}^{(k)} \mathbf{x} \le U_2 \|[D^{(k)}]^{-1} \mathbf{x}\|^2, \quad \forall \ \mathbf{x} \in \mathbb R^{\varsigma^{(k)}},
\]
where the diagonal matrices $D^{(k)}$ have been defined before.
\end{assump}

\begin{assump} \label{asp3}  For $x_i \in \mcal{T}_2$, there exist constants $0< L_3, U_3 < \infty$, independent of $i$ and $h$ such that
\[
L_3 \le (\mathbf{A}_{22})_{j_ij_i}^{-1} \Delta_i \le U_3,
\]
where $(\mathbf{A}_{22})_{j_ij_i}$ is the diagonal element of $\mathbf{A}_{22}$ associated with $x_i$, and $\Delta_i$ is as defined in \eqref{Ai}.
\end{assump}


The following result is an easy consequence of Assumption \ref{asp1}.
\begin{mylemma} \label{SplitStiff} Let $x=(\xi^T,\eta^T)^T \in \mathbb R^{\zeta_1+\zeta_2}$ where $\xi \in \mathbb R^{\zeta_1}$ and $\eta \in \mathbb R^{\zeta_2}$.
Then there exist positive constants $L_1$ and $U_1$, independent of $h$, such that
\[
L_1 \left[ \xi^T \mathbf{A}_{11} \xi + \eta^T \mathbf{A}_{22} \eta \right] \le x^T \mathbf{A} x \le U_1 \left[ \xi^T \mathbf{A}_{11} \xi + \eta^T \mathbf{A}_{22} \eta \right],
\]
where $\mathbf{A}$, $\mathbf{A}_{11}$ and $\mathbf{A}_{22}$ are matrices defined before.
\end{mylemma}

\noindent \emph{Proof:} Let $\xi = (\xi_i)_{x_i \in \mcal{T}_1}$ and $\eta = (\eta_i)_{x_i \in \mcal{T}_2}$. Consider $v_1 = \sum_{x_i\in \mcal{T}_1} \xi_i N_i \in \mcal{S}_1$ and $v_2 = \sum_{x_i \in \mcal{T}_2} \eta_i N_i \locphibar{1}{i} \in \overline{\mcal{S}}_2$. Then $B(v_1+v_2, v_1+v_2) = x^T \mathbf{A} x$, $B(v_1,v_1)= \xi^T \mathbf{A}_{11} \xi$, and $B(v_2,v_2) = \eta^T \mathbf{A}_{22} \eta$. The desired result is now immediate from Assumption \ref{asp1}. \myqed

\begin{mytheorem} \label{ThCondNo}
Suppose the Assumptions \ref{asp1}, \ref{asp2}, and \ref{asp3} are satisfied. Let $\mathbf{A}$ be the stiffness matrix of the GFEM with $\mcal{S}=\mcal{S}_1 + \overline{\mcal{S}}_2$. Then
\[
\frac{L_1}{U_1} \, \mathfrak{K}(\mathbf{A}_{11}) \le \mathfrak{K}(\mathbf{A}) \le \mathfrak{K}(\mathbf{A}_{11})\, \frac{U_1}{L_1} \, \frac{\max \big\{1,U_2 U_3 / \lambda_{max}(\widehat{\mathbf{A}}_{11}) \big\}}{\min \big\{1,L_2 L_3 / \lambda_{min}(\widehat{\mathbf{A}}_{11}) \big\}}\, ,
\]
where $\lambda_{min}(\widehat{\mathbf{A}}_{11})$, $\lambda_{max}(\widehat{\mathbf{A}})$ are the smallest and largest eigenvalues, respectively, of the matrix $\widehat{\mathbf{A}}_{11}$ defined before.
\end{mytheorem}

\begin{myremark} \upshape This result shows that under the Assumptions \ref{asp1}, \ref{asp2},  and \ref{asp3}, the scaled condition numbers of the stiffness matrices of the GFEM with $\mcal{S}=\mcal{S}_1 + \overline{\mcal{S}}_2$ and the basic part of the GFEM are of the same order. Thus the GFEM with $\mcal{S}=\mcal{S}_1 + \overline{\mcal{S}}_2$ is indeed an SGFEM.
\end{myremark}

\noindent \emph{Proof:} Let $\mathbf{z} = (\mathbf{z}_1,\mathbf{z}_2)^T \in \mathbb R^{\zeta_1+\zeta_2}$, where $\mathbf{z}_1\in \mathbb R^{\zeta_1}$ and $\mathbf{z}_2 \in \mathbb R^{\zeta_2}$. Then from the definition of $\widehat{\mathbf{A}}$ (see \eqref{ScldStiff}), we have
$
\mathbf{z}^T \widehat{\mathbf{A}} \mathbf{z} = \mathbf{z}^T \mathbf{D A D} \mathbf{z} = (\mathbf{Dz})^T \mathbf{A} (\mathbf{Dz})
$,
and since $\mathbf{Dz} = \left[(\mathbf{D}_1\mathbf{z}_1)^T, (\mathbf{D}_2\mathbf{z}_2)^T\right]^T$, from Lemma \ref{SplitStiff} we get
\begin{eqnarray}
L_1 \left[ (\mathbf{D}_1\mathbf{z}_1)^T \mathbf{A}_{11} (\mathbf{D}_1\mathbf{z}_1) + (\mathbf{D}_2\mathbf{z}_2)^T \mathbf{A}_{22} (\mathbf{D}_2\mathbf{z}_2) \right] &\le& \mathbf{z}^T \widehat{\mathbf{A}} \mathbf{z} \nonumber\\
&&\hspace{-7cm} \le U_1\left[ (\mathbf{D}_1\mathbf{z}_1)^T \mathbf{A}_{11} (\mathbf{D}_1\mathbf{z}_1) + (\mathbf{D}_2\mathbf{z}_2)^T \mathbf{A}_{22} (\mathbf{D}_2\mathbf{z}_2) \right] .\label{ThConNo:1}
\end{eqnarray}

Let $\mathbf{z}_2 = (f_1,f_2,\cdots,f_{\zeta_2})^T$ and consider $\mathbf{D}_2 = \mbox{diag}(d_1,d_2,\cdots,d_{\zeta_2})$ with $d_i= (\mathbf{A}_{22})_{ii}^{-1/2}$ as defined before. Then $\mathbf{D}_2 \mathbf{z}_2 = (d_1f_1,d_2f_2,\cdots,d_{\zeta_2}f_{\zeta_2})^T$. Recall that $d_{j_i}$ is associated with $x_i \in \mcal{T}_2$. Consequently, $d_{j_i}f_{j_i}$ is associated with $x_i \in \mcal{T}_2$.

Consider an element $\tau_k \in \mcal{K}_{enr}$. Following the notation given after \eqref{assembly}, let $\mathbf{\overline{z}}_2^{(k)}:= (\mathbf{D}_2 \mathbf{z}_2)^{(k)} \in \mathbb R^{\varsigma^{(k)}}$ such that the components of $\mathbf{\overline{z}}_2^{(k)}$ are the components of $\mathbf{D}_2\mathbf{z}_2$ corresponding to the vertices of $\tau_k$ associated with $A_{22}^{(k)}$. Now from \eqref{assembly} and using Assumption \ref{asp2}, we have
\begin{equation}
(\mathbf{D}_2\mathbf{z}_2)^T \mathbf{A}_{22} (\mathbf{D}_2\mathbf{z}_2) = \sum_{\tau_k \in \mcal{K}_{enr}} {\mathbf{\overline{z}}_2^k}^T A_{22}^{(k)}\mathbf{\overline{z}}_2^k \ge L_2\sum_{\tau_k \in \mcal{K}_{enr}} \|[D^{(k)}]^{-1}\mathbf{\overline{z}}_2^k\|^2 .\label{MainThm:eq1}
\end{equation}
We note that if $D^{(k)} = \mbox{diag}\{\delta_1^{(k)},\delta_2^{(k)}\}$, then \[
\|[D^{(k)}]^{-1}\mathbf{\overline{z}}_2^k\|^2 = [\delta_1^{(k)}]^{-2} \, [d_{j(k)-1}]^2 \, [f_{j(k)-1}]^2 + [\delta_2^{(k)}]^{-2}\, [d_{j(k)}]^2 \, [f_{j(k)}]^2,
\]
where $\mathbf{\overline{z}}_2^k = [d_{j(k)-1}, f_{j(k)-1}]^T$ following the notation given after \eqref{assembly}. Similarly, if $D^{(k)} = [\delta_1^{(k)}]$, then $\|[D^{(k)}]^{-1}\mathbf{\overline{z}}_2^k\|^2 = [\delta_1^{(k)}]^{-2}[d_{j(k)-1}]^2 [f_{j(k)-1}]^2$, and if $D^{(k)} = [\delta_2^{(k)}]$, then $\|[D^{(k)}]^{-1}\mathbf{\overline{z}}_2^k\|^2 = [\delta_2^{(k)}]^{-2}[d_{j(k)}]^2 [f_{j(k)}]^2$.

Now, it is important to note that if $\mcal{J}_1 := \{d_{j(k)-1} f_{j(k)-1}, d_{j(k)} f_{j(k)}\}_{\tau_k \in \mcal{K}_{enr}}$ (where the repeated elements appear only once) and $\mcal{J}_2 := \{d_{j_i}f_{j_i}\}_{x_i \in \mcal{T}_2}$, then $\mcal{J}_1=\mcal{J}_2$. Thus from \eqref{MainThm:eq1}, we have
\begin{equation}
(\mathbf{D}_2\mathbf{z}_2)^T \mathbf{A}_{22} (\mathbf{D}_2\mathbf{z}_2) \ge L_2 \sum_{x_i \in \mcal{T}_2} \Delta_i d_{j_i}^2 f_{j_i}^2 \ge L_2L_3\|\mathbf{z}_2\|^2, \label{ThConNo:1a}
\end{equation}
where we used Assumption \ref{asp3} to get the last inequality.
Similarly, we can show that
\[
(\mathbf{D}_2\mathbf{z}_2)^T \mathbf{A}_{22} (\mathbf{D}_2\mathbf{z}_2) \le U_2 U_3 \|\mathbf{z}_2\|^2.
\]
Therefore from \eqref{ThConNo:1} and using the definition of $\widehat{\mathbf{A}}_{11}$, we get
\begin{equation}
L_1 \left[ \mathbf{z}_1^T \widehat{\mathbf{A}}_{11} \mathbf{z}_1 + L_2 L_3 \|\mathbf{z}_2 \|^2 \right] \ \le \ \mathbf{z}^T \widehat{\mathbf{A}} \mathbf{z}
\le U_1\left[ \mathbf{z}_1^T \widehat{\mathbf{A}}_{11} \mathbf{z}_1 + U_2 U_3 \|\mathbf{z}_2\|^2 \right] .\label{ThConNo:2}
\end{equation}

Now from the lower bound of $\mathbf{z}^T \widehat{\mathbf{A}} \mathbf{z}$ in \eqref{ThConNo:2}, we have
\begin{eqnarray*}
\mathbf{z}^T \widehat{\mathbf{A}} \mathbf{z} &\ge& L_1 \left[ \lambda_{min}(\widehat{\mathbf{A}}_{11})\|\mathbf{z}_1\|^2 + L_2 L_3 \|\mathbf{z}_2\|^2 \right], 
\end{eqnarray*}
and therefore,
\begin{equation}
\lambda_{min}(\widehat{\mathbf{A}}) \ge L_1 \, \lambda_{min}(\widehat{\mathbf{A}}_{11}) \, \min \big\{1,L_2 L_3/ \lambda_{min}(\widehat{\mathbf{A}}_{11}) \big\}. \label{ThConNo:3}
\end{equation}
Similarly, using the upper bound of $\mathbf{z}^T \widehat{\mathbf{A}} \mathbf{z}$ in \eqref{ThConNo:2}, we can show
\begin{equation}
\lambda_{max}(\widehat{\mathbf{A}}) \le U_1 \, \lambda_{max}(\widehat{\mathbf{A}}_{11}) \, \max \big\{1,U_2 U_3 / \lambda_{max}(\widehat{\mathbf{A}}_{11}) \big\}. \label{ThConNo:4}
\end{equation}
Thus from \eqref{ThConNo:3} and \eqref{ThConNo:4}, we have
\begin{equation}
\mathfrak{K}(\mathbf{A}) = \frac{\lambda_{max}(\widehat{\mathbf{A}})}{\lambda_{min}(\widehat{\mathbf{A}})}
\le \mathfrak{K}(\mathbf{A}_{11})\, \frac{U_1}{L_1} \, \frac{\max \big\{1,U_2 U_3/ \lambda_{max}(\widehat{\mathbf{A}}_{11}) \big\}}{\min \big\{1,L_2 L_3/ \lambda_{min}(\widehat{\mathbf{A}}_{11}) \big\}} \, ,
\label{ThConNo:5}
\end{equation}
which the required upper bound. The required lower bound could be obtained by following the exact arguments in Proposition \ref{CondNoEx2} and \eqref{ThConNo:2}. Thus we get the desired result. \myqed

We mention that the notions and notations developed leading to Theorem \ref{ThCondNo} can also be extended to higher dimensions. An element will have $n_e$ vertices, e.g., $n_e$ could be $3$ or $4$ in 2-d. And the element stiffness matrices $A_{22}^{(k)}$ could be at most $n_e \times n_e$. The assembly argument \eqref{assembly} could be easily generalized to higher dimensions. For a given vertex $x_i$ and an enriched element $\tau_k$ in the FE star associated with $x_i$, the index $l(i,k)$ will again represent the local index of the vertex $x_{l(i,k)}^{(k)}$ of $\tau_k$ that coincides with $x_i$, i.e., $x_{l(i,k)}^{(k)}= x_i$. The expressions for $\Delta_i$, $(\mathbf{A}_{22})_{ii}$ and the Assumpsions \ref{asp1}, \ref{asp2}, \ref{asp3}, are exactly same in higher dimensions. Using these notions, the proof of Theorem \ref{ThCondNo} can be easily extended to higher dimensions. The approach presented here can also be extended for elasticity equations etc. We note however, the notations become a little more involved if $n_i>1$.

\begin{myremark} \label{AssumpComment} \upshape
We now make comments on the assumptions. The Assumption \ref{asp1} is always satisfied in 1-d. Let $B_0(u,v):= \int_\Omega u^\prime\, v^\prime \, dx$. Since $\locphibar{j}{i}(x_k)=0$ for $k=i-1,i,i+1$, it can be easily shown that $B_0(v_1,v_2)=0$ for all $v_1 \in \mcal{S}_1$ and $v_2 \in \overline{\mcal{S}}_2$. Therefore,
\begin{eqnarray*}
&&B(v_1+v_2,v_1+v_2) \ge \alpha B_0(v_1+v_2,v_1+v_2) \\
&&= \alpha [B_0(v_1,v_1)+ B_0(v_2,v_2)] \ge \frac{\alpha}{\beta} [ \|v_1\|_{\mcal{E}(\Omega)}^2 + \|v_2\|_{\mcal{E}(\Omega)}^2].
\end{eqnarray*}
Similarly, we can show that
\[
B(v_1+v_2,v_1+v_2) \le \frac{\beta}{\alpha} [ \|v_1\|_{\mcal{E}(\Omega)}^2 + \|v_2\|_{\mcal{E}(\Omega)}^2],
\]
and thus Assumption \ref{asp1} is satisfied with $L_1 = \frac{\alpha}{\beta}$ and $L_2 = \frac{\beta}{\alpha}$. In higher dimensions, this assumption has to be checked.

Assumption \ref{asp2} is equivalent to $L_2\|\mathbf{y}\|^2 \le \mathbf{y}^T \hat{A}_{22}^{(k)} \mathbf{y} \le U_2\|\mathbf{y}\|^2$ for all $\mathbf{y} \in \mathbb R^{\zeta^{(k)}}$. Thus $\hat{A}_{22}^{(k)}$ is uniformly positive definite in $k$ and its eigenvalues are uniformly bounded.

It is always possible to choose the diagonal matrix $D^{(k)}$ such that Assumption \ref{asp3} is satisfied. For example, it is easy to check that Assumption \ref{asp3} is satisfied with $L_3=U_3=1$ by choosing $D^{(k)} = \mbox{diag}\{\delta_1^{k},\delta_2^{(k)}\}$ with $\delta_2^{(k)} = (b_{jj}^{(k)})^{-1/2}$. The Assumption \ref{asp3} is trivially satisfied with $L_3=U_3=1$ when $D^{(k)}$ is a $1 \times 1$ matrix.
\myend
\end{myremark}

\begin{myremark} \upshape
As shown in the Appendix, the implementation of the SGFEM does not require scaling the stiffness matrix, i.e., the linear system involving the stiffness matrix $\mathbf{A}$, and not scaled version $\widehat{\mathbf{A}}$, is solved. The scaling was used only to define $\mathfrak{K}(\mathbf{A})$ and to study its order through Theorem \ref{ThCondNo}. We will show in the Appendix that $\mathfrak{K}(\mathbf{A})$ is an indicator of the loss of accuracy in the computed solution of the linear system associated with FEM, GFEM, and SGFEM. \myend
\end{myremark}

%% file: sec5p1_appl.tex
\section{Applications:} \label{APPLIC}

In this section we will present the SGFEM, when applied to three specific applications. We will primarily address in detail the scaled condition number of the  stiffness matrix of the method and show that the assumptions presented in the last section hold. The SGFEM, applied to each of these applications, will based on the uniform mesh $\{\tau_k\}_{k \in I\backslash \{0\}}$ with the set of vertices $\mcal{T}$, defined before.

\subsection{Interface Problems} \label{Interface} Let $a(x)$ in \eqref{VarProb} be a piecewise constant function and let $f$ be smooth. We will consider two situations, namely, $a(x)=a_1(x)$ and $a(x)=a_2(x)$, where
\[
 a_1(x) = \left\{ \begin{array}{ll}
                        \frac{1}{2}, & 0 \le x < b^*\\
                        1, & b^* \le x \le 1 \end{array} \right.
\ \
\mbox{and }
\ \
 a_2(x) = \left\{ \begin{array}{ll}
                        1, & 0 \le x < b_1^*\\
                        \frac{1}{2}, & b_1^* \le x < b_2^*\\
                        1, & b_2^* \le x \le 1 \end{array} \right.
\]
We note that the solution $u$ of \eqref{VarProb} does not belong to $H^2(\Omega)$.

We first consider $a(x)=a_1(x)$. We consider the set $\mcal{T}_2 \subset \mcal{T}$ as before. There exists an $m$ such that $b^* \in \open{\tau}_{m+1}=(x_{m},x_{m+1})$ and therefore, $b^* \in \omega_{m} \cap \omega_{m+1}$. For $x_i \in \mcal{T}$, we consider $V_i = \mbox{span}\{1, \locphi{1}{i}=\int_{x_{i-1}}^x (1/a_1(t))dt \}$.
Clearly, for $i\ne m,m+1$, we have $V_i = \mbox{span}\{1,(x-x_{i-1})\}$. Therefore recalling that $\overline{V}_i = \mbox{span}\{ \locphibar{1}{i} \}$,
where $\locphibar{1}{i} = \locphi{1}{i} - \intp{\omega_i} \locphi{1}{i}$, we get
$\overline{V}_i=\{0\}$ for $i \ne m,m+1$. We set $\mcal{T}_2=\{x_m,x_{m+1}\} \subset \mcal{T}$ and from the definition of $\overline{\mcal{S}}_2$, we have
\[
\overline{S}_2 = \sum_{x_i \in \mcal{T}_2} N_i \overline{V}_i = N_{m}\overline{V}_{m} + N_{m+1} \overline{V}_{m+1}\,.
\]
We further note that $\locphi{1}{m}$ is linear on $\tau_{m}$ and therefore, $\locphibar{1}{m}|_{\tau_{m}}=0$. Similarly, $\locphibar{1}{m+1}|_{\tau_{m+2}}=0$. Therefore
$\tau_{m+1}$ is the only enriched element, i.e., $\mcal{K}_{enr} = \{\tau_{m+1}\}$, and $\mathbf{A}_{22}=A_{22}^{(m+1)}$. Also, we can easily show that $\locphibar{1}{m}|_{\tau_{m+1}} = \locphibar{1}{m+1}|_{\tau_{m+1}}$. Let $b^* = x_m + \beta h$ with $0 < \beta < 1$. Then from a direct computation, we have
\begin{equation}
A_{22}^{(m+1)} = \left[ \begin{array}{cc}
                h \beta (1-\beta)^2(\frac{3}{2}+\beta-2\beta^2)/3 & h \beta^2 (1-\beta)^2(1+4\beta)/6 \\
                h \beta^2 (1-\beta)^2(1+4\beta)/6 & h \beta^2(1-\beta)(1+2\beta^2)/3
                \end{array} \right]. \label{eq:5.0}
\end{equation}
Clearly, $A_{22}^{(m+1)}$ is associated with the vertices $x_m$, $x_{m+1}$.
We choose the diagonal matrix $D^{(m+1)} = \mbox{diag}\{\delta_1^{(m+1)}, \delta_2^{(m+1)}\}$, where
\begin{equation}
\delta_1^{(m+1)}=h^{-1/2}\beta^{-1/2}(1-\beta)^{-1},\quad \delta_2^{(m+1)} = h^{-1/2}\beta^{-1}(1-\beta)^{-1/2}. \label{eq:5.0a}
\end{equation}
Then
\begin{eqnarray}
&&\hspace{-1cm}\hat{A}_{22}^{(m+1)} = D^{(m+1)} A_{22}^{(m+1)} D^{(m+1)} \nonumber \\[1ex]
&&=\left[ \begin{array}{cc}
                (\frac{3}{2}+\beta-2\beta^2)/3 & \beta^{1/2} (1-\beta)^{1/2}(1+4\beta)/6 \\
                \beta^{1/2} (1-\beta)^{1/2}(1+4\beta)/6 & (1+2\beta^2)/3
                \end{array} \right]. \label{eq:5.1}
\end{eqnarray}
The diagonal elements of $\hat{A}_{22}^{(m+1)}$ are $O(1)$ for all $0 < \beta < 1$. Also the eigenvalues of $\hat{A}_{22}^{(m+1)}$ are $\lambda_1 = (2-\beta)/6$ and $\lambda_2 = (1+\beta)/2$. Therefore, recalling Remark \ref{AssumpComment}, we have
\begin{equation}
\frac{1}{6}\|[D^{(m+1)}]^{-1}\mathbf{x}\|^2 \le \mathbf{x}^T A_{22}^{(m+1)} \mathbf{x} \le \|[D^{(m+1)}]^{-1}\mathbf{x}\|^2, \quad \forall \ \mathbf{x} \in \mathbf{R}^2, 
\end{equation}
and hence, Assumption \ref{asp2} is satisfied with $L_2 = \frac{1}{6}$ and $U_2 = 1$.

We set $\mathbf{D}_2 = \mbox{diag}\{d_1,d_2\}$ with $d_i = (\mathbf{A}_{22})_{ii}^{-1/2}$. Clearly, the diagonal elements of $\widehat{\mathbf{A}}_{22}=\mathbf{D}_2 \mathbf{A}_{22}\mathbf{D}_2$ are equal to 1. Recall that $\mcal{T}_2 = \{x_m,x_{m+1}\}$ and $\mcal{K}_m = \mcal{K}_{m+1} = \{\tau_{m+1}\}$. Therefore, $l(m,m+1) = 1$ and $l(m+1,m+1)=2$, where the index $l(i,k)$ for $x_i \in \mcal{T}_2$ and $\tau_k \in \mcal{K}_i$ was defined just before \eqref{Ai}. We also have $\mcal{K}_m^* = \mcal{K}_{m+1}^* = \{\tau_{m+1}\}$. Therefore from \eqref{Ai}, we have
$
\Delta_m = [\delta_1^{(m+1)}]^{-2} \mbox{ and }\ \Delta_{m+1} =  [\delta_2^{(m+1)}]^{-2}
$. Also the vertices $x_m,x_{m+1} \in \mcal{T}_2$ are associated with the diagonal elements $(\mathbf{A}_{22})_{j_mj_m},(\mathbf{A}_{22})_{j_{m+1}j_{m+1}}$, respectively, of $\mathbf{A}_{22}$, where $j_m=1,j_{m+1}=2$. It is easy to check that
\[
1 < (\mathbf{A}_{22})_{11}^{-1}\Delta_m,\ (\mathbf{A}_{22})_{22}^{-1}\Delta_{m+1} \le 6
\]
and the Assumption \ref{asp3} is satisfied with $L_3 = 1$ and $U_3=6$.

We have shown in Remark \ref{AssumpComment} that the Assumption \ref{asp1} is always satisfied in 1-d; in this case $L_1=\frac{1}{2}$ and $U_1=2$. Therefore, from Theorem \ref{ThCondNo}, we have that $\mathfrak{K}(\mathbf{A}) = \mcal{O}(h^{-2})$, and thus the GFEM with $\mcal{S}=\mcal{S}_1 + \overline{\mcal{S}}_2$ is indeed an SGFEM. We further note that Assumptions \ref{asp1}, \ref{asp2}, \ref{asp3} are satisfied for any $0< \beta < 1$, i.e., the constants $L_1, U_1, L_2, U_2, L_3$ and $U_3$ are independent of $\beta$. Therefore $\mathfrak{K}(\mathbf{A}) = \mcal{O}(h^{-2})$ even when $\beta \approx 0$ or $\beta \approx 1$, i.e., when the point of discontinuity $b^*$ of $a_1(x)$ is close to the one of the vertices $x_i$ (see also Remark \ref{Rem5.1beta}).


We next consider the \eqref{VarProb} with $a(x) = a_2(x)$. We again choose $V_i = \mbox{span}\{1,\locphi{1}{i} = \int_{x_{i-1}}^x (1/a_2(t)) dt \}$. If the points of discontinuity $b_1^*,\, b_2^*$ of $a_2(x)$ are separated, e.g.,  $b_1^* \in \open{\tau}_l$ and $b_2^* \in \open{\tau}_{l^*}$ with $|l-l^*| \ge 2$, then we can again show that the GFEM with $\mcal{S}=\mcal{S}_1 + \overline{\mcal{S}}_2$ is an SGFEM, based on the arguments given above.

Suppose there is an $m$ such that $b_1^* \in \open{\tau}_{m}$ and $b_2^* \in \open{\tau}_{m+1}$. Moreover, suppose $b_1^* = x_{m-1}+h/2$ and $b_2^* = x_m + \beta h$ with $0 < \beta < 1$. Note that $b_1^*$ is away from the vertices, whereas, $b_2^*$ could be close to either $x_m$ ($\beta \approx 0$) or $x_{m+1}$ ($\beta\approx 1$). As before, let $\overline{V}_i = \mbox{span}\{\locphibar{1}{i}\}$; clearly, $\overline{V}_i = \{0\}$ for $i \ne m-1,m,m+1$. Therefore $\mcal{T}_2=\{x_{m-1},x_m,x_{m+1}\}$ and
\[
\overline{S}_2 = \sum_{i=m-1,m,m+1} N_i\, \overline{V}_i\, .
\]
We further note that $\locphibar{1}{m-1} \big|_{\tau_{m-1}} = \locphibar{1}{m+1} \big|_{\tau_{m+2}} = 0$. Also it can be shown that $\locphibar{1}{m-1} \big|_{\tau_{m}} = \locphibar{1}{m} \big|_{\tau_{m}}$ and $\locphibar{1}{m} \big|_{\tau_{m+1}} = \locphibar{1}{m+1} \big|_{\tau_{m+1}}$. Therefore $\mcal{K}_{enr}=\{\tau_m,\tau_{m+1}\}$ (i.e., $\tau_m,\tau_{m+1}$ are the only enriched elements), and hence $\mathbf{A}_{22}$ is assembled from local stiffness matrices $A_{22}^{(m)}$ and $A_{22}^{(m+1)}$.

From direct computation, we get
\[
A_{22}^{(m)} = \left[ \begin{array}{cc}
                 \frac{h}{16} & \frac{h}{32} \\
                \frac{h}{32} & \frac{h}{32}
                \end{array} \right],
\]
and it is associated with the vertices $x_{m-1}$ and $x_m$. The matrix $A_{22}^{(m+1)}$ is same as in \eqref{eq:5.0} and is associated with $x_{m}$ and $x_{m+1}$. We choose $D^{(m)} = \mbox{diag}(\delta_1^{(m)},\ \delta_2^{(m)})$ with
$\delta_1^{(m)}= \delta_2^{(m)}=h^{-1/2}$ and $D^{(m+1)} = \mbox{diag}(\delta_1^{(m+1)},\ \delta_2^{(m+2)})$ with $\delta_1^{(m+1)},\, \delta_2^{(m+2)}$ as given in \eqref{eq:5.0a}.
Then
\[
\hat{A}_{22}^{(m)} := D^{(m)} A_{22}^{(m)} D^{(m)} =
\left[ \begin{array}{cc}
                 \frac{1}{16} & \frac{1}{32} \\
                \frac{1}{32} & \frac{1}{16}
                \end{array} \right].
\]
Clearly the diagonal elements of $\hat{A}_{22}^{(m)}$ are $O(1)$. The eigenvalues of $\hat{A}_{22}^{(m)}$ are $\lambda_1 = 1/32$ and  $\lambda_2=3/32$ and therefore (recall Remark \ref{AssumpComment}),
\begin{equation}
\frac{1}{32}\|[D^{(m)}]^{-1} \mathbf{x} \|^2 \le \mathbf{x}^T A_{22}^{(m)} \mathbf{x} \le \frac{3}{32}\|[D^{(m)}]^{-1} \mathbf{x} \|^2. \label{eq:5.3}
\end{equation}

Next, the matrix $\hat{A}_{22}^{(m+1)}:= D^{(m+1)} A_{22}^{(m+1)} D^{(m+1)}$ is same as the matrix given in \eqref{eq:5.1}. The diagonal elements of $\hat{A}_{22}^{(m+1)}$ are $O(1)$ and its eigenvalues are $\lambda_1=(2-\beta)/6$ and $\lambda_2=(1+\beta)/2$. Therefore
\begin{equation*}
\frac{1}{6}\|[D^{(m+1)}]^{-1}\mathbf{x}\|^2 \le \mathbf{x}^T A_{22}^{(m+1)} \mathbf{x} \le \|[D^{(m+1)}]^{-1}\mathbf{x}\|^2, \quad \forall \ \mathbf{x} \in \mathbf{R}^2.
\end{equation*}
Thus the above inequality together with \eqref{eq:5.3} implies that Assumption \ref{asp2} is satisfied with $L_2= \frac{1}{32}$ and $U_2=1$ for all $0 < \beta < 1$.

The matrix $\mathbf{A}_{22}$ is assembled from the matrices $A_{22}^{(m)}$, $A_{22}^{(m+1)}$ and is given by
\[
\mathbf{A}_{22} = \left[ \begin{array}{ccc}
    \frac{h}{16} & \frac{h}{32} & 0 \\[1.0ex]
    \frac{h}{16} & \frac{h}{32} + \frac{h \beta (1-\beta)^2}{3}\,(\frac{3}{2}+\beta-2\beta^2) & \frac{h \beta^2 (1-\beta)^2}{3}\,(\frac{1}{2}+2\beta) \\[1.0ex]
    0 & \frac{h \beta^2 (1-\beta)^2}{3}\,(\frac{1}{2}+2\beta)& \frac{h\beta^2 (1-\beta)}{3}\,(1+2\beta^2)
                \end{array} \right].
\]
We choose
$
\mathbf{D}_2 = \mbox{diag}(d_1,\ d_2,\ d_3) \mbox{ with } d_i = (\mathbf{A}_{22})_{ii}^{-1/2}.
$
Clearly the diagonal elements of $\widehat{\mathbf{A}}_{22}:= \mathbf{D}_2 \mathbf{A}_{22} \mathbf{D}_2$ are equal to 1. Consider the vertex $x_m \in \mcal{T}_2$. Then $\mcal{K}_m = \{\tau_m,\tau_{m+1}\}$ and $l(m,m) =2$, $l(m,m+1)=1$. Also in this case, $\mcal{K}_i^* = \mcal{K}_i$. Therefore, from \eqref{Ai}, we have $\Delta_m = [\delta_2^{(m)}]^{-2} + [\delta_1^{(m+1)}]^{-2}$. Similarly, we can show that $\Delta_{m-1} = [\delta_1^{(m)}]^{-2}$ and $\Delta_{m+1} = [\delta_2^{(m+1)}]^{-2}$. We also note that the vertices $x_{m-1},x_m,x_{m+1} \in \mcal{T}_2$ are associated with the diagonal elements $(\mathbf{A}_{22})_{j_{m-1}j_{m-1}}, (\mathbf{A}_{22})_{j_{m}j_{m}},(\mathbf{A}_{22})_{j_{m+1}j_{m+1}}$, respectively, of $\mathbf{A}_{22}$, where $j_{m-1}=1,j_m=2,j_{m+1}=3$. An easy calculation yields
\[
1 \le (\mathbf{A}_{22})_{11}^{-1}\Delta_{m-1},\ (\mathbf{A}_{22})_{22}^{-1} \Delta_m,\ (\mathbf{A}_{22})_{33}^{-1}\Delta_{m+1} \le 16.
\]
Thus Assumption \ref{asp3} is satisfied with $L_3 = 1$ and $U_3 = 16$ for all $0 < \beta < 1$. We have shown before that Assumption \ref{asp1} is always satisfied in 1-d. Therefore from Theorem \ref{ThCondNo}, we infer that $\mathfrak{K}(\mathbf{A}) = \mcal{O}(h^{-2})$; the result is true even when $\beta \approx 0$ or $\beta \approx 1$. Thus the GFEM with $\mcal{S}=\mcal{S}_1+\overline{\mcal{S}}_2$ is indeed an SGFEM.

We remark that for $a(x)=a_1(x)$ or $a(x)=a_2(x)$, we can show that there exists $\bar{\xi}^i \in \overline{V}_i$ such that $\|u - \intp{\omega_i}u - \bar{\xi}^i \|_{\mcal{E}(\omega_i)} = O(h)$ for each $x_i \in \mcal{T}_2$. Thus using the standard interpolation estimates and using Theorem \ref{ApproxThmSGFEM1d}, we have $\|u - u_h \|_{\mcal{E}(\Omega)} = O(h)$, where $u_h$ is the SGFEM solution.

\begin{myremark} \label{Rem5.1beta} \upshape Note that $A_{22}^{(m+1)}$
and thus $\mathbf{A}_{22}$, $\mathbf{A}$ degenerate as $\beta \to 0$ or $\beta \to 1$. Let $\epsilon_0$ be small, say, $\epsilon_0 = 10^{-14}$. We adjust the implementation when $\beta \le \epsilon_0$ or $1-\beta \le \epsilon_0$ by setting $\beta = \epsilon_0$ or $1-\beta=\epsilon_0$, respectively. We emphasize that $\mathfrak{K}(\mathbf{A})$ is bounded
independently of $\beta$. \myend
\end{myremark} 

%% file: sec5p2_appl.tex
\subsection{Problems with singular solutions}

Let $a(x)=1$ in \eqref{VarProb} and suppose $f(x)$ be such that the solution $u$ of \eqref{VarProb}-\eqref{Buv} is of the form $u = x^\alpha + u_0$, where $\frac{1}{2} < \alpha < \frac{3}{2}$, $\alpha \ne 1$, and $u_0$ is smooth with $u_0(0)=0$. Clearly $u \notin H^2(\Omega)$. Let $0 < D < 1$ and set $\Omega_l:=(0,D)$, $\Omega_r:=(D,1)$. Then $u \in H^2(\Omega_r)$ and $|u|_{H^2(\Omega_r)} \le C[\, |x^\alpha|_{H^2(\Omega_r)} + |u_0|_{H^2(\Omega_r)}]$. Clearly, $|u|_{H^2(\Omega_r)}$ depends on $D$ and is extremely large for $D\approx 0$.

We consider $\mcal{T}_1 \subset \mcal{T}$ as before. Let $\mcal{T}_2 := \{x_i \in \mcal{T}:\, \omega_i \cap \Omega_l \ne \emptyset \}$, where the patches $\omega_i$ have been defined before. Clearly, $x_0,x_1 \in \mcal{T}_2$. Let $k^* \in I$ be the largest index such that $x_i \in \mcal{T}_2$ for $0\le i \le k^*-1$. For $x_i \in \mcal{T}_2$, let
\[
V_i = \mbox{span}\{1, \locphi{1}{i}=(x-x_i), \locphi{2}{i} = x^\alpha|_{\omega_i} \},
\]
and for $x_i \in \mcal{T}\backslash \mcal{T}_2$, let $V_i = \mbox{span}\{1, \locphi{1}{i}=(x-x_i) \}$. Clearly, $\overline{V}_i = \{0\}$ for $x_i \in \mcal{T}\backslash \mcal{T}_2$. For $x_i \in \mcal{T}_2$, we have
\[
\overline{V}_i = \mbox{span}\{\locphibar{2}{i} = \sigma_i \} \ne 0,
\]
where $\sigma_i := (x^\alpha - \intp{\omega_i}x^\alpha)|_{\omega_i}$; recall $\intp{\omega_i}x^\alpha$ is the piecewise linear polynomial that interpolates $x^\alpha$ at the vertices $\{x_{i-1},x_i,x_{i+1}\}$ of $\omega_i$ for $i\ne 0$, and  $\intp{\omega_0}x^\alpha$ interpolates $x^\alpha$ at $\{x_0,x_1\}$. For an element $\tau_k \subset \overline{\omega}_i$ (with $x_i \in \mcal{T}_2$), we define $\sigma^{(k)}:= (x^\alpha - I_kx^\alpha)|_{\tau_k}$, where $I_kx^\alpha \in \mcal{P}^1(\tau_k)$  interpolates $x^\alpha$ at $x_{k-1}$, $x_k$. Clearly, $\intp{\omega_i}x^\alpha  = I_kx^\alpha$ on $\tau_k \subset \overline{\omega}_i$. It is also clear that $[\sigma^{(k)}]^\prime \not \equiv 0$ on $\tau_k\subset \overline{\omega}_i$.

We define $\overline{\mcal{S}}_2 = \sum_{i=0}^{k^*-1} N_i \overline{V}_i$ and we consider the GFEM based $\mcal{S} = \mcal{S}_1 + \overline{S}_2$. We first address the convergence of the GFEM solution $u_h$. It is easy to show that for $x_i \in \mcal{T}_2$, there exists $\bar{\xi}_i \in V_i$ such that
\[
\|u - \intp{\omega_i}u - \bar{\xi}_i\|_{\mcal{E}(\omega_i)} \le Ch|u_0|_{H^2(\omega_i)}\, .
\]
Also for $x_i \in \mcal{T}\backslash \mcal{T}_2$, from standard interpolation result we have
\[
\|u - \intp{\omega_i} u\|_{\mcal{E}(\omega_i)} \le Ch|u|_{H^2(\omega_i)} \le Ch[|x^\alpha|_{H^2(\omega_i)} + |u_0|_{H^2(\omega_i)}].
\]
Therefore, from Theorem \ref{ApproxThmSGFEM1d}, there exists $v \in \mcal{S}_1 + \overline{\mcal{S}}_2$ such that
\begin{eqnarray*}
&& \hspace{-.8cm} \|u - v\|_{\mcal{E}(\Omega)} \le C h \big[ \sum_{x_i \in \mcal{T}_2} |u_0|_{H^2(\omega_i)}^2 + \sum_{x_i \in \mcal{T}\backslash \mcal{T}_2} \{ |x^\alpha|_{H^2(\omega_i)}^2 + |u_0|_{H^2(\omega_i)}^2\} \big]^{1/2}\\
&& \hspace{1.2cm} \le Ch \big[|u_0|_{H^2(\Omega)}^2 + |x^\alpha|_{H^2(\Omega_r)}^2 \big]^{1/2}.
\end{eqnarray*}
Thus we have $\|u - u_h\|_{\mcal{E}(\Omega)} \le Ch$, where $u_h$ is the GFEM solution; note that $C$ depends on $|x^\alpha|_{H^2(\Omega_r)}$ and thus on $D$.

We note that $\Omega_l$ is independent of $h$. However, if $D=h^\gamma$, $\gamma < 1$ (i.e., $|\Omega_l|=h^\gamma$), then one can show that $\|u - u_h\|_{\mcal{E}(\Omega)} =O(h^{1-\gamma})$. Thus if we enrich only a fixed number of patches in the neighborhood of the singularity, we loose the optimal order of convergence.

We now address the scaled condition number of the stiffness matrix $\mathbf{A}$ of the GFEM. We note that the matrix $\mathbf{A}_{22}$ is assembled from element stiffness matrices $A_{22}^{(k)}$ for the element $\tau_k$, where $\tau_k$ is enriched. We note that the set of enriched elements is given by $\mcal{K}_{enr}:=\{\tau_k \in \{\tau_l\}_{l \in I \backslash \{0\}}:\, x_k \in \mcal{T}_2\}$. We further note that if $\tau_k \in \mcal{K}_{enr}$, then $\tau_j \in \mcal{K}_{enr}$ for $1 \le j \le k$. Also from the definition of $k^*$,  it is clear that $\tau_{k^*} \in \mcal{K}_{enr}$ and $\tau_j \notin \mcal{K}_{enr}$ for $j\ge k^*+1$. Now, for $\tau_k \in \mcal{K}_{enr}$, $k \ne k^*$, the matrices $A_{22}^{(k)}$ are of the form $A_{22}^{(k)} = \{b_{lm}^{(k)}\}_{l,m=1}^2$; the entries $b_{lm}^{(k)}$ are as given by
\[
b_{lm}^{(k)} = \int_{\tau_k} (N_{k-2+l} \sigma)^\prime (N_{k-2+m}\sigma)^\prime \, dx\, .
\]
Also, since $x_{k^*} \notin \mcal{T}_2$, we have $A_{22}^{(k^*)}=[b_{11}^{(k^*)}]$ (an $1\times 1$ matrix), where $b_{11}^{(k^*)}$ is given by the above expression.

\begin{mylemma} \label{EntriesLemm} The entries of the matrix $A_{22}^{(k)}$ are as follows:
\begin{eqnarray*}
&&b_{11}^{(k)} = \int_{\tau_k} N_{k-1}^2 {\sigma^\prime}^2\,dx, \quad
b_{22}^{(k)} = \int_{\tau_k} N_{k}^2 {\sigma^\prime}^2\,dx ,\\
&& b_{12}^{(k)} = b_{21}^{(k)} = \int_{\tau_k} N_{k-1} N_{k}  {\sigma^\prime}^2\,dx\, .
\end{eqnarray*} \end{mylemma}
The proof is easy and we do not present it here. \myqed

It is clear from above that for $\tau_k \in \mcal{K}_{enr}$ and $k\ne k^*$, the diagonal elements $b_{11}^{(k)},b_{22}^{(k)} >0$ and therefore $A_{22}^{(k)}$ is associated with $x_{k-1}$, $x_k$. Also $b_{11}^{(k^*)}>0$ and thus $A_{22}^{(k^*)}$ is associated with $x_{k^*-1}$. A simple observation yields that the size of $\mathbf{A}_{22}$ is $k^* \times k^*$.

Let $\tau_k \in \mcal{K}_{enr}$ and set $x_{k-1/2} := (k-\frac{1}{2})h$; $x_{k-1/2}$ is the mid-point of $\tau_k$. We define
\[
G_k = \big|[x^\alpha]^{\prime \prime}(x_{k-1/2})\big| = |\alpha (\alpha-1) (k-\textstyle\frac{1}{2})^{\alpha-2}h^{\alpha-2}|.
\]
Note that for $1 \le j \le k^*-1$, $\tau_{j+1} \in \mcal{K}_{enr}$ implies $\tau_j \in \mcal{K}_{enr}$, and we have
\begin{equation}
1 \le \frac{G_{j}}{G_{j+1}} = \Big(\frac{j+\frac{1}{2}}{j-\frac{1}{2}}\Big)^{2-\alpha} \le 3^{2-\alpha} . \label{BndOnRatio}
\end{equation}


We now obtain a few results, which will be used to establish that the GFEM with $\mcal{S}=\mcal{S}_1+\overline{\mcal{S}}_2$ is an SGFEM.

\begin{mylemma} \label{BndsOnNormSigmaPrime} For $x_k \in \mcal{K}_{enr}$, there exist positive constants $C_1^*, C_2^*$, independent of $k$ and $h$ but may depend on $\alpha$, such that
\[
C_1^*h^{3/2} \le \frac{\|[\sigma^{(k)}]^\prime\|_{L^2(\tau_k)}}{G_k} \le C_2^*h^{3/2}.
\] \end{mylemma}

\medskip

\noindent \emph{Proof:} (a) First let $2 \le k \le k^*$ and let $g(x) = [\sigma^{(k)}]^\prime (x)$ for $x\in \tau_k$. Then
\begin{eqnarray}
\max |g^\prime(x)| &=& \max_{x \in \tau_k} |[\sigma^{(k)}]^{\prime \prime}(x)| = |\alpha (\alpha -1)|x_{k-1}^{\alpha -2} \nonumber \\
&=& |\alpha(\alpha-1)|(k-{\textstyle\frac{1}{2}})^{\alpha-2}h^{\alpha-2} \Big( \frac{k-1}{k-\frac{1}{2}} \Big)^{\alpha-2} \nonumber \\
&=& G_k \Big( \frac{k-\frac{1}{2}}{k-1} \Big)^{2-\alpha} \le G_k \Big(\frac{3}{2}\Big)^{2-\alpha} := M_k .\label{gPrimeMax}
\end{eqnarray}
Similarly,
\begin{equation}
\min |g^\prime(x)| = G_k \Big( \frac{k-\frac{1}{2}}{k} \Big)^{2-\alpha}
\ge G_k \Big(\frac{3}{4}\Big)^{2-\alpha} := m_k .\label{gPrimeMin}
\end{equation}

We next note that $g^\prime$ does not change sign in $\tau_k$ and thus $g$ is monotonic in $\tau_k$. Also, since
$
\int_{\tau_k} g\, dx = \sigma^{(k)}\big|_{x_{k-1}}^{x_{k}} = 0,
$
there exists a unique $x_k^* = x_{k-1} + \zeta h \in \tau_k$ with $0 < \zeta < 1$ such that $g(x_k^*)=0$ and $x_k^*$ is characterized by
\begin{equation}
\int_{x_{k-1}}^{x_k^*} |g|\, dx = \int_{x_k^*}^{x_{k}} |g|\, dx \, . \label{IntegralsEqual}
\end{equation}
We now obtain bounds on $\zeta$, independent of $k$. Since $g(x_k^*)=0$, it is clear from the mean value theorem, \eqref{gPrimeMax}, and \eqref{gPrimeMin} that
\begin{eqnarray}
m_k |x-x_{k}^*| &\le& \min |g^\prime| |x-x_{k}^*| \le |g(x)| \nonumber \\
 &\le& \max |g^\prime| |x-x_{k}^*| \le M_k |x-x_{k}^*|, \quad x \in \tau_k \, . \label{BoundOng}
\end{eqnarray}
Consequently,
\begin{equation}
m_k\, \frac{\zeta^2h^2}{2} \le \int_{x_{k-1}}^{x_k^*} |g|\, dx \le M_k \frac{\zeta^2h^2}{2} \label{LeftIntegralBound}
\end{equation}
and
\begin{equation}
m_k\, \frac{(1-\zeta)^2h^2}{2} \le \int_{x_k^*}^{x_{k}} |g|\, dx \le M_k \frac{(1-\zeta)^2h^2}{2}\, . \label{RightIntegralBound}
\end{equation}
Now from \eqref{IntegralsEqual}, \eqref{LeftIntegralBound}, and \eqref{RightIntegralBound}, we have
\[
m_k\, \frac{\zeta^2h^2}{2} \le \int_{x_{k-1}}^{x_k^*} |g|\, dx = \int_{x_k^*}^{x_{k}} |g|\, dx \le M_k \frac{(1-\zeta)^2h^2}{2}\, ,
\]
and thus
\[
\zeta \le \frac{\sqrt{M_k}}{\sqrt{M_k} + \sqrt{m_k}} = \frac{(3/2)^{(2-\alpha)/2}}{(3/2)^{(2-\alpha)/2} + (3/4)^{(2-\alpha)/2}}\, ,
\]
where we used the definition of $m_k$ and $M_k$ given in \eqref{gPrimeMin} and \eqref{gPrimeMax} respectively.

Using a similar argument we obtain a lower bound of $\zeta$; we summarize the bounds of $\zeta$ as
\begin{eqnarray}
&& \zeta_l \le \zeta \le \zeta_r, \quad \mbox{where} \nonumber \\
&&\zeta_l = \frac{(3/4)^{(2-\alpha)/2}}{(3/2)^{(2-\alpha)/2} + (3/4)^{(2-\alpha)/2}} \nonumber \\
\mbox{ and }  && \zeta_r = \frac{(3/2)^{(2-\alpha)/2}}{(3/2)^{(2-\alpha)/2} + (3/4)^{(2-\alpha)/2}} \, .\label{EstimateOnZeta}
\end{eqnarray}

Finally, from \eqref{BoundOng} and using the definition of $M_k$ (given in \eqref{gPrimeMax}), we get
\begin{eqnarray*}
\int_{\tau_k} |g|^2\, dx &=& \int_{x_{k-1}}^{x_k^*} |g|^2\, dx + \int_{x_k^*}^{x_{k}} |g|^2\, dx  \\
&\le& M_k^2 \int_{x_{k-1}}^{x_k^*} (x_k^* -x)^2\, dx + M_k^2  \int_{x_k^*}^{x_{k}} (x-x_k^*)^2\, dx  \\
&=& \frac{M_k^2h^3}{3} [ \zeta^3 + (1-\zeta)^3] \\
&\le& \frac{G_k^2 (3/2)^{2(2-\alpha)}h^3}{3} [\zeta_r^3 + (1-\zeta_l)^3],
\end{eqnarray*}
and similarly, we have
\begin{eqnarray*}
\int_{\tau_k} |g|^2\, dx &\ge& \frac{m_k^2h^3}{3} [ \zeta^3 + (1-\zeta)^3] \\
&\ge& \frac{G_k^2 (3/4)^{2(2-\alpha)}h^3}{3} [\zeta_l^3 + (1-\zeta_r)^3].
\end{eqnarray*}
Thus
\[
\bar{C}_1^*h^{3/2} \le \frac{\|[\sigma^{(k)}]^\prime\|_{L^2(\tau_k)}}{G_k} \le \bar{C}_2^* h^{3/2},\quad \mbox{for } 2 \le i \le N,
\]
where $\bar{C}_1^* = (3/4)^{2-\alpha}\sqrt{[\zeta_l^3 + (1-\zeta_r)^3]/3}$ and
$\bar{C}_2^* = (3/2)^{2-\alpha} \sqrt{[\zeta_r^3 + (1-\zeta_l)^3]/3}$.

\medskip

(b) We now consider $k=1$. We note that on $\tau_1=(0,h)$,
\[
[\sigma^{(1)}]^\prime(x) = \alpha x^{\alpha - 1} - h^{\alpha-1}.
\]
By a direct computation, we get
\[
\int_0^h |[\sigma^{(1)}]^\prime|^2\, dx = \frac{(\alpha-1)^2h^{2\alpha-1}}{2\alpha-1}\, .
\]
Therefore,
\[
\frac{\int_0^h |[\sigma^{(1)}]^\prime|^2\, dx}{G_1^2} = \frac{(\alpha-1)^2h^{2\alpha-1}}{(2\alpha-1)\alpha^2 (\alpha-1)^2h^{2\alpha - 4}2^{4-2\alpha} }
= \frac{h^3}{(2\alpha-1) \, \alpha^2\, 2^{4-2\alpha}}:=\bar{C}^*.
\]
Hence we get the desired result with $C_1^*=\min (\bar{C}_1^*, \bar{C}^*)$ and $C_2^*= \max (\bar{C}_2^*, \bar{C}^*)$. \myqed

\begin{mylemma} \label{LowBndOnSigpLp} Suppose $\tau_k \in \mcal{K}_{enr}$ and let $l_k(x)$ be a linear function, defined on $\tau_k$, such that $l_k(x_{k-1}) = y_{1}$ and $l_k(x_{k}) = y_{2}$. Then there exists a positive constant $C_3^*$, independent of $k$ and $h$ but may depend on $\alpha$, such that
\[
\|{[\sigma^{(k)}]^\prime} l_k\|_{L^2(\tau_k)} \ge C_3^* G_k h^{3/2} (y_{1}^2+y_2^2)^{1/2}.
\]
\end{mylemma}

\noindent \emph{Proof:} (a) Let $2 \le k \le k^*$ and define $g(x)= [\sigma^{(k)}]^\prime (x)$ for $x\in \tau_k$. On $\tau_k$, we have seen in the proof of Lemma \ref{BndsOnNormSigmaPrime} that $g(x_k^*)=0$ where $x_k^* = x_{k-1}+ \zeta h$ and $0 < \zeta_l \le \zeta \le \zeta_r <1$. We have also seen that
\begin{equation}
m_k|x-x_k^*| \le |g(x)| \le M_k|x-x_k^*|, \quad \forall \ x \in \tau_k \, , \label{BoundOngLem5}
\end{equation}
where
\[
m_k = G_k (3/4)^{2-\alpha}, \quad M_k = G_k (3/2)^{2-\alpha}.
\]

Let
\[ \bar{x}_k \equiv x_{k-1} + \zeta_rh + \frac{1-\zeta_r}{2}\, h.
\]
Then
\begin{equation}
|\bar{x}_k - x_k^*| = |x_{k-1} + \zeta_rh + \frac{1-\zeta_r}{2}\, h -x_{k-1} -\zeta h| \ge \frac{1-\zeta_r}{2}\, h .\label{smallerint}
\end{equation}
Also from the definition of $\bar{x}_k$, it is clear that $g(x) \ne 0$ in $(\bar{x}_k,x_k)$ and thus from \eqref{BoundOngLem5} and \eqref{smallerint}, we have
\begin{equation}
|g(x)| \ge m_k |\bar{x}_k - x_k^*|  \ge \frac{1-\zeta_r}{2}\, m_k h.
\end{equation}
Therefore,
\begin{equation}
\int_{x_{k-1}}^{x_k} |g|^2|l_k|^2\, dx \ge \int_{\bar{x}_{k}}^{x_k} |g|^2|l_k|^2\, dx \ge m_k^2 h^2 \frac{(1-\zeta_r)^2}{4} \int_{\bar{x}_{k}}^{x_k}|l_k|^2\, dx \, .\label{LowerBnd}
\end{equation}

We make the change of variable $ y = \frac{x-\bar{x}_k}{h} \frac{2}{1-\zeta_r}$. Then
\[
F(y_1,y_2):=\frac{\int_{\bar{x}_{k}}^{x_i}|l_k|^2\, dx}{y_1^2 + y_2^2}
= \frac{(1-\zeta_r)h \int_0^1|\tilde{l}(y)|^2\, dy}{2(y_1^2 + y_2^2)}\, ,
\]
where
\[
\tilde{l}(y) = l_k(\bar{x}_k +\frac{(1-\zeta_r)hy}{2}) = y_1\, \frac{1-\zeta_r}{2}(1-y) +
y_2\, (\frac{1+\zeta_r}{2}+\frac{1-\zeta_r}{2}y).
\]
Thus $F(y_1,y_2)$ is independent of $k$. We next note that $F(y_1,y_2)$ is a continuous function and $F(\beta y_1, \beta y_2) = F(y_1,y_2)$. It is well known that the minimum of $F(y_1,y_2)$ is attained on the compact set $y_1^2+y_2^2=1$. Hence there is a constant $C_{min}$, independent of $k$ but may depend on $\zeta_r$, such that
\[
0 < C_{min}^2 \frac{(1-\zeta_r)h}{2}\le F(y_1,y_2) = \frac{\int_{\bar{x}_{k}}^{x_k}|l|^2\, dx}{y_1^2 + y_2^2}\, .
\]
Thus from \eqref{LowerBnd}, we have
\begin{eqnarray*}
\int_{x_{k-1}}^{x_k} |g|^2|l_k|^2\, dx &\ge& C_{min}^2 m_k^2 h^3 \frac{(1-\zeta_r)^3}{8}(y_1^2+y_2^2) \\
&=&{B_1^*}^2 G_k^2 h^3(y_1^2 +y_2^2),
\end{eqnarray*}
where ${B_1^*}^2 = (3/4)^{2(2-\alpha)}C_{min}^2(1-\zeta_r)^3/8$.

(b) We now consider $k=1$. On $\tau_1=(0,h)$, we have $g(x) = \alpha x^{\alpha -1}-h^{\alpha-1}$. It is easy to see that $g(x_1^*)=0$, where $x_1^* = \zeta h$ with $\zeta =\zeta(\alpha)= \alpha ^{1/1-\alpha}$. Since $\zeta(\alpha)$ is increasing for $\frac{1}{2} < \alpha < \frac{3}{2}$ (with $\zeta$ redefined for $\alpha=1$), we have $\zeta \le \zeta^* \equiv \zeta(3/2)=(2/3)^2$.

Set $\bar{x}_1 = \zeta^* h + (1-\zeta^*)h/2$. Since $|g(x)|$ is increasing in $(x_1^*,h)$,  we have $|g(x)| \ge \bar{g}_{min} \equiv  |g(\bar{x}_1)|$ on $(\bar{x}_1,h)$. Therefore,
\begin{eqnarray*}
&&\int_0^h |g|^2 |l_1|^2\, dx > \int_{\bar{x}_1}^h |g|^2 |l_1|^2\, dx
\ge {\bar{g}_{min}}^2 \int_{\bar{x}_1}^h |l_1|^2\, dx \\
&&\hspace{1cm} =\frac{\bar{g}_{min}^2 G_1^2}{\alpha^2|\alpha-1|^2h^{2(\alpha-2)}/2^{2(\alpha-2)}}\,  \int_{\bar{x}_1}^h |l_1|^2\, dx
= C^2 G_1^2 h^2 \int_{\bar{x}_1}^h |l_1|^2\, dx\, ,
\end{eqnarray*}
where
\[
C^2 =\frac{ 2^{4-2\alpha}\bar{g}_{min}^2}{\alpha^2(\alpha-1)^2 h^{2(\alpha-1)}}
= \frac{ 2^{4-2\alpha}[\alpha \frac{(1+\zeta^*)}{2} -1]^2 }{ \alpha^2(\alpha-1)^2 }\, .
\]
As before, we can show that
\[
\int_{\bar{x}_1}^h |l_1|^2\, dx \ge C_{min} \frac{(1-\zeta^*)h}{2}(y_{1}^2 + y_2^2)^{1/2},
\]
and therefore,
\[
\int_0^h |g|^2 |l_1|^2\, dx \ge {B_2^*}^2 G_1^2 h^3 (y_{1}^2 + y_2^2)^{1/2},
\]
where ${B_2^*}^2 = C^2 C_{min}(1-\zeta^*)h/2$.
Finally, defining ${C_3^*} = \min (B_1^*, B_2^*)$ and recalling that $g=[\sigma^{(k)}]^\prime$, we get the desired result. \myqed

\medskip

Now, for $k\ne k^*$, consider the diagonal matrix $D^{(k)}=\mbox{diag}(\delta_1^{(k)},\delta_2^{(k)})$ with $\delta_1^{(k)}=\delta_2^{(k)} = G_k^{-1}h^{-3/2}$ and set $\hat{A}_{22}^{(k)} = D^{(k)} A_{22}^{(k)} D^{(k)}$. The diagonal elements of $\hat{A}_{22}^{(k)}$ (see Lemma \ref{EntriesLemm}) are
\[
\bar{b}_{11}^{(k)} = \frac{1}{G_k^2 h^3} \int_{\tau_k} N_{k-1}^2 {\sigma^\prime}^2\,dx\ , \quad \bar{b}_{22}^{(k)} = \frac{1}{G_k^2 h^3} \int_{\tau_k} N_{k}^2 {\sigma^\prime}^2\,dx \, .
\]
Using Lemmas \ref{LowBndOnSigpLp} and \ref{BndsOnNormSigmaPrime}, it is clear that
\[
C_3^* \le \bar{b}_{11}^{(k)} \le \frac{1}{G_k^2 h^3} \int_{\tau_k} {\sigma^\prime}^2\,dx \le C_2^*,
\]
where $C_2^*,\, C_3^*$ are independent of $k$ and $h$. Similarly,
\[
C_3^* \le \bar{b}_{22}^{(k)} \le \frac{1}{G_k^2 h^3} \int_{\tau_k} {\sigma^\prime}^2\,dx \le C_2^*.
\]
We let $D^{(k^*)} = [\delta_1^{(k^*)}]$ with $\delta_1^{(k^*)}=G_{k^*}^{-1}h^{-3/2}$.  Using similar arguments we show the $C_3^* \le \bar{b}_{11}^{(k^*)} \le C_2^*$, where $\hat{A}_{22}^{(k^*)} = D^{(k^*)} A_{22}^{(k^*)} D^{(k^*)} = [\bar{b}_{11}^{(k^*)}]$. Thus the diagonal elements of $\widehat{A}_{22}^{(k)}$ are $O(1)$ for all $\tau_k \in \mcal{K}_{enr}$.

\medskip

We next show that the element matrices $A_{22}^{(k)}$ satisfy the Assumption \ref{asp2}.

\begin{myprop} \label{SingA2} For $\tau_k \in \mcal{K}_{enr}$, the matrices $A_{22}^{(k)}$ satisfies Assumption \ref{asp2}.
\end{myprop}

\noindent \emph{Proof:} Suppose $k \ne k^*$ and let $\mathbf{x} = (x_1,x_2)^T \in \mathbb R^2$. Then, using Lemma \ref{EntriesLemm}, we have
\begin{eqnarray*}
\mathbf{x}^T A_{22}^{(k)} \mathbf{x} &=& b_{11}^{(k)} x_1^2 +2 b_{12}^{(k)} x_1 x_2 + b_{22}^{(k)} x_2^2\\
&&\hspace{-.5in}=  \int_{\tau_k} [x_1^2 N_{k-1}^2 + 2x_1 x_2 N_{k-1}N_{k} + x_2^2 N_k^2] {[\sigma^{(k)}]^\prime}^2\, dx\\
&&\hspace{-.5in}= \int_{\tau_k} [x_1 N_{k-1} + x_2 N_{k}]^2 {[\sigma^{(k)}]^\prime}^2\, dx \le 2 (x_1^2 + x_2^2)\int_{\tau_k} {[\sigma^{(k)}]^\prime}^2\, dx \, ,
\end{eqnarray*}
and using Lemma \ref{BndsOnNormSigmaPrime}, we have
\[
\mathbf{x}^T A_{22}^{(k)} \mathbf{x} \le C^*_2 h^3 G_k^2 \|\mathbf{x}\|^2 = C^*_2 \|[D^{(k)}]^{-1}\mathbf{x}\|^2.
\]
Next from Lemma \ref{LowBndOnSigpLp}, it is immediate that
\begin{eqnarray*}
\mathbf{x}^T A_{22}^{(k)} \mathbf{x} &=& \int_{\tau_k} [x_1 N_{k-1} + x_2 N_{k}]^2 {[\sigma^{(k)}]^\prime}^2\, dx \\
&\ge& C^*_1 G_k^2 h^3 (y_1^2 + y_2^2) = C_1^* \|[D^{(k)}]^{-1}\mathbf{x}\|^2 .
\end{eqnarray*}
Similar bounds fox  $A_{22}^{(k^*)}x^2$ for all $x \in \mathbb R$ also hold.
Thus Assumption \ref{asp2} is satisfied with $L_2 = C_1^*$ and $U_2=C_2^*$. \myqed

\medskip

Next, recalling that $\mathbf{A}_{22}$ is $k^* \times k^*$, we choose the diagonal matrix $\mathbf{D}_2=\mbox{diag}(d_1,d_2,\cdots,d_{k^*})$ with $d_j = (\mathbf{A}_{22})_{jj}^{-1/2}$. Clearly, the diagonal elements of $\widehat{\mathbf{A}}_{22} = \mathbf{D}_2\mathbf{A}_{22}\mathbf{D}_2$ are equal to 1. Note that $(\mathbf{A}_{22})_{jj}$, $1\le j \le k^*$, is associated with the vertex $x_{j-1} \in \mcal{T}_2$. Also, $(\mathbf{A}_{22})_{11}=b_{11}^{(0)}$ and $(\mathbf{A}_{22})_{jj} = b_{22}^{(j-1)}+b_{11}^{(j)}$ for $2 \le j \le k^*$.

Now, for $x_i \in \mcal{T}_2$ and $i \ne 0,k^*-1$ (recall $x_i \notin \mcal{T}_2$ for $k^*\le i \le N$), we have $\mcal{K}_i = \{\tau_i,\tau_{i+1}\}$, $l(i,i)=2$, $l(i,i+1)=1$ and $\mcal{K}_i^* = \mcal{K}_i$. Also  $\mcal{K}_0 = \{\tau_1\}$, $l(0,1)=1$, $\mcal{K}_0^* = \mcal{K}_0$ and $\mcal{K}_{k^*-1} = \{\tau_{k^*}\}$, $l(k^*-1,k^*)=1$, $\mcal{K}_{k^*-1}^* = \mcal{K}_{k^*-1}$. Therefore, from \eqref{Ai}, we have $\Delta_i = [\delta_2^{(i)}]^{-2}+[\delta_1^{(i+1)}]^{-2}$ for $x_i \in \mcal{T}_2$ and $i\ne 1,k^*-1$; also $\Delta_0 = [\delta_1^{(1)}]^{-2}$ and $\Delta_{k^*-1} = [\delta_1^{(k^*)}]^{-2}$.

We now show that Assumption \ref{asp3} is satisfied.

\begin{myprop}\label{SingA3}
Let $(\mathbf{A}_{22})_{jj}$, $1 \le j \le k^*$, be the diagonal elements of $\mathbf{A}_{22}$ and consider $\Delta_i$ for $x_i \in \mcal{T}_2$, defined above. Then Assumption \ref{asp3} is satisfied.
\end{myprop}

\noindent \emph{Proof:} Let $x_i \in \mcal{T}_2$ and $i\ne 0,k^*-1$; $x_i$ is associated with $(\mathbf{A}_{22})_{j_ij_i}$, where $j_i = i+1$. Therefore using the definition of $(\mathbf{A}_{22})_{i+1,i+1}$, $\delta_1^{(i+1)}$, and $\delta_2^{(i)}$, we have
\begin{equation}
(\mathbf{A}_{22})_{j_ij_i}^{-1} \Delta_i = (\mathbf{A}_{22})_{i+1,i+1}^{-1} \Delta_i = \frac{G_i^2h^3}{b_{22}^{(i)}+b_{11}^{(i+1)}} + \frac{G_{i+1}^2h^3}{b_{22}^{(i)}+b_{11}^{(i+1)}} \, .\label{PrfA3-1}
\end{equation}
Now using Lemmas \ref{LowBndOnSigpLp}, \ref{BndsOnNormSigmaPrime}, and \eqref{BndOnRatio}, it is immediate that
\begin{eqnarray*}
{C_3^*}^2 &\le& \frac{b_{22}^{(i)}}{G_{i}^2h^3} \le {C_2^*}^2 , \\
\frac{{C_3^*}^2}{3^{2(2-\alpha)}} \le \frac{b_{11}^{(i+1)}}{G_{i+1}^2h^3}\frac{G_{i+1}^2}{G_{i}^2} &=& \frac{b_{11}^{(i+1)}}{G_{i+1}^2h^3} \le {C_2^*}^2,
\end{eqnarray*}
and therefore,
\begin{equation}
\frac{1}{2{C_2^*}^2} \le \frac{G_{i}^2h^3}{b_{22}^{(i)}+b_{11}^{(i+1)}} \le \frac{3^{2(2-\alpha)}}{{C_3^*}^2(1+3^{2(2-\alpha)})} \, .\label{PrfA3-2}
\end{equation}
Similarly, we get
\[
\frac{1}{{C_2^*}^2(1+3^{2(2-\alpha)})} \le \frac{G_{i+1}^2h^3}{b_{22}^{(i)}+b_{11}^{(i+1)}} \le \frac{1}{2{C_2^*}^2},
\]
and combining \eqref{PrfA3-1},\eqref{PrfA3-2}, we infer that there exist constants $L_3$, $U_3$, such that
\[
L_3 \le (\mathbf{A}_{22})_{j_ij_i}^{-1} \Delta_i \le U_3,
\]
where
\begin{eqnarray*}
&&L_3 = \frac{1}{2{C_2^*}^2} + \frac{1}{{C_2^*}^2(1+3^{2(2-\alpha)})} \, ,\\
&&U_3 = \frac{1}{2{C_3^*}^2} + \frac{3^{2(2-\alpha)}}{{C_3^*}^2(1+3^{2(2-\alpha)})} \, .
\end{eqnarray*}
Thus Assumption \ref{asp3} hold for $x_i \in \mcal{T}_2$, $i \ne 0,k^*-1$. The proofs for $x_i \in \mcal{T}_2$, $i = 0,k^*-1$ are simpler and we do not include them here, \myqed

Based on Propositions \ref{SingA2}, \ref{SingA3}, it is clear that Assumptions \ref{asp2} and \ref{asp3} are satisfied. Assumption \ref{asp1} always hold in 1-d. Thus from Theorem \ref{ThCondNo}, $\mathfrak{K}(\mathbf{A}) = O(h^{-2})$ and the GFEM with $\mcal{S}=\mcal{S}_1 + \overline{\mcal{S}}_2$ is an SFEM.

%% file: sec5p3_appl.tex
\subsection{Problems with discontinuous solutions} \label{DiscontProb} We now address a problem, which is different from \eqref{VarProb}. Let $\Omega = (0,1)$ and set $\Omega_l = (0,c)$ and $\Omega_r=(c,1)$, where $0 < c < 1$ is fixed. Consider
\begin{eqnarray*}
&&H(\Omega):= \left\{v \in L_2(\Omega):\, v(0)=v(1)=0,\, \int_{\Omega_l}{v^\prime}^2 dx < \infty,\right. \\
&& \hspace{5cm} \left. \mbox{ and } \int_{\Omega_l}{v^\prime}^2 dx < \infty \right\}.
\end{eqnarray*}
Then $(H(\Omega), \|\cdot\|_H)$ is a Hilbert space, where
\[
\|v\|_H^2 := |v|_{H^1(\Omega_l)}^2 + |v|_{H^1(\Omega_r)}^2.
\]
We note that $H^1_0(\Omega) \subset H(\Omega)$ and functions in $H(\Omega)$ may have jump discontinuity at $x=c$.

For $f \in L_2(\Omega)$, we consider the problem
\begin{equation}
u \in H(\Omega),\ B(u,v) = F(v),\quad \forall \ v \in H(\Omega), \label{VarProb2}
\end{equation}
where
\[
B(u,v):= \int_{\Omega_l}u^\prime v^\prime \, dx + \int_{\Omega_r}u^\prime v^\prime \, dx \mbox{ and } F(v):= \int_{\Omega} fv\, dx.
\]
The bilinear form $B(\cdot,\cdot)$ is coercive and bounded in $H(\Omega)$. Also $F(\cdot)$ is a bounded linear functional on $H(\Omega)$. Thus the problem \eqref{VarProb2} has a unique solution.

If $f$ and the solution $u\in H(\Omega)$ of \eqref{VarProb2} are smooth in $\Omega_l$ and $\Omega_r$, then $u$ is the solution of the boundary value problem
\begin{eqnarray*}
&& - u^{\prime \prime} = f  \mbox{ on } \Omega_l,\quad -u^{\prime \prime} = f \mbox{ on } \Omega_r ,\\
&& \ u(0) = u(1) = 0 \mbox{ and } u^\prime(c^-) = u^\prime(c^+) = 0.
\end{eqnarray*}
This problem mimics the problem with a crack in higher dimensions, where the solution is discontinuous along the crack line away from the crack-tip.

We now give a characterization of the solution of \eqref{VarProb2}. We will use the Heaviside function
\begin{equation}
H_c(x) = \left\{ \begin{array}{rc}
                1, & 0 \le x < c ;\\
                -1, & c \le x \le 1. \end{array} \right. \label{heavi}
\end{equation}

\begin{mylemma} \label{DiscontSolnLem}
Suppose $u \in H(\Omega)$ such that $u^\prime(c^-) = u^\prime(c^+) = 0$ and
\[
\int_{\Omega_l} (u^{\prime \prime})^2 dx <\infty,\quad \int_{\Omega_l} (u^{\prime \prime})^2 dx <\infty.
\]
Then
\[
u(x) = s(x) + \tilde{u}(x),
\]
where $s$ is a step function with discontinuity at $x=c$ and $\tilde{u} \in H^2(\Omega)$
\end{mylemma}

Proof: We first note that $u(c^-)$ and $u(c^+)$ are well defined. We define
\begin{equation}
\tilde{u} = u -\frac{u(c^-) - u(c^+)}{2}\, [H_c - 1], \label{utilda}
\end{equation}
where $H_c(x)$ is given in \eqref{heavi}. It is easy to check that
\begin{eqnarray*}
&&\tilde{u}\big|_{\Omega_l} = u\big|_{\Omega_l},\ \tilde{u}(c^-) = u(c^-), \\
&&\tilde{u}^\prime\big|_{\Omega_l} = u^\prime\big|_{\Omega_l},  \ \mbox{ and } \ \tilde{u}^\prime(c^-) = u^\prime(c^-).
\end{eqnarray*}
Similarly,
\begin{eqnarray*}
&&\tilde{u}\big|_{\Omega_r} = u\big|_{\Omega_r}+[u(c^-) - u(c^+)],\ \tilde{u}(c^+) = u(c^-), \\ &&\tilde{u}^\prime\big|_{\Omega_r} = u^\prime\big|_{\Omega_r},  \ \mbox{ and } \ \tilde{u}^\prime(c^+) = u^\prime(c^+).
\end{eqnarray*}
Note that $\tilde{u}(c^-)=\tilde{u}(c^+)= u(c^-)$. We define $\tilde{u}(c) = u(c^-)$ so that $\tilde{u}$ is continuous at $x=c$ and thus is continuous on $\Omega$. Also since $u^\prime(c^-) = u^\prime(c^+) =0$, it is clear from above that $\tilde{u}^\prime(c^-) = \tilde{u}^\prime(c^+)=0$. We define $\tilde{u}^\prime(c) = 0$ so that $\tilde{u}^\prime$ is continuous at $x=c$, and consequently $\tilde{u}^\prime$ is continuous in $\Omega$. Moreover,
\[
\int_\Omega (\tilde{u}^{\prime \prime})^2 dx = \int_{\Omega_l} (\tilde{u}^{\prime \prime})^2 dx + \int_{\Omega_r} (\tilde{u}^{\prime \prime})^2 dx = \int_{\Omega_l} ({u}^{\prime \prime})^2 dx + \int_{\Omega_r} ({u}^{\prime \prime})^2 dx \le \infty.
\]
Thus $\tilde{u} \in H^2(\Omega)$ and considering $s = \frac{u(c^-) - u(c^+)}{2}\, [H_c - 1]$ in \eqref{utilda}, we get the desired result. \myqed

Suppose $c \notin \mcal{T}$, i.e., $c$ is not a vertex of the mesh. Therefore, there exists an $m$ such that $ c \in \tau_{m+1}$ and hence, $c \in \omega_m \cap \omega_{m+1}$. We assume that $m \ne 1,N$; this is always achieved for $h$ small. Since $u(0)=u(1)=0$, we consider $\mcal{T}_1 = \mcal{T}\backslash \{x_0,x_N\}$ (see Remark \ref{RemT1T2}).

For $1 \le i \le N-1$, we consider $V_i = \mbox{span}\{1, \locphi{1}{i}=(x-x_i),\locphi{2}{i}=H_c(x)\}$ and we set $V_0 = \mbox{span}\{\locphi{1}{0}=(x-x_0)\}$ and $V_N = \mbox{span}\{\locphi{1}{N}=(x-x_N)\}$. Note that $V_i \in H(\omega_i)$ for $i\in I$ (i.e., for $x_i \in \mcal{T}$). Clearly, $\overline{V}_i = \{0\}$ for $i\in I\backslash \{ m, m+1\}$. We set $\mcal{T}_2=\{x_m,x_{m+1}\}\subset \mcal{T}$ and define
\[
\overline{\mcal{S}}_2 = N_m \overline{V}_m^1 + N_{m+1} \overline{V}_{m+1}^1.
\]
We consider the GFEM with $\mcal{S} = \mcal{S}_1 + \overline{\mcal{S}}_2$.

Since $\locphi{2}{m}=H_c$ is constant in $\tau_{m}$, we have $\locphibar{2}{m}|_{\tau_{m}} = 0$. Similarly, $\locphibar{2}{m+1}|_{\tau_{m+2}} = 0$. Therefore $\mcal{K}_{enr}=\{\tau_{m+1}\}$. Moreover, the functions $\locphibar{2}{m},\, \locphibar{2}{m+1}$ are discontinuous at $x=c$, their values are zero at $x=x_m, x_{m+1}$, and $\locphibar{2}{m}|_{\tau_{m+1}} = \locphibar{2}{m+1}|_{\tau_{m+1}}$.

We assume that $f$ is such that solution $u \in H(\Omega)$ of \eqref{VarProb2} satisfies the assumptions of Lemma \ref{DiscontSolnLem} and $u = s + \tilde{u}$, where $s$ is a step-function with a discontinuity at $x=c$ and $\tilde{u} \in H^2(\Omega)$. We now address the convergence of the GFEM solution. We first note that Theorem \ref{ApproxThmSGFEM1d} hold for $u\in H(\Omega)$ with $\mcal{E}(\Omega), \mcal{E}(\omega_i)$ replaced by $H(\Omega), H(\omega_i)$ and with $\overline{V}_i \in H(\omega_i)$. Now for $x_i \in \mcal{T}\backslash \mcal{T}_2$, we have $u \in H^2(\omega_i)$ and from the standard interpolation result
\begin{equation}
\|u-\intp{\omega_i}u \|_{H(\omega_i)} = \|u-\intp{\omega_i}u \|_{\mcal{E}(\omega_i)}\le Ch|u|_{H^2(\omega_i)} = Ch|\tilde{u}|_{H^2(\omega_i)}. \label{StepConv:1}
\end{equation}
For $x_m \in \mcal{T}_2$, it is easy to show that there exists $\bar{\xi}^m \in \overline{V}_m$ such that $u- \intp{\omega_m}u - \bar{\xi}^m = \tilde{u}-\intp{\omega_m}\tilde{u}$ on $\omega_m$.
Therefore, $\|u- \intp{\omega_m}u - \bar{\xi}^m\|_{L^2(\omega_m)} \le C|\omega_m|\,\|u- \intp{\omega_m}u - \bar{\xi}^m\|_{\mcal{E}(\omega_m)}$, and from the standard interpolation result, we have
\begin{equation}
\|u- \intp{\omega_m}u - \bar{\xi}^m\|_{H(\omega_m)} = \|\tilde{u} - \intp{\omega_m}\tilde{u}\|_{H^1(\omega_m)} \le Ch|\tilde{u}|_{H^2(\omega_m)}. \label{StepConv:2}
\end{equation}
Similarly, there exists $\bar{\xi}_{m+1} \in \overline{V}_{m+1}$ such that
\[
\|u- \intp{\omega_{m+1}}u - \bar{\xi}_{m+1}\|_{H(\omega_{m+1})} \le Ch|\tilde{u}|_{H^2(\omega_{m+1})}. \label{StepConv:2a}
\]
Therefore combining \eqref{StepConv:1}, \eqref{StepConv:2}, \eqref{StepConv:2a} and using the Theorem \ref{ApproxThmSGFEM1d} with modifications as mentioned above, we infer that there exists $v \in \mcal{S}=\mcal{S}_1+\overline{\mcal{S}}_2$ such that
\[
\|u - v\|_{H(\Omega)} \le Ch|\tilde{u}|_{H^2(\Omega)}.
\]
Therefore, $\|u-u_h\|_{H(\Omega)} = O(h)$, where $u_h$ is the GFEM solution.

We next address the scaled condition number of the stiffness matrix of the GFEM. Since $\tau_{m+1}$ is the only element in $\mcal{K}_{enr}$, we have
$\mathbf{A}_{22} = A_{22}^{(m+1)}$.

Let $c = x_m + \beta h$ with $0 < \beta < 1$. A direct computation yields that
\[
A_{22}^{(m+1)} = \frac{4}{h}\left[ \begin{array}{cc}
            (4 - 9\beta+6\beta^2)/3 & -\frac{1}{3}+2\beta - 2\beta^2 \\
            -\frac{1}{3}+2\beta - 2\beta^2 & (1-3\beta + 6\beta^2)/3
            \end{array} \right].
\]
Thus $A_{22}^{(m+1)}$ is associated with vertices $x_m$ and $x_{m+1}$. We choose the diagonal matrix $D^{(m+1)} = \mbox{diag}\{\delta_1^{(m+1)},\delta_2^{(m+1)}\}$ with $\delta_1^{(m+1)}=\delta_2^{(m+1)}=h^{1/2}/2$.
Then
\[
\hat{A}_{22}^{(m+1)} = D^{(m+1)} A_{22}^{(m+1)} D^{(m+1)}
= \left[ \begin{array}{cc}
    (4 - 9\beta+6\beta^2)/3 & -\frac{1}{3}+2\beta - 2\beta^2 \\
            -\frac{1}{3}+2\beta - 2\beta^2 & (1-3\beta + 6\beta^2)/3
            \end{array} \right].
\]
The diagonal elements of $\hat{A}_{22}^{(m+1)}$ are $O(1)$ for any $0 < \beta < 1$. The eigenvalues of $\hat{A}_{22}^{(m+1)}$ are $\lambda_1 = \frac{5}{6}-2\beta+2\beta^2 -T$ and $\lambda_2 = \frac{5}{6}-2\beta+2\beta^2 +T$, where $T = \frac{1}{6} \sqrt{13-84\beta+228\beta^2-288\beta^3+144\beta^4}$ (obtained from \emph{MAPLE}). It can be shown that
\[
\frac{1}{6} \le \lambda_1 \le \frac{5}{6} - \frac{\sqrt{13}}{6},\quad
\frac{1}{2} \le \lambda_2 \le \frac{5}{6} + \frac{\sqrt{13}}{6}.
\]
Thus, as before, we have
\[
\frac{1}{6}\|[D^{(m+1)}]^{-1}\mathbf{x}\|^2 \le \mathbf{x}^T A_{22}^{(m+1)} \mathbf{x} \le \Big( \frac{5}{6} + \frac{\sqrt{13}}{6} \Big)\|[D^{(m+1)}]^{-1}\mathbf{x}\|^2,\quad \forall \ \mathbf{x}\in \mathbb R^2,
\]
and the Assumption \ref{asp2} is satisfied with $L_2 = \frac{1}{6}$ and $U_2 = \frac{5}{6} + \frac{\sqrt{13}}{6}$.

We choose $\mathbf{D}_{2} = \mbox{diag}\{d_1,d_2\}$ with $d_i = (\mathbf{A}_{22})_{ii}^{-1/2}$. Clearly, the diagonal elements of $\widehat{\mathbf{A}}_{22} = \mathbf{D}_{2} \mathbf{A}_{22} \mathbf{D}_{2}$ are equal to 1. As in the first example (i.e., when $a(x)=a_1(x)$) in Section \ref{Interface}, we have $\mcal{T}_2=\{x_m,x_{m+1}\}$ and $\mcal{K}_m=\mcal{K}_{m+1}=\{\tau_{m+1}\}$. Therefore, $l(m,m+1)=1$ and $l(m+1,m+1)=2$. Also $\mcal{K}_m^*=\mcal{K}_{m+1}^*=\{\tau_{m+1}\}$ and therefore from \eqref{Ai}, we have $\Delta_m=[\delta_1^{(m+1)}]^{-2}$ and $\Delta_{m+1}=[\delta_2^{(m+1)}]^{-2}$. Also $x_m,x_{m+1} \in \mcal{T}_2$ are associated with $(\mathbf{A}_{22})_{y_my_m}, (\mathbf{A}_{22})_{y_{m+1}y_{m+1}}$, respectively, where $y_m=1,$ and $y_{m+1}=2$. It is easy to check that
\[
\frac{3}{4} \le (\mathbf{A}_{22})_{11}^{-1}\Delta_m,\ (\mathbf{A}_{22})_{22}^{-1}\Delta_{m+1} \le \frac{24}{5}.
\]
Thus Assumption \ref{asp3} is satisfied with $L_3 = 3/4$ and $U_3=24/5$. Therefore from Theorem \ref{ThCondNo}, we have $\mathfrak{K}(\mathbf{A}) = O(h^{-2})$, where $\mathbf{A}$ is the stiffness matrix, for all $0 < \beta < 1$.


%% file: sec_conclu.tex
\section{Conclusion}\label{conclusion}
The GFEM uses special enrichment functions, based on the available (or extracted) information on the unknown solution of the underlying variational problem. The use of special enrichment functions gives rise to the excellent convergence properties of the GFEM. In fact, for a given problem, it is possible to choose several classes of enrichment functions such that the GFEM, employing each of these enrichment classes, will yield excellent convergence properties. However, GFEM employing some of these classes of enrichments could be ill-conditioned, i.e., there could be severe loss of accuracy in the computed solution of the linear system associated with the GFEM. The loss of accuracy could be much more than that experienced in a standard FEM. In this paper, we have presented and analyzed a modification of the GFEM -- the stable GFEM (SGFEM), which does not have the problem with severe loss of accuracy. SGFEM has all the advantages of the GFEM and is also very robust with respect to the parameters of the enrichments (e.g., the parameter $\beta$ in Sections \ref{Interface} and \ref{DiscontProb}). The loss of accuracy is characterized by the scaled condition number and is expressed through Hypothesis H, which was validated based on various examples.

The abstract framework developed in this paper has been applied to a one-dimensional problem for the clarity of exposition. This framework could also be applied to higher dimensional problems, which will be reported in a forthcoming paper. 

\medskip

\noindent \textbf{Acknowledgement:}  We thank Professor C. Armando Duarte of the Department of Civil and Environmental Engineering, University of Illinois at Urbana-Champaign, and Professor John E. Osborn of the Department of Mathematics, University of Maryland at College Park, for fruitful discussions on certain aspects of this paper. We also thank Mr. Karl Schulz of PECOS, ICES, University of Texas at Austin, for his help on the use of \emph{superLU} and \emph{MUMPS} on the \emph{Lonestar system}, Texas Advanced Computing Center, and also for various illuminating discussions and interpretation of the results, presented in this paper.   

%% file: sec6_appendix.tex
\section{Appendix} \label{Appendix}

\emph{Validation and Verification} (V \& V) is a fairly new field and is still in its developing stage (\cite{ASME, RoachBook, BabOden, OberRoyBook}). Suppose a
mathematical model of some ``Reality'' (e.g., a physical, chemical or biological system or process), formulated for a particular goal or purpose, is given. The objective of V \& V is to assess whether the predictions based on the computed solution of a mathematical model are reliable enough so that they could be the basis for certain decisions related to the goal.

Validation is the process of building confidence on the mathematical model (\cite{ASME, RoachBook, OberRoyBook}). The process is of course is constrained by the cost, available time, and skills, as explicitly underlined in \cite{SISO}. It is based on a set of properly selected problems and their mathematical models for which experimental data is available. These problems are called validation problems and they are chosen with varying level of complexity; more complex problems are closer to the ``Reality''. Of course, obtaining the experimental data for the validation problems with increasing complexity is increasingly costly. The prediction based on the computed solution of these problems is then compared with the experimental data. The assessment of the difference is based on a specified tolerance and a suitably selected metric (could be more than one) relative to the specific goal. If the measure of the difference is larger than the tolerance for any validation problem, the mathematical model is rejected. If none of the validation models are rejected, then one could have confidence that the mathematical model realistically describes the ``Reality'', with respect to the goal, beyond the scope of the chosen validation problems. The level of confidence will be based on the tolerance as well as the number and the selection of the validation problems. We mention that the set of the validation problems is finite, their selection has a large subjective component, and a philosophical question about the justification of the confidence in the mathematical model could certainly be raised (see \cite{KleinEtal}).

Numerical algorithms and their properties obtained from the mathematical analysis are always based on various assumptions that are not satisfied when the algorithm is implemented on a computer. For example, infinite precision arithmetic is often assumed while describing a numerical algorithm or stating an inference about the algorithm obtained from the analysis. However, this assumption is always violated by a computer working with finite precision arithmetic. The output from the computer implementation of the algorithm may also depend, for example, on the package in which the algorithm have been implemented, the compiler, the processor, the computing platform with single or multiple processors, among other factors. Consequently, the output may vary even when the same outcome is predicted by the mathematical analysis for two different algorithms. For example, suppose the solution of the linear system $Ax=b$ is sought using algorithms of the form $P_i^TAP_iz = P_i^Tb$, $x=P_iz$, where $P_i$ is a permutation matrix for $i=1,2$. Both the algorithms should yield that solution $x$, however, the computed solutions could be different (see Problems 1a,b, below). Thus the implementation of a numerical algorithm in a computer is analogous to a ``Reality''; the goal is to obtain a particular quantity of interest for a particular purpose (related to a decision). The mathematical model of this ``Reality'' is the inference obtained from the mathematical analysis, or other statements based on the inference, about obtaining the quantity of interest from the algorithm. Therefore, the process of validation of the inference has to be performed to have confidence in the inference or a statement based on the inference.

We have briefly formulated the following hypothesis in the Introduction. Let $Ax=b$, $x,b \in \mathbb R^n$ be a linear system, where the $n \times n$ matrix $A$ belongs to a class of sparse matrices that include the stiffness matrices associated with FEM, GFEM, or SGFEM. Let $\hat{x}$ be the computed solution of the linear system, obtained from an elimination method, e.g., some variant of Gaussian elimination. Moreover, $\hat{x}$ is computed in finite precision arithmetic with machine precision $\epsilon$. Let $H = DAD$ where $D$ is a diagonal matrix with $D_{ii}=A_{ii}^{-1/2}$; clearly, $H_{ii}=1$. Recall the scaled condition number $\mathfrak{K}(A)$ of $A$ is given by
$
\mathfrak{K}(A):= \kappa_2(H),
$
where $\kappa_2(H)=\|H\|_2\|H^{-1}\|_2$ is the condition number of $H$ based on the $\|\cdot \|_2$ vector norm. Also recall $\eta :=\|x-\hat{x}\|_2/\|x\|_2$.

\begin{quote}
\textbf{Hypothesis H:} For $n$, not small,
\begin{equation}
\eta \approx Cn^\beta \mathfrak{K}(A) \epsilon; \ \beta \approx 0, \label{HypH1}
\end{equation}
where $\hat{x}$ has been computed in an computing environment satisfying the IEEE standard for floating point arithmetic (with the guard digit), there is no overflow or underflow during the computation of $\hat{x}$, and
$C$, $\beta$ do not depend on $n$ as well as other factors mentioned before. \end{quote}
The $\approx$ in \eqref{HypH1} means there exist $0<\bar{C}_1,\,\bar{C}_2$ and $0<\bar{\beta}$ small, such that $\eta = Cn^\beta \mathfrak{K}(A) \epsilon$ with $\bar{C}_1 \le C \le \bar{C}_2$ and $ |\beta| \le \bar{\beta}$. Also this hypothesis addresses the range $N$ for which not (almost) all digits of accuracy is lost (see Problem 3a).

Hypothesis H is based on certain mathematical inferences (results), which we will discuss later. The validation of \eqref{HypH1} with respect to the tolerance $\tau=\{\tau_1,\tau_2\}$ means that $ \bar{C}_2 / \bar{C}_1 \le \tau_1$ and $ \bar{\beta} \le \tau_2$. Note that $\tau_2$ is primary and should be small for confidence in \eqref{HypH1}, however $\tau_1$ could be allowed to be larger. The set of validation problems consists of stiffness matrices of FEM, GFEM, SGFEM, and other similar matrices, e.g., arising in finite difference method, applied to solve various linear elliptic variational problems of increasing complexity. For confidence in the Hypothesis, we require that \eqref{HypH1} is not rejected for any of the validation problems relative to the given tolerance $\tau$. We note that it is possible to select a tolerance such that the hypothesis is not rejected, however, the tolerance have to be admissible (e.g., reasonably small) for the decision making process. In our case, the decision will be whether to accept the SGFEM over the standard GFEM. We note that the class of matrices for which the hypothesis will be validated is not precisely defined, similar to a class of complex physical or engineering problem.

We now give a theoretical rationale for \eqref{HypH1}. There is a lot of literature available on the accuracy of the computed solutions of the linear system $Ax=b$. We particularly mention the classic \cite{WilkinsonBook} and a modern book \cite{HighamBook} with an excellent survey of the theoretical results in the area. Typically, the loss of accuracy in the numerical solution due to round-offs is analyzed by the backward error analysis. This analysis shows that the computed solution is the exact solution of a perturbed linear system, and it provides estimates of the perturbations in terms of the data of the linear system. A bound on the loss of accuracy in the computed solution, measured by $\eta$ defined before, is then obtained using the perturbation estimates.

It is well known from a standard perturbation argument that for a full matrix $A$,
\begin{equation}
\eta \le f(n)\kappa_2(A)\epsilon, \label{BndOnEta}
\end{equation}
where $\epsilon$ is the machine precision and $f(n)$ depends on the algorithm used to solve $Ax=b$ (see e.g., \cite{GolubVL, KinCheney}). In Hypothesis H, we hypothesize that $\kappa_2(A)$ is replaced by $\mathfrak{K}(A)$. We also hypothesize that  $\bar{C}_1n^\beta \le f(n) \le \bar{C}_3n^\beta$ and $ |\beta| \le \bar{\beta}$, where $\bar{C}_1$, $\bar{C}_2$, and $\bar{\beta}$ are as defined before. It is important to note that in the mathematical literature, only an upper bound of $\eta$ is available; in contrast, the Hypothesis H addresses both the upper and lower bounds of $\eta$.

Consider the linear system $DADz = Db$, where $D$ is a diagonal matrix with $D_{ii}=2^{g_i}$ in the rage of the floating point system. Clearly $x=Dz$. We now cite the following old result of F. L. Bauer (\cite{Bauer}):
\begin{mytheorem} \label{Bauer} Let $\hat{x}$, $\hat{z}$ be the computed solutions of the linear systems $Ax=b$ and $DADz=Db$, respectively, obtained from an elimination method with no pivoting. Furthermore, we assume that there is no overflow or underflow in the computation of $\hat{x}$, $\hat{z}$. Then all the digits of $\hat{x}$ and $\hat{x}_D$ are same, where $\hat{x}_D = D\hat{z}$.
\end{mytheorem}
We note that the result of the above Theorem is not true if the diagonal elements of $D$ are not binary. However in that situation, the quantities $\|x-\hat{x}\|$, $\|x-\hat{x}_D\|$, and $\|\hat{x}-\hat{x}_D\|$ are of the same order.

We next note that it is possible to find a diagonal matrix such that $\kappa_2(DAD)  \le \kappa_2(A)$. For example, for $\mu > 0$, let
$A = \left[ \begin{array}{cc} 1 & 1 \\ 0 & \mu \end{array} \right]$. Then $\kappa_2(A) = 1/\mu$ is large for $\mu$ small. Let $\mu = \chi 2^{-d}$, where $1/2 < \chi < 2$. Consider $D = \left[ \begin{array}{ll} 1 & 0 \\ 0 & 2^{d/2} \end{array} \right]$. Then $DAD = \left[ \begin{array}{ll} 1 & 2^{d/2} \\ 0 & \chi \end{array} \right]$, and $1 \le \kappa_2(DAD) < 2$ and consequently, $ \kappa_2(DAD) < \kappa_2(A)$ for $\mu$ small. Let $D^*$ be the diagonal matrix such that $\kappa_2(D^*AD^*) = \min_{D}\kappa_2(DAD)$ (minimum over all diagonal matrices $D$ with binary diagonal elements), then from the above theorem and \eqref{BndOnEta}, we have
\[
\eta \le f(n) \, \kappa_2(D^*AD^*) \epsilon \le f(n)\, \kappa_2(A) \epsilon.
\]
Thus $\kappa_2(D^*AD^*)$ provides more accurate information about $\eta$ than $\kappa_2(A)$. But in general, it is not easy to find either $D^*$ or $\kappa_2(D^*AD^*)$. In Hypothesis H, we used $\mathfrak{K}(A) = \kappa_2(DAD)$, where $D$ is a diagonal matrix with $D_{ii} = A_{ii}^{-1/2}$ and $D_{ii}$ may not be binary. We note, however, that not using a binary only influences $\bar{C}_1$, $\bar{C}_2$ by factors of $1/2$ and $2$ respectively. We also mention that in the literature (\cite{Demmel, HighamBook}), an upper bound of the form \eqref{BndOnEta} for $\eta$ is available with $f(n)=Cn^2 $ and $\kappa_2(A)$ replaced by $\mathfrak{K}(A)$ for symmetric positive definite linear systems solved by Cholesky decomposition. In Hypothesis $H$, we used $f(n)=Cn^\beta$, $\beta\approx 0$, based on our computational experience.

We now consider a set of validation problems, whose exact solution (experimental data) is known. The solution to these problems will be computed on various computers using double precision, i.e., with 16 digits of accuracy.
\begin{description}
\item[Problem 1:] We consider approximating the solution $u(x) = x$ of the problem $-u^{\prime \prime}(x) = 0, \ x\in (0,1), u(0)=0, u(1)=1,$ by the FEM using piecewise linear finite elements.
    \begin{description}
    \item[Problem 1a:] We use the FE mesh vertices $x_i=ih$ for $i=0,1,\cdots, N$ and $h=1/N$. The FE solution is same as the exact solution $u$ of the problem. Let the associated linear systems be $A^{(1)} x^{(1)} = b^{(1)}$. The exact solution vector $x^{(1)}$ is known, namely, $x^{(1)}_i = ih$, $i=1,2,\cdots,N$. We will solve the linear system by the standard LU decomposition algorithm for sparse matrices without partial pivoting.
    \item[Problem 1b:] We use the mesh vertices $x_i=(N-i+1)h$, $i=0,1,\cdots, N$. The FE solution is same as the exact solution $u$ and let the associated linear system be $A^{(2)} x^{(2)} = b^{(2)}$; it is known that $x^{(2)}_i = (N-i+1)h$, $i=1,\cdots, N$. Note that the elements of $x^{(2)}$ are the permuted elements of $x^{(1)}$ and thus $\|x^{(1)}\|_2 = \|x^{(2)}\|_2$. We will solve the linear system by the same algorithm as Problem 1a.
    \end{description}
    The computations are performed on a Dell Latitude PC with INTEL CORE(TM)2 CPU, 1.20GHZ.
\item[Problem 2:] We approximate the solution $u(x) = 1$ of the problem $-u^{\prime \prime}(x) = 0, \ x\in (0,1), u(0)=1, u(1)=1,$ by the piecewise linear FEM  based on the mesh vertices as in Problem 1a. Let the associated linear system be $Ax=b$. It is clear that the exact solution is given by $x_i=1$, $i=1,2,\cdots,N$.
    \begin{description}
    \item[Problem 2a:] The linear system is solved by a sparse matrix direct solver \emph{superLU} \cite{superlu} on a single processor.
    \item[Problem 2b:] The linear system is solved by a sparse matrix direct solver \emph{MUMPS} \cite{MUMPS} on a single processor.
    \item[Problem 2c:] The linear system is solved by \emph{MUMPS}, using parallel computation, on 128 processors.
    \end{description}
    The computations were performed on the Lonestar system at Texas Advanced Computing Center. Lonestar is a Linux based cluster comprised of 1888 compute nodes connected via high speed quad-data rate infiniband, with each compute node containing two hex-core socket (INTEL Xeon 5680 processors) for an aggregate system size of 22656 cores. Each core runs at a peak of 3.33GHZ.
\item[Problem 3:] We consider approximating the solution $u(x) = x^2$ of the problem $-u^{\prime \prime}(x) = -2, \ x\in (0,1)$, $u(0)=0,u^\prime(1)=2$, by the GFEM based on $\mcal{S}= \mcal{S}_1 + \mcal{S}_2$ (see \eqref{GFEMspace2}). We use $n_i = 1$ and $\locphi{1}{i}(x) = x^2$, $i=0,1,\cdots,N$. We order the shape functions as $N_0\locphi{1}{0}$, $N_1$, $N_1\locphi{1}{1}$, $N_2$,$\cdots$, $N_N$, $N_N\locphi{1}{N}$ and suppose the associated stiffness matrix is $Ax=b$, where $A$ is of the order $2N+1$. The GFEM solution is same as the exact solution $u$. It is easy to see that $x_{2i+1}=1$, $i=0,2,\cdot,N$ and $x_{2i}=0$, $i=1,2,\cdots,N$. The linear system is solved by the same algorithm and on the same platform as in Problem 1a.
\item[Problem 4:] We consider approximating the solution of the same problem in Problem 3 by the SGFEM based on $\mcal{S}= \mcal{S}_1 + \overline{\mcal{S}}_2$ (see \eqref{SGFEMspaceMain}) with $n_i =1$, $\mcal{T}_2 = \mcal{T}$, and $\locphibar{1}{i}= x^2-\intp{\omega_i} x^2$. We order the shape functions as $N_0\locphibar{1}{0}$, $N_1$, $N_1\locphibar{1}{1}$, $N_2$,$\cdots$, $N_N$, $N_N\locphibar{1}{N}$ and suppose the associated stiffness matrix is $Ax=b$, where $A$ is of the order $2N+1$. The GFEM solution is same as the exact solution $u$ and it is easy to see that $x_{2i+1}=1$, $i=0,2,\cdot,N$ and $x_{2i}=(ih)^2$, $i=1,2,\cdots,N$. The linear system is solved by the same method and on the same platform as in Problem 1a.
\end{description}

We will now validate Hypothesis H based on the validation problems described above. We will consider the tolerance $\tau=(\tau_1,\tau_2)$, with $\tau_1 = 400$ and $\tau_2=0$.

\begin{figure}[t]
\includegraphics[width=2.2in]{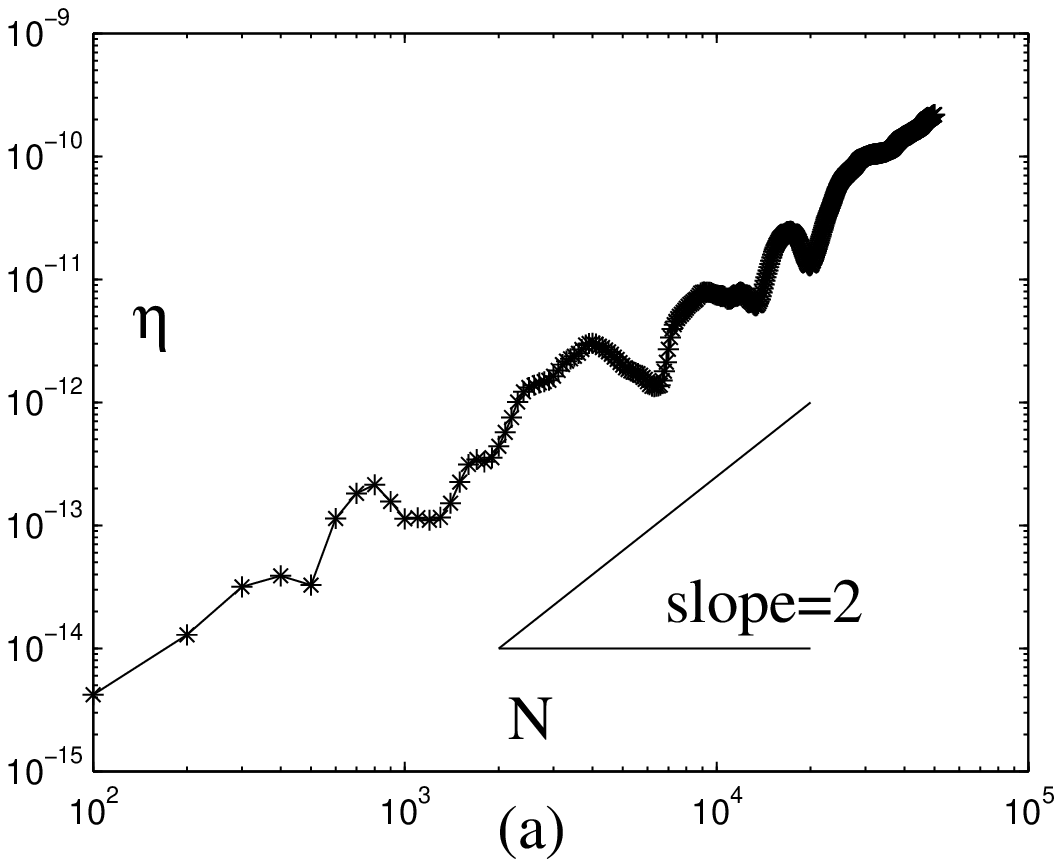}
\includegraphics[width=2.2in]{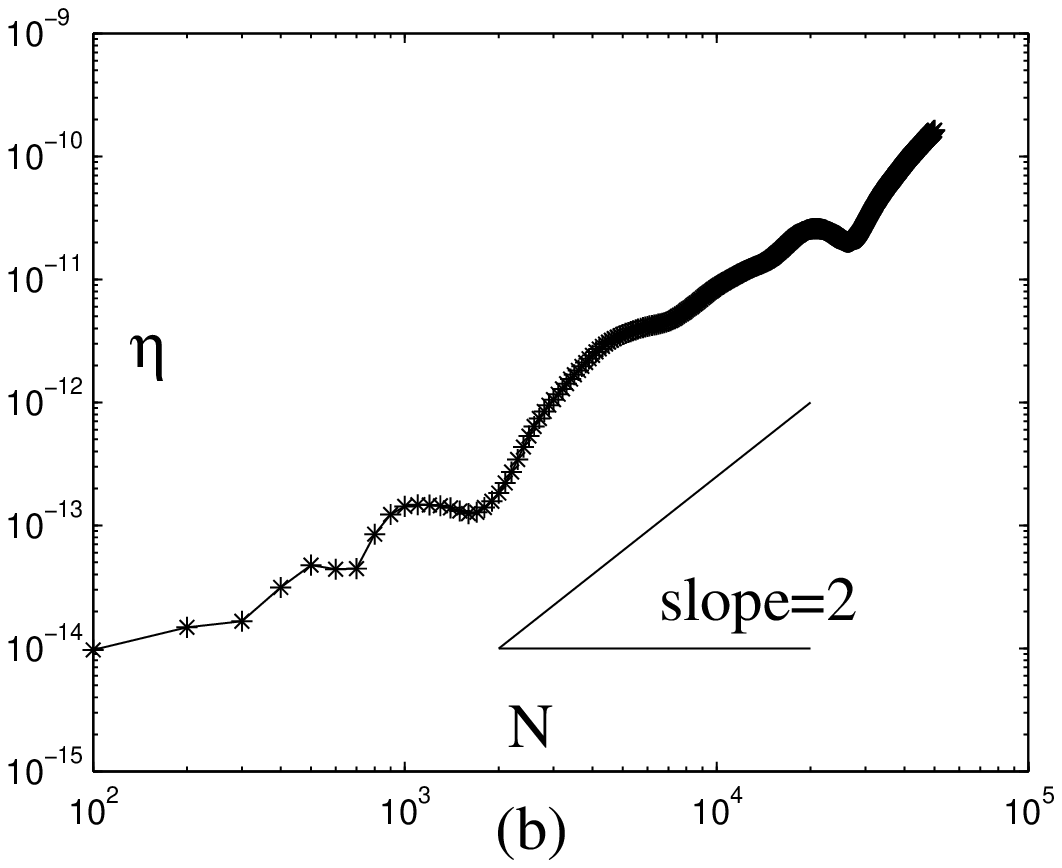}\\[-1ex]
\vspace{-.4cm}
\caption{\label{PlotRenumberNode} Log-log plots of $\eta^{(k)}=\|x^{(k)}-\hat{x}^{(k)}\|_2/\|x^{(k)}\|_2$ where $\hat{x}^{(k)}$ is the computed solution of $A^{(k)} x^{(k)} = b^{(k)}$, $k=1,2$, associated with FEM with vertices $x_i=ih$ and $x_i=(N-i)h$, $i=0,1,\cdots,N$, respectively. $\eta^{(1)}$, $\eta^{(2)}$ have been computed and presented in (a) and (b), respectively, for $N=100, 200, \cdots, 50000$}.
\end{figure}
Let $\hat{x}^{(1)}$ and $\hat{x}^{(2)}$ be the computed solutions of the linear systems $A^{(1)} x^{(1)} = b^{(1)}$ and $A^{(2)} x^{(2)} = b^{(2)}$ of \textbf{Problem 1a} and \textbf{Problem 1b}, respectively. It can be shown that for large $N$, $\mathfrak{K}(A^{(1)})=\mathfrak{K}(A^{(2)}) \approx 0.4\, N^2$. We have computed and presented the log-log plots of the relative errors $\eta^{(k)}=\|x^{(k)}-\hat{x}^{(k)}\|_2/\|x^{(k)}\|_2$, $k=1,2$, with respect to $N=100, 200, \cdots,50000$ in Figure \ref{PlotRenumberNode}. We have observed that
$\bar{C}_1^{(k)} [0.4 N^2] \le \eta^{(k)} \le \bar{C}_2^{(k)} [0.4 N^2]$ for $k=1,2$ with $\bar{C}_2 / \bar{C}_1 \le 120 < \tau_1$ (note $\tau_2 =0$). Thus we do not reject Hypothesis H. Note that we did not reject the hypothesis based only on the subset of meshes with the values of $N$, mentioned above. Moreover, it is interesting to note that the plots of $\eta^{(1)}$ and $\eta^{(2)}$ are quite different. Thus the computed solution is affected by changing the order of the FE mesh vertices, in spite of the fact that $\|x^{(1)}\|_2 = \|x^{(2)}\|_2$.

\begin{figure}[t]
\begin{center}
\includegraphics[width=2.2in]{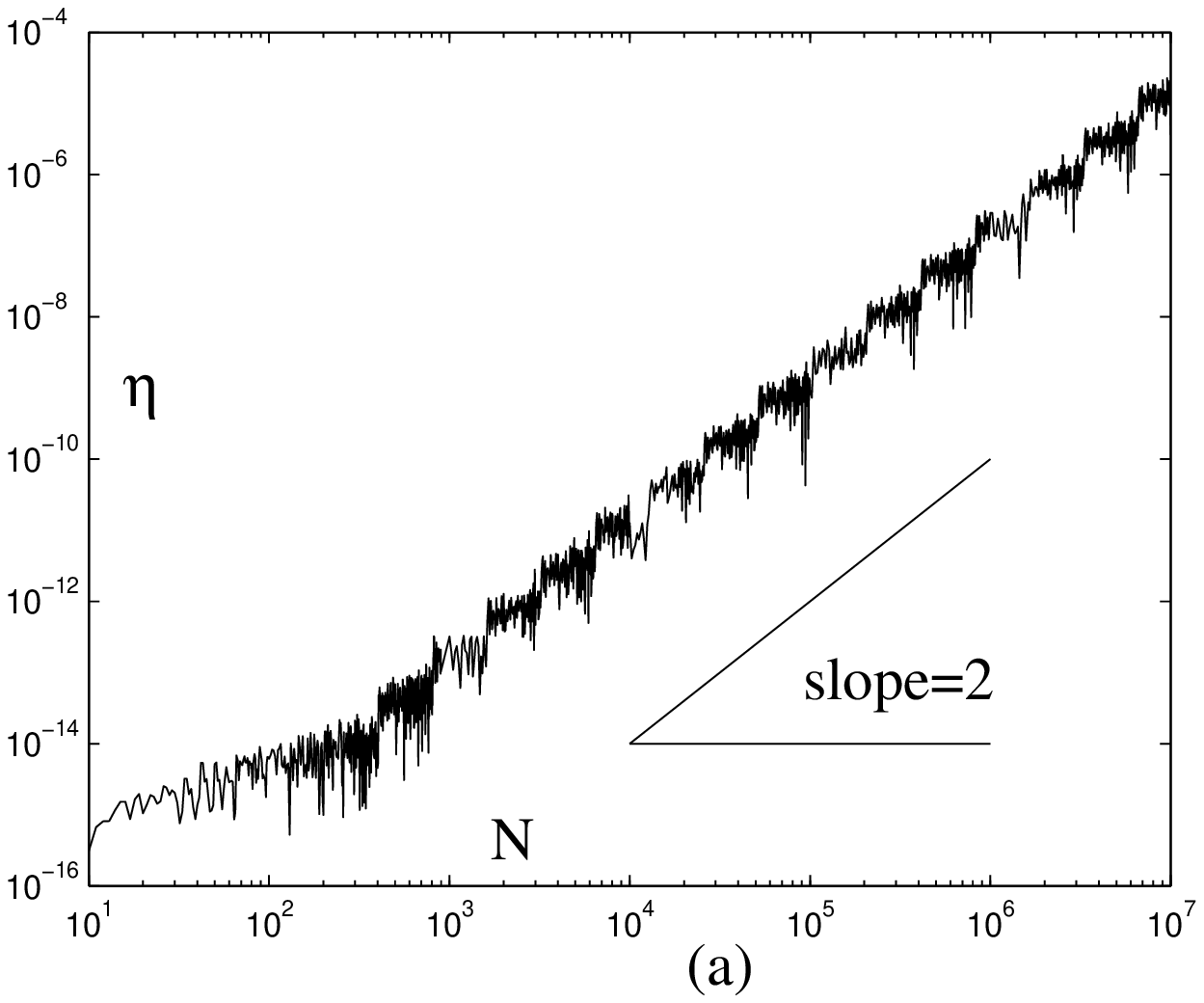}
\includegraphics[width=2.2in]{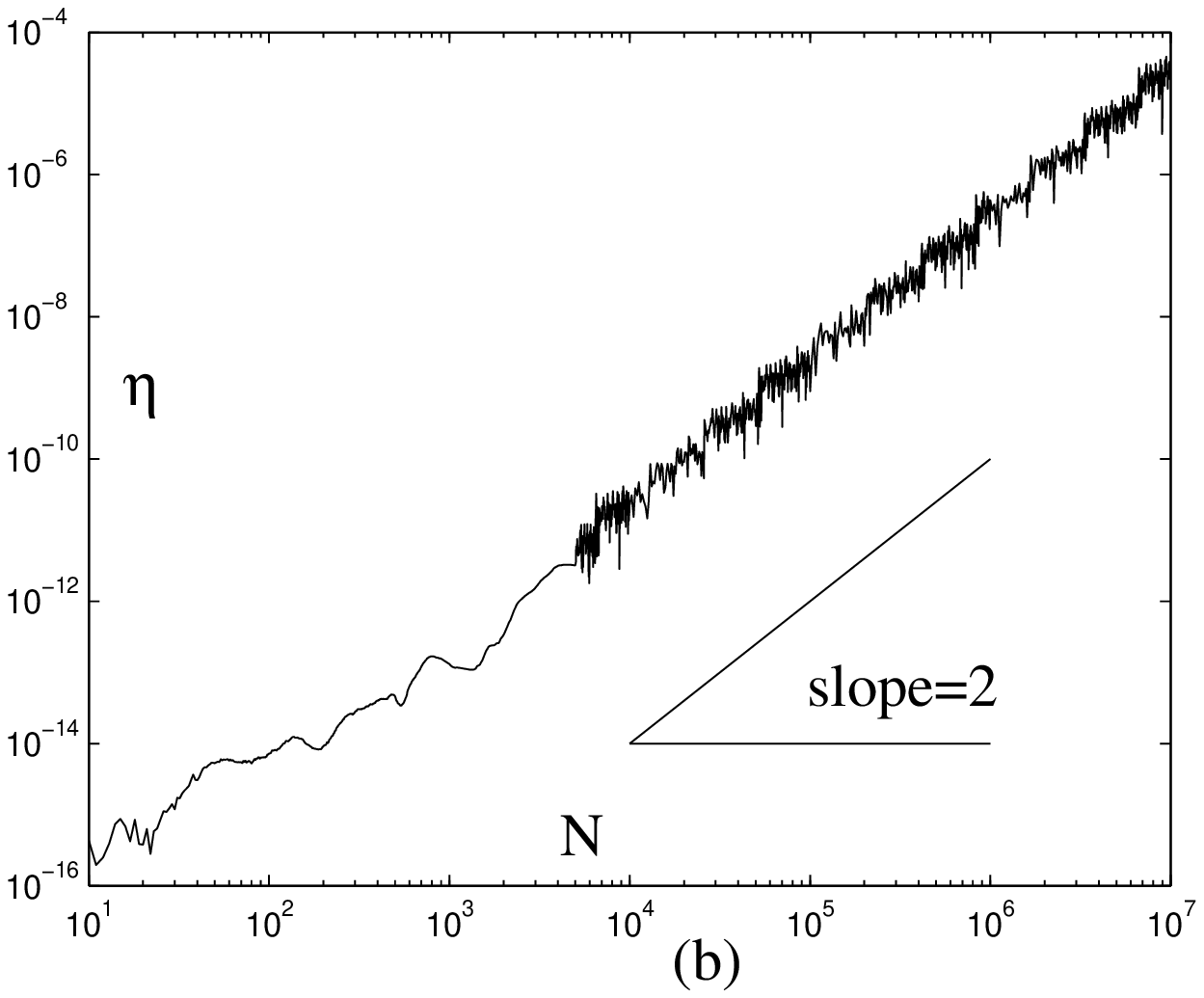}\\
\includegraphics[width=3in]{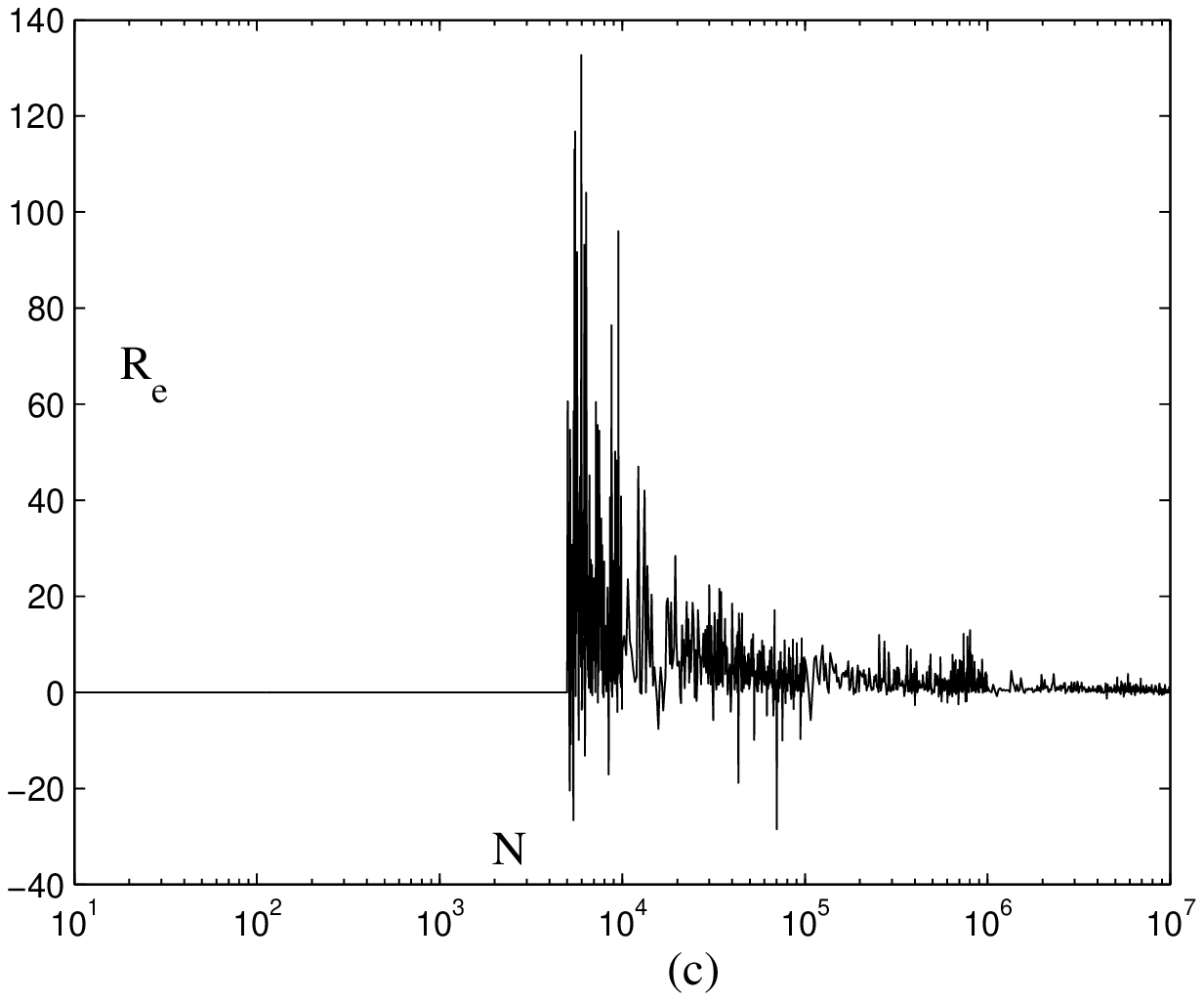}
\vspace{-.3cm}
\caption{\label{Karl1}  (a)  Log-log plot of $\eta^{(a)}=\|x-\hat{x}^{(a)}\|_2/\|x\|_2$ with respect to $N$, where $\hat{x}^{(a)}$ is the computed solution of Problem 2a (using \emph{superLU}). (b) Log-log plot of $\eta^{(b)}=\|x-\hat{x}^{(b)}\|_2/\|x\|_2$ with respect to $N$, where $\hat{x}^{(b)}$ is the computed solution of Problem 2b (using \emph{MUMPS}). (c) Semi-log plot of $100*(\eta^{(c)}-\eta^{(b)})/\eta^{(b)}$ with respect to $N$, where $\eta^{(c)}=\|x-\hat{x}^{(c)}\|_2/\|x\|_2$ and $\hat{x}^{(c)}$ is the computed solution of Problem 1c (using \emph{MUMPS} with 128 processors). Proportionally distributed 1931 values of $N$ in the interval $[10,10^7]$ are used in all the figures.
}
\end{center}
\end{figure}

In \textbf{Problem 2}, we solve the linear system $Ax=b$ using two different software \emph{superLU} and \emph{MUMPS}; we also implement \emph{MUMPS} on multiple processors. Let $\hat{x}^{(a)},\, \hat{x}^{(b)},\, \hat{x}^{(c)}$ be the computed solutions of Problems 2a, 2b, and 2c, respectively. These solutions were computed for $10 \le N \le 10^7$, with $90$ values of $N$ in the range $[10,10^2)$, with 400 values of $N$ in the range $[10^2,10^3)$, and 360 values of $N$ in the range $[10^i,10^{i+1})$, $i=3,4,5,6$, and with $N=10^7$. We presented the log-log plots of  $\eta^{(k)}=\|x-\hat{x}^{(k)}\|_2/\|x\|_2$, $k=a,b$, for the values of $N$ given before, in Figures \ref{Karl1}a and \ref{Karl1}b respectively. We observed that for $N\ge 100$, $\bar{C}_2 / \bar{C}_1 \le 200 < \tau_1$ for both the problems. Thus we do not reject the Hypothesis H for $N\ge 100$. Note that for $N\le 100$, Figures \ref{Karl1}a and b suggest that $\beta \approx -1$. It is also clear from Figure \ref{Karl1}b that the implementation of the algorithm in \emph{MUMPS} changes drastically for $N>5\times 10^3$; this is not the case with \emph{superLU}, as seen in Figure \ref{Karl1}a. Thus the computed solution depends on the software package, as mentioned before. For Problem 2c, we did not display the log-log plot of $\eta^{(c)}=\|x-\hat{x}^{(c)}\|_2/\|x\|_2$ as it would be very similar to the plot of $\eta^{(b)}$ in Figure \ref{Karl1}b. However, we computed $R_e\equiv 100(\eta^{(c)}-\eta^{(b)})/\eta^{(b)}$ -- the ``signed relative difference percent'' --- and presented the semi-log plot of $R_e$ in Figure \ref{Karl1}c for the same values of $N$, given before. For $N \le 5\times 10^3$, we see that $R_e \approx 0$ and values of $R_e$ starts to oscillate for $N > 5\times 10^3$. This indicates that the implementation in MUMPS changes drastically. Figure \ref{Karl1}c also suggests that $\eta^{(c)}$ is larger than $\eta^{(b)}$ for most values on $N$, and $\eta^{(c)}$ gets closer to $\eta^{(b)}$ as $N$ increases.

\begin{figure}
\begin{center}
\includegraphics[width=2.2in]{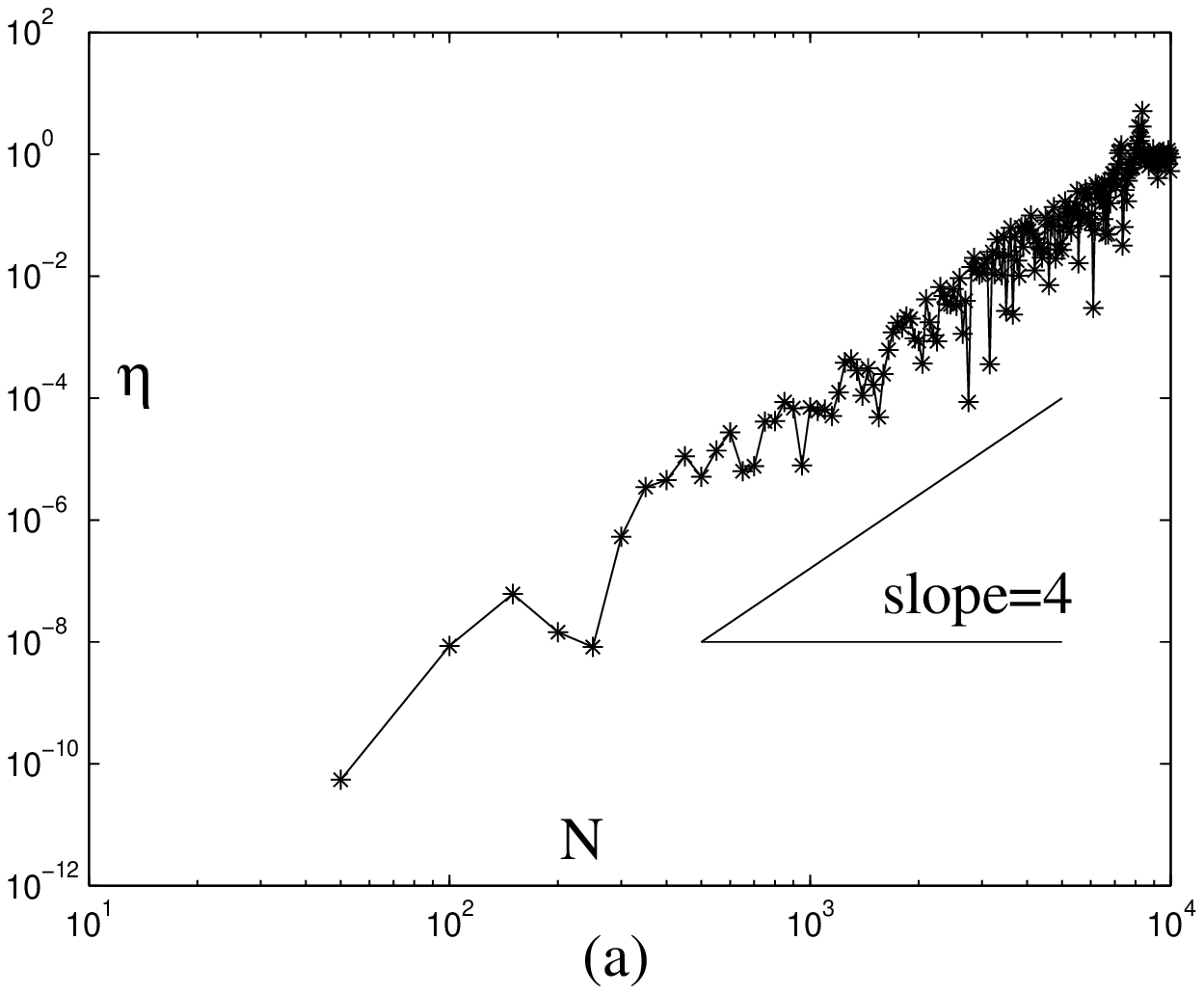}
\includegraphics[width=2.2in]{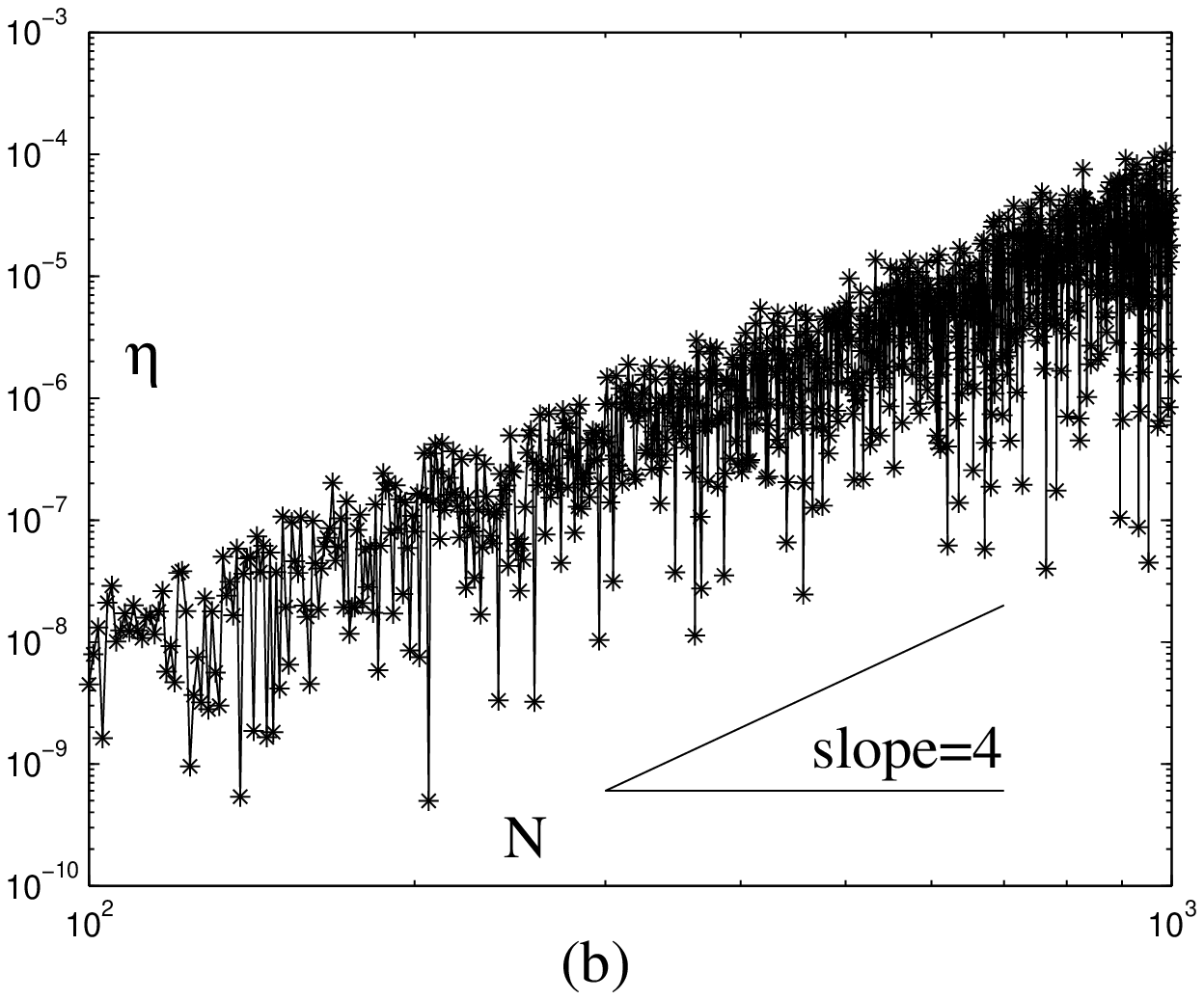}
\vspace{-.4cm}\caption{\label{GFEMplot}  Plots of $\eta = \|x-\hat{x}\|_2/\|x\|_2$ where $\hat{x}$ is the computed solution of the linear system $Ax=b$ of Problem 3, associated with the GFEM with vertices $x_i=ih$, $i=0,1,\cdots,N$. In (a), we used $N=50, 100, 150, \cdots,10000$, and in (b), we used every value of $N$ in the interval $[100,1000]$ to show the detail.    }
\end{center}
\end{figure}

Let $\hat{x}$ be the solution of the linear system $Ax=b$ of \textbf{Problem 3}. We have $\mathfrak{K}(A) = O(N^4)$ (see Section 3.1). The log-log plot of $\eta=\|x-\hat{x}\|_2/\|x\|_2$ with respect to $N=50, 100, 150, \cdots, 10000$ have been presented in Figure \ref{GFEMplot}a. In Figure \ref{GFEMplot}b, we show the details in the range $100 \le N \le 1000$, where we have presented the log-log plot of $\eta$ for every value of $N$ in this range. Based on both these data (i.e., the values of $\eta$ for every value of $N$ in the range $100 \le N \le 1000$ and for $N=1050, 1100, 1150, \cdots, 10000$),we have observed that $\bar{C}_1 N^4 \le \eta \le \bar{C}_2 N^4$ with $\bar{C}_2 / \bar{C}_1 \le 340 < \tau_1$ (note $\tau_2 =0$). Thus we do not reject the Hypothesis H, again based on the subset of meshes with the values of $N$ mentioned above. It is important to note that in Problem 3a, all the digits of accuracy were lost for $N\ge 9000$, and thus the Hypothesis H does not address the value of $N \ge 9000$. We also computed $\eta$ for every value of $N$ in the range $9000 \le N \le 11000$; $\eta$ was of the order $1$ and oscillated around $1$.

\begin{figure}
\begin{center}
\includegraphics[width=2.2in]{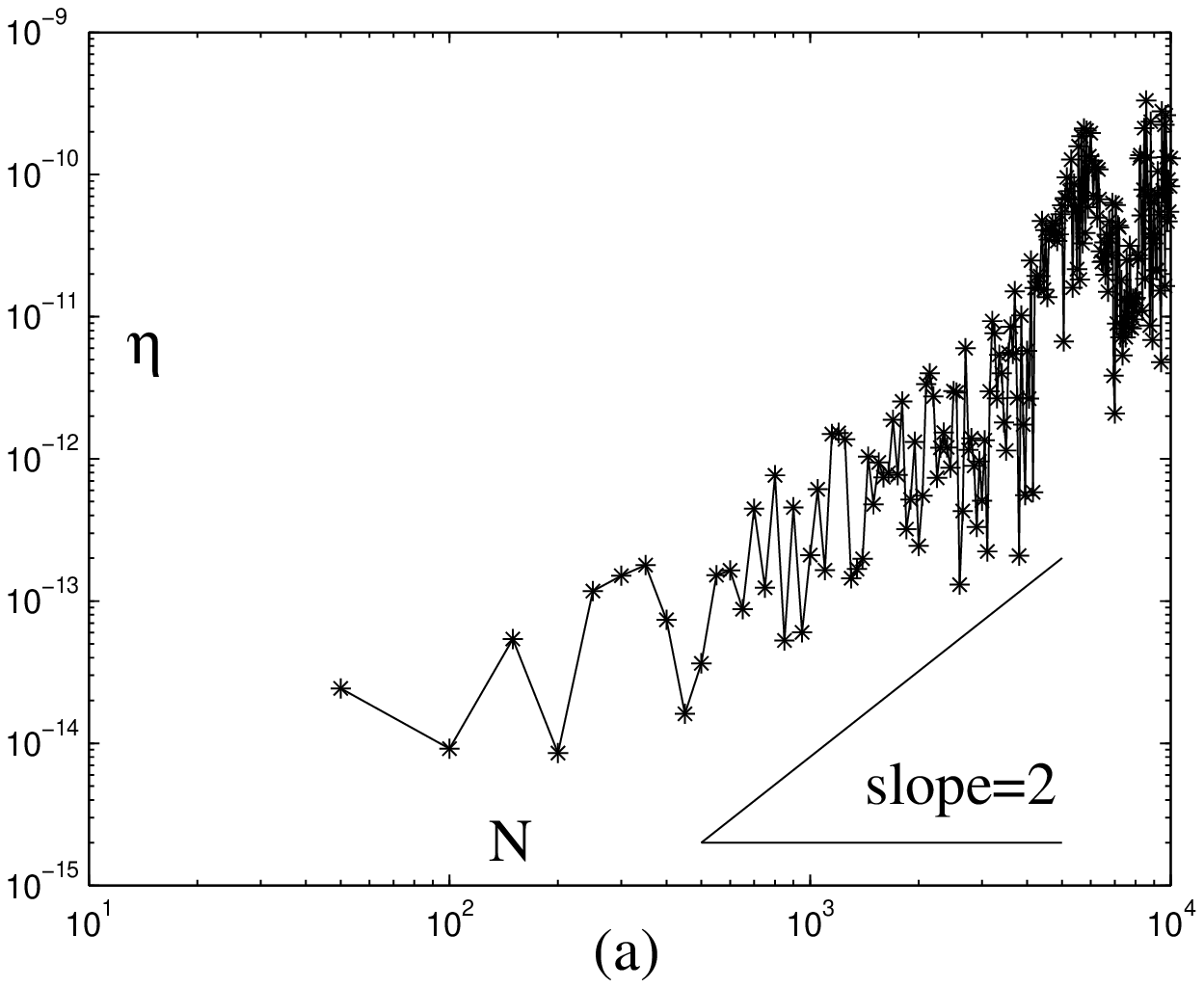}
\includegraphics[width=2.2in]{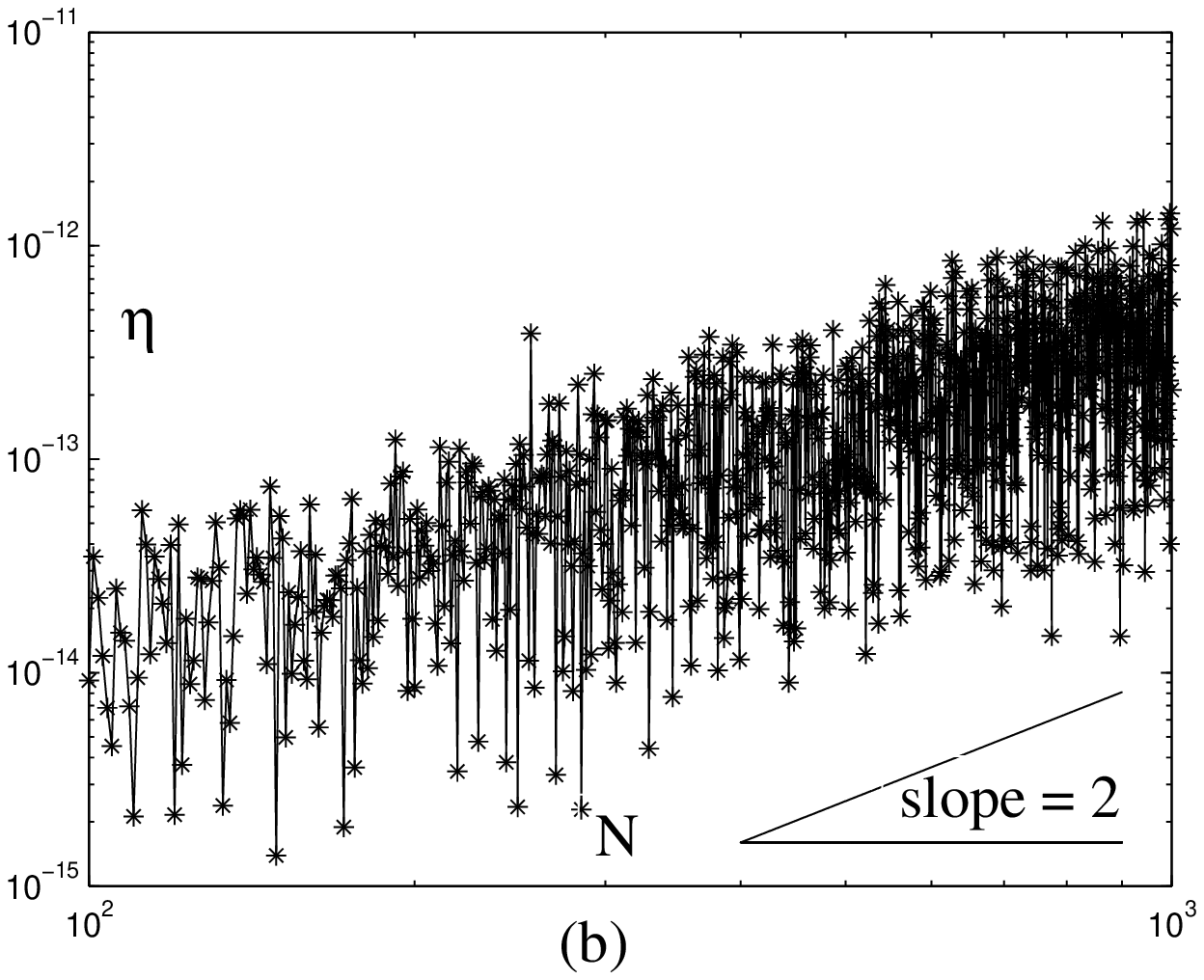}
\vspace{-0.4cm}\caption{\label{SGFEMplot}  Plots of $\eta = \|x-\hat{x}\|_2/\|x\|_2$ where $\hat{x}$ is the computed solution of the linear system $Ax=b$ of Problem 4, associated with the SGFEM with vertices $x_i=ih$, $i=0,1,\cdots,N$. In (a), we used $N=50, 100, 150, \cdots,10000$, and in (b), we used every value of $N$ in the interval $[100,1000]$ to show the detail. }
\end{center}
\end{figure}

Let $\hat{x}$ be the computed solution of the linear system $Ax=b$ of \textbf{Problem 4}. We have shown in this paper that $\mathfrak{K}(A) = O(h^2)$. We have presented the log-log plot of $\eta=\|x-\hat{x}\|_2/\|x\|_2$, with respect to $N=50, 100, 150, \cdots, 10000$ in Figure \ref{SGFEMplot}a, and for every value of $N$ in the range $100 \le N \le 1000$ in Figure \ref{SGFEMplot}b. Based on both these data (i.e., the values of $\eta$ for every value of $N$ in the range $100 \le N \le 1000$ and $N=1050, 1100, 1150, \cdots, 10000$), we observed that $\bar{C}_1 N^2 \le \eta \le \bar{C}_2 N^2$ with $\bar{C}_2/\bar{C}_1 \le 240 < \tau_1$ (note $\tau_2 = 0$). Thus we do not reject the Hypothesis H (based on meshes with these values of $N$).

Thus we did not reject the Hypothesis H for any validation problems with respect to the tolerance $\tau_1 = 400$ and $\tau_2=0$. But we would reject the Hypothesis H if we choose $\tau_1=300$, since $\bar{C}_2 / \bar{C}_1 \le 340 \not \le \tau_1$ in Problem 3. However, if the values of $\eta$ for every value of $N$ in the range $[100,1000]$ were not available (see Figure \ref{GFEMplot}b), then we will have $\bar{C}_2 / \bar{C}_1 \le 250 < \tau_1$, and we thus we would not reject Hypothesis H. Hence validation depends on the values of $N$, i.e., on the number of validation problems considered, since each value of $N$ (in each of Problems 1, 2, 3, and 4) constitutes a separate validation problem. But as mentioned before, the choice of the tolerance depends on the type of decision related to the goal. For example in this paper, we have to decide whether to accept SGFEM over the standard GFEM. In this case, we may allow $\tau_1$ to be bigger; in fact if $\tau_1 = 500$, we still accept SGFEM over GFEM since the value of $\eta$ for GFEM will be much larger than the $\eta$ of SGFEM for large $N$.

We summarize by stating that
\begin{quote} (a) we have confidence in Hypothesis H, based on the chosen validation problems (Problems 1--4). We underline that we have also considered other 2- and 3-dimensional validation problems for the Hypothesis H, which we do not present in this paper. We will present a more substantial validation of Hypothesis H in a future publication.

(b) Because of our confidence in Hypothesis H, we prefer the use of SGFEM over GFEM, since linear system of SGFEM is less prone to the loss of accuracy than the linear system of the GFEM, when solved using an elimination method.
\end{quote}
\begin{myremark}\upshape
As mentioned before, all the computations presented here were performed with $10^{-16}$ accuracy. However, all the figures, presented above, indicate that that the apparent accuracy is about $10^{-18}$. This is likely the effect of various cancelations. \myend
\end{myremark}